\documentclass[a4paper]{amsart}

\usepackage[utf8]{inputenc}			
\usepackage[T1]{fontenc}			
\usepackage[british]{babel} 		

\usepackage{
amsmath,			
amsthm,				
amssymb,			
mathtools,			
bbm,				
stmaryrd,			
calligra,			
enumitem,			
graphicx,			
subcaption, 		
float,				
color,				
tikz,				
microtype,			
hyphenat,			
hyperref,			
cleveref,			
}

\usetikzlibrary{
arrows.meta,
}

\newtheorem{teo}{Theorem}[section]
\newtheorem{lem}[teo]{Lemma}
\newtheorem{coro}[teo]{Corollary}
\newtheorem{fact}[teo]{Fact}
\newtheorem{quest}[teo]{Question}
\newtheorem{conj}[teo]{Conjecture}
\newtheorem{claim}{Claim}[teo]

\theoremstyle{definition}
\newtheorem{defi}[teo]{Definition}
\newtheorem{ejem}[teo]{Example}
\newtheorem{obs}[teo]{Remark}
\newtheorem*{nota}{Notation}

\crefname{teo}{Theorem}{Theorems}
\crefname{lem}{Lemma}{Lemmas}
\crefname{coro}{Corollary}{Corollaries}
\crefname{fact}{Fact}{Facts}
\crefname{claim}{Claim}{Claims}
\crefname{defi}{Definition}{Definitions}
\crefname{obs}{Remark}{Remarks}
\crefname{enumi}{}{}
\crefname{section}{Section}{Sections}

\setenumerate{itemsep=0pt, topsep=3pt}

\renewcommand\qedsymbol{\setbox0=\hbox{\ \ \ \footnotesize{\normalfont Q.E.D.}}\kern\wd0 \strut \hfill \kern-\wd0 \box0}
\newcommand*\claimqed{\renewcommand\qedsymbol{\setbox0=\hbox{\ \ \ {\normalfont \ensuremath{\square}}}\kern\wd0 \strut \hfill \kern-\wd0 \box0}\qedhere\renewcommand\qedsymbol{\setbox0=\hbox{\ \ \ \footnotesize{\normalfont Q.E.D.}}\kern\wd0 \strut \hfill \kern-\wd0 \box0}}

\newcommand*\N{\mathbb{N}}
\newcommand*\Z{\mathbb{Z}}
\newcommand*\Tt{\mathbb{T}}

\newcommand*\tbullet{\text{\raisebox{0.25pt}{\scalebox{0.6}{$\bullet$}}}}
\newcommand*\tq{\mathrel{:}}

\newcommand*\tp{\mathrm{tp}}
\newcommand*\qftp{\mathrm{qftp}}
\newcommand*\id{\mathrm{id}}
\newcommand*\acl{\mathrm{acl}}
\newcommand*\lang{\mathtt{L}}
\newcommand*\Th{\mathrm{Th}}
\newcommand*\Age{\mathrm{Age}}
\newcommand*\Aut{\mathrm{Aut}}

\newcommand*\llangle{\mathopen{\mathord{\langle}\!\mathord{\langle}}}
\newcommand*\rrangle{\mathclose{\mathord{\rangle}\!\mathord{\rangle}}}
\newcommand*\meet{\mathbin{\text{\raisebox{0.25pt}{\scalebox{0.6}{$\wedge$}}}}}
\newcommand*\cone{\mathrm{cone}}
\newcommand*\Rr{\mathrm{R}}

\catcode`@=11 
\def\indep@#1#2{{
\setbox0=\hbox{$#1 x$}
\setbox2=\hbox{$#1 \mathchoice{\scriptstyle\in@xp}{\scriptstyle\in@xp}{\scriptscriptstyle\in@xp}{\scriptscriptstyle\in@xp}$}
\setbox4=\hbox{$#1 \mathchoice{\scriptstyle\inb@se}{\scriptstyle\inb@se}{\scriptscriptstyle\inb@se}{\scriptscriptstyle\inb@se}$}
#1\mathchoice
{\setbox6=\hbox{$\mathop{{}}\limits^{\hbox{\kern\wd0}}_{\hbox{\kern-0.5\wd4\copy4\hss}}$}
\setbox8=\hbox{$\mathop{\hbox{\kern0.9pt\copy2\hss}}\limits^{\hbox{\kern\wd0}}_{\hbox{\kern-0.5\wd4\copy4\hss}}$}
\hbox{\kern\wd6\hbox to 0pt{\hss$#1\mid$\hss}\lower.9\ht0\hbox to 0pt{\hss$#1\smile$\hss}\raise1\ht0\hbox to 0pt{\kern0.9pt\copy2\hss}\lower1\ht0\hbox to 0pt{\hss\lower\ht4\hbox{\copy4}\hss}\hbox to 0pt{}\kern\wd8}}
{\setbox6=\hbox{$\mathop{\hbox{\kern0.9pt\copy2\hss}}\limits^{\hbox{\kern\wd0}}_{\hbox{\kern0.6\wd0 \copy4\hss}}$}
\hbox{\kern\wd0\hbox to 0pt{\hss$#1\mid$\hss}\lower.9\ht0\hbox to 0pt{\hss$#1\smile$\hss}\raise1\ht0\hbox to 0pt{\kern0.9pt\copy2\hss}\kern0.6\wd0\lower0.2\ht0\hbox{\hss\lower\ht4\copy4\hss}\kern-0.6\wd0\kern-\wd4\hbox to 0pt{}\kern\wd6}}
{\setbox6=\hbox{$\mathop{\hbox{\kern0.9pt\copy2\hss}}\limits^{\hbox{\kern\wd0}}_{\hbox{\kern0.6\wd0 \copy4\hss}}$}
\hbox{\kern\wd0\hbox to 0pt{\hss$#1\mid$\hss}\lower.9\ht0\hbox to 0pt{\hss$#1\smile$\hss}\raise1\ht0\hbox to 0pt{\kern0.9pt\copy2\hss}\kern0.6\wd0\lower0.2\ht0\hbox{\hss\lower\ht4\copy4\hss}\kern-0.6\wd0\kern-\wd4\hbox to 0pt{}\kern\wd6}}
{\setbox6=\hbox{$\mathop{\hbox{\kern0.9pt\copy2\hss}}\limits^{\hbox{\kern\wd0}}_{\hbox{\kern0.6\wd0 \copy4\hss}}$}
\hbox{\kern\wd0\hbox to 0pt{\hss$#1\mid$\hss}\lower.9\ht0\hbox to 0pt{\hss$#1\smile$\hss}\raise1\ht0\hbox to 0pt{\kern0.9pt\copy2\hss}\kern0.6\wd0\lower0.2\ht0\hbox{\hss\lower\ht4\copy4\hss}\kern-0.6\wd0\kern-\wd4\hbox to 0pt{}\kern\wd6}} 
}}
\def\indep#1#2{\def\in@xp{#1}\def\inb@se{#2}\mathop{\mathpalette\indep@\relax}}
\catcode`@=12
\newcommand*\indepe[2][{}]{\indep{#1}{#2}}

\catcode`@=11 
\def\notindep@#1#2{{
\setbox0=\hbox{$#1 x$}
\setbox2=\hbox{$#1 \mathchoice{\scriptstyle\notin@xp}{\scriptstyle\notin@xp}{\scriptscriptstyle\notin@xp}{\scriptscriptstyle\notin@xp}$}
\setbox4=\hbox{$#1 \mathchoice{\scriptstyle\notinb@se}{\scriptstyle\notinb@se}{\scriptscriptstyle\notinb@se}{\scriptscriptstyle\notinb@se}$}
#1\mathchoice
{\setbox6=\hbox{$\mathop{{}}\limits^{\hbox{\kern\wd0}}_{\hbox{\kern-0.5\wd4\copy4\hss}}$}
\setbox8=\hbox{$\mathop{\hbox{\kern0.9pt\copy2\hss}}\limits^{\hbox{\kern\wd0}}_{\hbox{\kern-0.5\wd4\copy4\hss}}$}
\hbox{\kern\wd6\hbox to 0pt{\mathchardef\nn=12854\hss$#1 \nn$\kern1.4\wd0\hss}\hbox to 0pt{\hss$#1\mid$\hss}\lower.9\ht0\hbox to 0pt{\hss$#1\smile$\hss}\raise1\ht0\hbox to 0pt{\kern0.9pt\copy2\hss}\lower1\ht0\hbox to 0pt{\hss\lower\ht4\hbox{\copy4}\hss}\hbox to 0pt{}\kern\wd8}}
{\setbox6=\hbox{$\mathop{\hbox{\kern0.9pt\copy2\hss}}\limits_{\hbox{\kern0.6\wd0 \copy4\hss}}$}
\hbox{\kern\wd0\hbox to 0pt{\mathchardef\nn=12854\hss$#1 \nn$\kern1.4\wd0\hss}\hbox to 0pt{\hss$#1\mid$\hss}\lower.9\ht0\hbox to 0pt{\hss$#1\smile$\hss}\raise1\ht0\hbox to 0pt{\kern0.9pt\copy2\hss}\kern0.6\wd0\lower0.2\ht0\hbox{\hss\lower\ht4\copy4\hss}\kern-0.6\wd0\kern-\wd4\hbox to 0pt{}\kern\wd6}}
{\setbox6=\hbox{$\mathop{\hbox{\kern0.9pt\copy2\hss}}\limits_{\hbox{\kern0.6\wd0 \copy4\hss}}$}
\hbox{\kern\wd0\hbox to 0pt{\mathchardef\nn=12854\hss$#1 \nn$\kern1.4\wd0\hss}\hbox to 0pt{\hss$#1\mid$\hss}\lower.9\ht0\hbox to 0pt{\hss$#1\smile$\hss}\raise1\ht0\hbox to 0pt{\kern0.9pt\copy2\hss}\kern0.6\wd0\lower0.2\ht0\hbox{\hss\lower\ht4\copy4\hss}\kern-0.6\wd0\kern-\wd4\hbox to 0pt{}\kern\wd6}}
{\setbox6=\hbox{$\mathop{\hbox{\kern0.9pt\copy2\hss}}\limits_{\hbox{\kern0.6\wd0 \copy4\hss}}$}
\hbox{\kern\wd0\hbox to 0pt{\mathchardef\nn=12854\hss$#1 \nn$\kern1.4\wd0\hss}\hbox to 0pt{\hss$#1\mid$\hss}\lower.9\ht0\hbox to 0pt{\hss$#1\smile$\hss}\raise1\ht0\hbox to 0pt{\kern0.9pt\copy2\hss}\kern0.6\wd0\lower0.2\ht0\hbox{\hss\lower\ht4\copy4\hss}\kern-0.6\wd0\kern-\wd4\hbox to 0pt{}\kern\wd6}} 
}}
\def\notindep#1#2{\def\notin@xp{#1}\def\notinb@se{#2}\mathop{\mathpalette\notindep@\relax}}
\catcode`@=12

\title{Automorphisms of the Rado meet{\hyp}tree}
\author{Itay Kaplan}
\address{Einstein Institute of Mathematics, Hebrew University of
Jerusalem, 91904, Jerusalem Israel.}
\email{kaplan@math.huji.ac.il}

\author{Binyamin Riahi}
\address{Einstein Institute of Mathematics, Hebrew University of
Jerusalem, 91904, Jerusalem Israel.}
\email{binyamin.riahi@mail.huji.ac.il}
\author{Arturo Rodr\'{\i}guez Fanlo}
\address{Einstein Institute of Mathematics, Hebrew University of
Jerusalem, 91904, Jerusalem Israel.}
\email{Arturo.Rodriguez@mail.huji.ac.il}

\thanks{Rodriguez Fanlo supported by the Israel Academy of Sciences and Humanities \& Council for Higher Education Excellence Fellowship Program for International Postdoctoral Researchers.}
\thanks{Kaplan would like to thank the Israel Science Foundation (ISF) for their
support of this research (grants no. 1254/18 and 804/22).}

\begin{document}
\maketitle
\begin{abstract}
We prove that the group of automorphisms of the generic meet{\hyp}tree expansion of an infinite non{\hyp}unary free Fra\"{\i}ss\'{e} limit over a finite relational language is simple. As a prototypical case, the group of automorphism of the Rado meet{\hyp}tree (i.e. the Fra\"{\i}ss\'{e} limit of finite graphs which are also meet{\hyp}trees) is simple.
\end{abstract}

\section*{Introduction}

There is a long history studying the normal subgroups of groups of automorphisms of countable homogeneous structures. Starting with the countable set \cite{onofri1929teoria,schreier1933uber} and the countable dense linear order without end points \cite{holland1963lattice,lloyd1964lattice}, and continuing with, for example, the Rado graph \cite{truss1989group}, the countable $2${\hyp}homogeneous trees \cite{droste1989automorphism}, the countable universal homogeneous partially ordered set \cite{glass1993automorphism}, the countable universal homogeneous $k${\hyp}uniform hypergraph \cite{lovell2007automorphism}, and $\Aut(\mathbb{C}/\mathbb{Q}^{\mathrm{alg}})$ \cite{lascar1992automorphismes,lascar1997group}. 

A general method for this line of results was provided by Macpherson and Tent \cite{macpherson2011simplicity}, consequently improved by Tent and Ziegler \cite{tent2013isometry}. The strategy in \cite{tent2013isometry} has two components: the existence of a stationary independence relation (i.e. a ternary relation on finite subsets satisfying various properties; see \cite[Definition 2.1]{tent2013isometry}) and automorphisms which move maximally \cite[Definition 2.5]{tent2013isometry}. This technique was used in \cite{macpherson2011simplicity} to show that the group of automorphisms of a transitive non{\hyp}trivial free Fra\"{\i}ss\'{e} limit is simple. In \cite{tent2013isometry}, it shows that the group of isometries of the Urysohn space modulo the normal subgroup of bounded isometries is simple. 

An involved variation of this technique was used by Calderoni, Kwiatkowska and Tent in \cite{calderoni2021simplicity} to show that generic linear order expansions (i.e. expansions given by adding a generic dense linear order) of countable transitive non{\hyp}trivial free Fra\"{\i}ss\'{e} limits on relational languages have simple groups of automorphisms. Independently, Li further generalised the Macpherson{\hyp}Tent{\hyp}Ziegler technique in \cite{li2019automorphism}. In Li's approach, stationary independence relations were replaced by stationary weak independence relations  and moving maximally was replaced by ``almost''\footnote{In this paper, we omit the word ``almost'' as we only use this notion of moving maximally; see \cref{d:move maximally}.}  moving (right/left) maximally. She used this improved version to give an alternative proof of Calderoni{\hyp}Kwiatkowska{\hyp}Tent result and more in \cite{li2020automorphism}.

An ordered tree is a semilinear partially ordered set. A meet{\hyp}tree is an ordered tree which is also a lower semilattice. Trees and meet{\hyp}trees are very useful in model theory. In particular, they depict some interesting pathologies of NIP theories, see \cite{05-KaSh975,06-KaSh946}.

While most of the studied countable homomogeneous structures have very few normal subgroups of the group of automorphisms, homogeneous trees have the maximal number (i.e. $2^{2^{\aleph_0}}$) of normal subgroups \cite[Theorem 1.4]{macpherson2011simplicity}. In particular, this applies to the universal dense meet{\hyp}tree, i.e. the Fra\"{\i}ss\'{e} limit of finite meet{\hyp}trees. So, it is natural to wonder whether the Calderoni{\hyp}Kwiatkowska{\hyp}Tent's and Li's result mentioned above can be adapted to the case of the universal dense meet{\hyp}tree.

In this paper, we positively answer this question. Namely, we show that the groups of automorphisms of generic meet{\hyp}tree expansions of countable non{\hyp}unary free Fra\"{\i}ss\'{e} limits over relational languages are simple (\cref{c:main}). As a prototypical example, think of the \emph{Rado meet{\hyp}tree} $\Tt^{\Rr}$, i.e. the Fra\"{\i}ss\'{e} limit of finite graphs which are also meet{\hyp}trees. A direct consequence of our result is that $\Aut(\Tt^{\Rr})$ is simple. 

Our approach is very similar to Li's. However, in our case, we follow a divide{\hyp}and{\hyp}conquer strategy. We divide our group of automorphisms to those that setwise fix a branch and those that fix a point. While the universal dense meet{\hyp}tree does not admit a stationary weak independence relation \cite[Corollary 6.10]{kaplan2020automorphism}, these two expansions (the one given by adding a predicate for the branch and the one given by adding a constant) do admit stationary weak independence relations, allowing us to apply the results of \cite{li2019automorphism}.

This paper is organised as follows. In \cref{s:section 1} we recall some basic notions about meet{\hyp}trees. In \cref{s:section 2}, we introduce the above mentioned stationary weak independence relations (\cref{d:cone-independence,d:semibranch-independence,l:cone-independence,c:semibranch-independence}). In \cref{s:section 3}, we define the notion of generic meet{\hyp}tree expansion in \cref{d:generic meet-tree expansion} and prove its existence for strong Fra\"{\i}ss\'{e} limits over relational languages in \cref{l:generic meet-tree expansions}. We start to manipulate automorphisms in \cref{s:section 4}. This section consists of long back{\hyp}and{\hyp}forth arguments and it is the technical core of this paper. It allows us to transform non{\hyp}trivial automorphisms to have strong non{\hyp}structure properties, reducing the general case to two main cases (\cref{l:dicotomia,c:fan strictly increasing,c:semibranch strictly increasing}). In \cref{s:section 5}, we show that, after a commutation, these two main cases move maximally (\cref{l:alpha ordering,c:branched alpha ordering,l:move maximally}), getting our main results, \cref{t:main,c:main}, by applying \cite[Corollary 2.16]{li2019automorphism}. 

\subsection*{Acknowledgements}
We would like to thank Dugald Macpherson for telling us about \cite{droste1989automorphism}.

\begin{nota} Throughout the paper we deal with several languages $\lang_\star$. In each case, we write $\langle\tbullet\rangle_\star$, $\tp_\star$, $\qftp_\star$ and $\equiv^\star$ when working in $\lang_\star$. \end{nota}

\section{Preliminaries on meet{\hyp}trees} \label{s:section 1}
Recall the definition of ordered trees and meet{\hyp}trees from the introduction. An \emph{ordered tree} is a semilinear partially ordered set, i.e. a partially ordered set $(T,\leq)$ such that, for any $a\in T$, $T_{<a}$ with the induced order is linearly ordered. A \emph{meet{\hyp}tree} is an ordered tree which is also a lower semilattice, i.e. any two elements $a,b\in T$ have a largest common smaller or equal element. The \emph{meet operation} on $T$ is the binary operation defined as $\meet:\ a,b\mapsto a\meet b\coloneqq \max\{d\in T\tq d\leq a\mathrm{\ and\ }d\leq b\}$. 

A \emph{branch} is a maximal chain. A \emph{(closed) semibranch} $\Gamma$ of a meet{\hyp}tree $T$ is a downwards closed chain such that, for any $a\in T$, there is a largest element in $\Gamma$ smaller or equal than $a$. Given a semibranch $\Gamma$ of a meet{\hyp}tree $T$, the \emph{semibranch{\hyp}projection} is the function $\pi_\Gamma:\ T\rightarrow \Gamma$ given by $a\mapsto \pi_\Gamma(a)\coloneqq \max\{b\in\Gamma\tq b\leq a\}$. A \emph{semibranched meet{\hyp}tree} is a meet{\hyp}tree with a distinguished semibranch and its semibranch{\hyp}projection. Given a meet{\hyp}tree $T$ and a semibranch $\Gamma$, we write $T_\Gamma$ for the corresponding semibranched meet{\hyp}tree. For a semibranch $\Gamma$, we write $\Gamma^\circ\coloneqq \{a\in \Gamma\tq \mathrm{there\ is\ }b\in \Gamma\mathrm{\ with\ }a<b\}$.

A \emph{pointed meet{\hyp}tree} is a meet{\hyp}tree with a distinguished point. Given a meet{\hyp}tree $T$ and a point $\gamma$, we write $T_\gamma$ for the corresponding pointed meet{\hyp}tree. The \emph{(open) cone equivalence relation above $\gamma$} is the equivalence relation $\sim_\gamma$ on $T_{> \gamma}$ given by $a\sim_\gamma b$ if and only if $a\meet b>\gamma$. The \emph{(open) cones} above $\gamma$ are the equivalence classes under the cone equivalence relation above $\gamma$. For $a>\gamma$, we write $\cone_\gamma(a)$ for the cone above $\gamma$ containing $a$. For $a\not >\gamma$, set $\cone_\gamma(a)=\emptyset$. For a subset $A$ of $T$, write $\cone_\gamma(A)\coloneqq\bigcup_{a\in A}\cone_\gamma(a)$.\smallskip

The following results are elementary:
\begin{lem} \label{l:semibranches} Let $T$ be a meet{\hyp}tree and $\Gamma\subseteq T$. Then, $\Gamma$ is a semibranch if and only if $\Gamma$ is a branch or $\Gamma=T_{\leq \gamma}$ where $\gamma=\max\Gamma$ and $\pi_\Gamma(a)=a\meet \gamma$ for all $a\in T$.
\end{lem} 
\begin{coro} \label{c:semibranches tree} Let $T$ be a meet{\hyp}tree and $T^{\infty}\coloneqq\{\Gamma\subseteq T\tq \Gamma\mathrm{\ semibranch}\}$. Then, $(T^{\infty},\subseteq,\cap)$ is a meet{\hyp}tree and every element of $T^{\infty}$ is bounded by a maximal element, where the maximal elements of $T^{\infty}$ are the branches of $T$. Furthermore, the canonical map $\gamma\mapsto T_{\leq \gamma}$ is an embedding of meet{\hyp}trees from $T$ into $T^{\infty}$. In the particular case of $T^{\infty}$ in place of $T$, the canonical embedding of $T^{\infty}$ into $(T^{\infty})^{\infty}$ is an isomorphism.
\end{coro} 

The first{\hyp}order language of ordered trees is the first{\hyp}order language of ordered sets $\lang_{<}$, i.e. the one{\hyp}sorted language with equality and a binary relation symbol for the order. The first{\hyp}order language of meet{\hyp}trees $\lang_{\wedge}$ is the expansion of $\lang_<$ given by adding a binary function symbol for the meet operation. The first{\hyp}order language of semibranched meet{\hyp}trees $\lang_{\Gamma}$ is the expansion of $\lang_{\wedge}$ given by adding a unary relation symbol for the semibranch and a unary function symbol for the semibranch{\hyp}projection. The first{\hyp}order language of pointed meet{\hyp}trees is the expansion $\lang_\gamma$ of $\lang_{\wedge}$ given by adding a constant symbol for the point.



The following lemma is elementary:
\begin{lem}\label{l:generated} Let $T$ be a meet{\hyp}tree, $\Gamma$ a semibranch, $\gamma$ a point and $X$ a subset. Then:
\begin{enumerate}[label={\rm{(\arabic*)}}, ref={\rm{\arabic*}}, wide]
\item \label{itm:generated 1} $\langle X\rangle_{\gamma}=\langle X\,\gamma\rangle_{\wedge}$.
\item \label{itm:generated 2} If $X$ is finite and $x\in X$, $\bigwedge X=\min\{x\meet y\tq y\in X\}$. 
\item \label{itm:generated 3} $\langle X\rangle_{\wedge}=\{x\meet y\tq x,y\in X\}$.
\item \label{itm:generated 4} If $X$ is finite non{\hyp}empty and $x\notin X$,  
$\langle X\, x\rangle_{\wedge}=\langle X\rangle_{\wedge}\cup \{e,x\}$,
where $e\coloneqq\max\{x\meet y\tq y\in X\}$. It follows by induction that $|\langle X\rangle_{\wedge}|\leq 2|X|-1$.
\item \label{itm:generated 5} For any $x,y\in X$, we have $\pi_\Gamma(x\meet y)=\pi_\Gamma(x)\meet \pi_\Gamma(y)=\pi_\Gamma(x)\meet y=\min\{\pi_\Gamma(x),\pi_\Gamma(y)\}$.
\item \label{itm:generated 6} $\langle X\rangle_{\Gamma}=\langle X\rangle_{\wedge}\cup \{\pi_{\Gamma}(x)\tq x\in X\}$ and $\langle X\rangle_\Gamma\cap \Gamma=\{\pi_\Gamma(x)\tq x\in X\}$.
\item \label{itm:generated 7} If $X$ is finite, 
$\langle X\rangle_{\Gamma}=\langle X\,\gamma_X\rangle_{\wedge}$ where $\gamma_X=\max\{\pi_\Gamma(x)\tq x\in X\}$. By {\rm{(\cref{itm:generated 4})}}, $|\langle X\rangle_\Gamma|\leq 2|X|+1$. Moreover, for any $\gamma_0\in\Gamma$, we have 
$\langle X\,\gamma_0\rangle_{\wedge}=\langle X\,\gamma_X\rangle_{\wedge}\cup\{\gamma_0\}$.
\item \label{itm:generated 8} $\cone_\gamma(X)=\cone_\gamma(\langle X\rangle_\gamma)=\cone_\gamma(\langle X\rangle_\Gamma)$.
\end{enumerate}
\end{lem}

\begin{lem} \label{l:quantifier free type} Let $\lang$ be a language and $T$ a theory on $\lang$. Let $\lang_{\mathrm{rel}}$ be the relational reduct, i.e. the reduct consisting on the relation symbols of $\lang$. Assume that for every function symbol of $\lang$, there is a formula in $\lang_{\mathrm{rel}}$ defining it in $T$. Let $a$ and $b$ enumerate two models of $T$. Then, $\qftp_{\lang}(a)=\qftp_{\lang}(b)$ if and only if $\qftp_{\lang_{\mathrm{rel}}}(a)=\qftp_{\lang_{\mathrm{rel}}}(b)$. In other words, every $\lang_{\mathrm{rel}}${\hyp}isomorphism between two models of $T$ is an $\lang${\hyp}isomorphism.
\begin{proof}
Note that, by the hypotheses, any $\lang_{\mathrm{rel}}${\hyp}isomorphism of $T${\hyp}models is an $\lang${\hyp}isomorphism.
\end{proof}
\end{lem}

In particular, \cref{l:quantifier free type} applies to meet{\hyp}trees and semibranched meet{\hyp}trees. We get then the following corollaries:
\begin{coro}\label{c:quantifier free type tree}
Let $T$ be a meet{\hyp}tree, $\Gamma$ a semibranch and $\gamma$ a point. Let $a=(a_i)_{i<n}$ and $b=(b_i)_{i<n}$ be finite tuples in $T$. Then:
\begin{enumerate}[label={\rm{(\arabic*)}}, ref={\rm{\arabic*}}, wide]
\item \label{itm:quantifier free type tree 1} $\qftp_{\gamma}(a)=\qftp_{\gamma}(b)$ if and only if $\qftp_{\wedge}(a,\gamma)=\qftp_{\wedge}(b,\gamma)$.
\item \label{itm:quantifier free type tree 2} $\qftp_{\wedge}(a)=\qftp_{\wedge}(b)$ if and only if, for any $i,j,k<n$, we have that $a_i\meet a_j\leq a_k$ if and only if $b_i\meet b_j\leq b_k$.
\item \label{itm:quantifier free type tree 3} Pick $\gamma_a=\max\{\pi_{\Gamma}(a_i)\}_{i<n}$ and $\gamma_b=\max\{\pi_{\Gamma}(b_i)\}_{i<n}$. Then, $\qftp_{\Gamma}(a)=\qftp_{\Gamma}(b)$ if and only if $\qftp_{\wedge}(a,\gamma_a)=\qftp_{\wedge}(b,\gamma_b)$. Furthermore, for any $\gamma\in \Gamma$ with $\gamma>\gamma_a$ and $\gamma>\gamma_b$, $\qftp_{\Gamma}(a)=\qftp_{\Gamma}(b)$ if and only if $\qftp_\gamma(a)=\qftp_\gamma(b)$.
\end{enumerate}
\end{coro}

\begin{coro} \label{c:cones are indiscernible} Let $T$ be a meet{\hyp}tree, $\gamma$ a point, $C$ a finite subset of $T$ and $a,b$ single elements above $\gamma$. Suppose $a\notin\cone_{\gamma}(C)$ and $b\notin \cone_\gamma(C)$. Then, $\qftp_\gamma(a/C)=\qftp_\gamma(b/C)$.
\begin{proof} By \cref{l:generated}(\ref{itm:generated 4}), we have that $\langle a,C\rangle_\gamma=\langle C\rangle_\gamma\cup \{a\}$ and $\langle b,C\rangle_\gamma=\langle C\rangle_\gamma\cup\{b\}$. Furthermore, for $x\in\langle C\rangle_\gamma$, $x\leq a$ if and only if $x\leq \gamma$, if and only if $x\leq b$. Thus, by \cref{l:quantifier free type}, we conclude that $\qftp_\gamma(a/C)=\qftp_\gamma(b/C)$.
\end{proof}
\end{coro} 

\begin{coro}\label{c:semibranch intervals are indiscernible}
Let $T$ be a meet{\hyp}tree, $\Gamma$ a semibranch, $C$ a finite subset of $T$ and $a,b$ single elements in $\Gamma$. Write $\pi(C)=\{\pi_{\Gamma}(c)\tq c\in C\}$. Then, $\qftp_\Gamma(a/C)=\qftp_\Gamma(b/C)$ if and only if $\qftp_{<}(a/\pi(C))=\qftp_{<}(b/\pi(C))$.
\begin{proof} The ``only if'' implication is obvious, let's prove the ``if'' one. Set $\bar{C}=\langle C\rangle_{\Gamma}$. By \cref{l:generated}(\ref{itm:generated 7}), $\langle a\, C\rangle_\Gamma=\bar{C}\cup \{a\}$ and $\langle b\, C\rangle_\Gamma=\bar{C}\cup \{b\}$. By \cref{l:quantifier free type}, $\qftp_\Gamma(b/C)=\qftp_\Gamma(a/C)$ if and only if $\qftp_{<}(b/\bar{C})=\qftp_{<}(a/\bar{C})$. Pick $c\in\bar{C}$ arbitrary. If $c\leq a$, then $c\in \Gamma$ so, by \cref{l:generated}(\ref{itm:generated 6}), $c\in\pi(C)$ and, by hypothesis, $c\leq b$. If $a\leq c$, as $a\in \Gamma$, then $a\leq \pi_{\Gamma}(c)\leq c$ and, by \cref{l:generated}(\ref{itm:generated 6}), $\pi_{\Gamma}(c)\in \pi(C)$, so $b\leq \pi_{\Gamma}(c)\leq c$ by hypothesis. Symmetrically, $c\leq b$ implies $c\leq a$ and $b\leq c$ implies $a\leq c$.
\end{proof}
\end{coro}

The following are very useful lemmas:
\begin{lem}[Three points]\label{l:three elements} Let $T$ be a meet{\hyp}tree and $a,b,c$ three points. At least two of the three possible pairwise meets between $a,b,c$ equal $a\meet b\meet c$. In other words, if $a\meet b\neq a\meet b\meet c$, then $a\meet c=b\meet c=a\meet b\meet c$.    
\begin{proof} See \cite[Fact 2.19]{kaplan2020automorphism}.
\end{proof}
\end{lem}
\begin{lem}[Four points]\label{l:four elements}
Let $T$ be a meet{\hyp}tree and $a,b,c,d$ four points. If $a\meet c=b\meet c$ and $a\meet d\neq b\meet d$, then $a\meet c=d\meet c$. 
\begin{proof} By \cref{l:three elements} with $a,b,d$, either $a\meet b=a\meet d<b\meet d$ or $a\meet b= b\meet d<a\meet d$. Since the statement is symmetric on $a$ and $b$, we can assume without loss of generality that $b\meet d<a\meet d$. By \cref{l:three elements} with $a,b,c$, we have $a\meet c=b\meet c\leq a\meet b$. Therefore, $a\meet c<a\meet d$, concluding $a\meet c=d\meet c$ by \cref{l:three elements} with $a,c,d$. 
\end{proof}
\end{lem}
\begin{figure}[ht]
\begin{center}
\begin{tikzpicture}[scale=0.8]
\node (xL) at (-6,-3) {};  \node (yL) at (-6,3) {};  
\node (xR) at (6,-3) {}; \node (yR) at (6,3) {};


\node (bot) at (1,-2) {}; 
\node (d) [draw, shape=circle, fill, scale=0.15, label=right:$d$] at (-3,2) {}; 
\node (a) [draw, shape=circle, fill, scale=0.15, label=right:$a$] at (-1,2) {}; 
\node (b) [draw, shape=circle, fill, scale=0.15, label=right:$b$] at (1,2) {}; 
\node (c) [draw, shape=circle, fill, scale=0.15, label=right:$c$] at (3,2) {};
\node (ad) [draw, shape=circle, fill, scale=0.15, label={left:$a\meet d$}] at (-2,1) {};
\node (bd) [draw, shape=circle, fill, scale=0.15, label=left:$b\meet d$] at (-1,0) {}; 
\node (ac) [draw, shape=circle, fill, scale=0.15, label=left:$a\meet c$] at (0,-1) {}; 

\draw (d) -- (bot);
\draw (a) -- (ad); 
\draw (b) -- (bd);
\draw (c) -- (ac); 
\end{tikzpicture}
\caption{Four points Lemma.}
\end{center}
\end{figure}
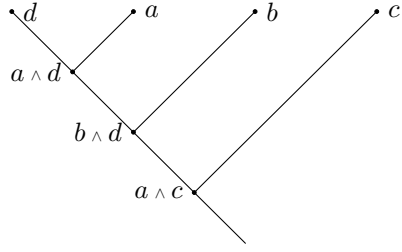

The class of finite meet{\hyp}trees forms a uniformly locally finite strong Fra\"{\i}ss\'{e} class in the language of meet{\hyp}trees. Its Fra\"{\i}ss\'{e} limit is the \emph{(countable) universal dense meet{\hyp}tree}, which we will always denote by $\Tt$. Thus, $\Tt$ is $\omega${\hyp}categorical, homogeneous and has quantifier elimination in the language of meet{\hyp}trees. See \cite[Fact 2.16]{kaplan2020automorphism}. 

The \emph{branched (countable) universal dense meet{\hyp}tree} is the (countable) universal dense meet{\hyp}tree $\Tt$ with a distinguished branch $\Gamma$ and its branch{\hyp}projection. The \emph{pointed (countable) universal dense meet{\hyp}tree} is the (countable) universal dense meet{\hyp}tree $\Tt$ with a distinguished point $\gamma$. Their uniqueness is justified by the following straightforward lemma:
\begin{lem}\label{l:fraisse trees} 
The class of finite pointed/semibranched meet{\hyp}trees is a uniformly locally finite strong Fra\"{\i}ss\'{e} class in the language of pointed/semibranched meet{\hyp}trees. Moreover, its Fra\"{\i}ss\'{e} limit is (isomorphic to) the branched/pointed universal dense meet{\hyp}tree. In particular, these structures are $\omega${\hyp}categorical, homogeneous and have quantifier elimination in their corresponding languages.
\begin{proof} The fact that these classes are strong Fra\"{\i}ss\'{e} classes is straightforward. The ``moreover'' part is a particular case of \cref{l:generic meet-tree expansions} below.
\end{proof}
\end{lem}
\begin{obs} In light of \cref{l:fraisse trees}, the pointed/branched universal dense meet{\hyp}tree is actually the universal dense pointed/branched meet{\hyp}tree.
\end{obs}

\section{Cone-independence} \label{s:section 2}
\begin{defi}[Cone{\hyp}independece] \label{d:cone-independence}
Let $T$ be a meet{\hyp}tree, $\gamma$ a point and $A,B,C$ subsets. We say that $A$ is \emph{$\gamma${\hyp}cone{\hyp}independent} of $B$ over $C$, denoted $A\indepe[\gamma]{C}B$, if 
\begin{enumerate}[label={\rm{(\roman*)}}, wide]
\item \label{itm:def cone-independence i} for any $a\in\langle AC\rangle_\gamma$ and $b\in\langle BC\rangle_\gamma$ with $b\leq a$, there is $c\in \langle C\rangle_\gamma$ such that $b\leq c\leq a$; and
\item \label{itm:def cone-independence ii} for any $a\in\langle AC\rangle_\gamma$ and $b\in \langle BC\rangle_\gamma$ with $a\meet b\notin\langle C\rangle_\gamma$ there is $c\in \langle C\rangle_\gamma$ such that $a\meet c\neq b\meet c$.
\end{enumerate}
\end{defi}
We apply the independence relation to tuples in the usual way.
\begin{obs} \label{o:reflexivity}  Obviously, $A\indepe[\gamma]{C}B$ if and only if $\langle AC\rangle_\gamma \indepe[\gamma]{\langle C\rangle_\gamma}\langle BC\rangle_\gamma$. Also, note that $A\indepe[\gamma]{C}B$ implies $\langle AC\rangle_\gamma\cap\langle BC\rangle_\gamma= \langle C\rangle_\gamma$ by point \ref{itm:def cone-independence i} of \cref{d:cone-independence}.
\end{obs}
\begin{figure}[ht]
\begin{center}
\begin{tikzpicture}[scale=0.8]
\node (xL) at (-6,-3) {};  \node (yL) at (-6,3) {};  
\node (xR) at (6,-3) {}; \node (yR) at (6,3) {};


\node (gamma) [draw, shape=circle, fill, scale=0.15, label={right:$\gamma$}] at (0,0) {}; 
\node (bot) at (0,-3) {}; 
\node (a1) [draw, shape=circle, fill, scale=0.15, label=right:$a_1$] at (-2.5,1.5) {}; 
\node (b1) [draw, shape=circle, fill, scale=0.15, label=right:$b_1$] at (0.5,2) {}; 
\node (b2) [draw, shape=circle, fill, scale=0.15, label=left:$b_2$] at (0,-1) {}; 
\node (a2) [draw, shape=circle, fill, scale=0.15, label=left:$a_2$] at (0,-1.5) {};
\node (pic) [draw, shape=coordinate, label={left:$c\meet \gamma$}] at (0,-2) {};
\node (b4) [draw, shape=circle, fill, scale=0.15, label=left:$b_4$] at (0,-2.5) {}; 
\node (b3) [draw, shape=circle, fill, scale=0.15, label=right:$b_3$] at (1,-0.5) {}; 
\node (a3) [draw, shape=circle, fill, scale=0.15, label=right:$a_3$] at (2.5,-1) {};
\node (a3meetb3) [draw, shape=coordinate, label={below:$a_3\meet b_3$}] at (1,-1.75) {};
\node (c) [draw, shape=circle, fill, scale=0.15, label={above:$c$}] at (1.6,-1.45) {};

\draw (a1) -- (gamma);
\draw (b1) -- (gamma); 
\draw (gamma) -- (bot);
\draw (pic) -- (a3meetb3); 
\draw (a3meetb3) -- (a3); 
\draw (a3meetb3) -- (b3);
\end{tikzpicture}
\caption{Example of $\bar{a}\indepe[\gamma]{c}\bar{b}$.}
\end{center}
\end{figure}
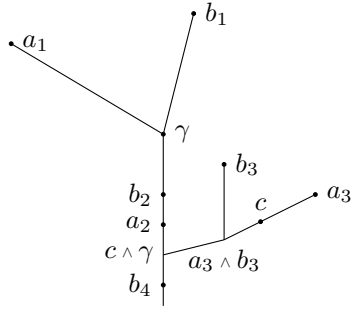
 
The name $\gamma${\hyp}cone{\hyp}independence is justified by the following lemma:
\begin{lem}\label{l:cones and cone independence} Let $T$ be a meet{\hyp}tree, $\gamma$ a point and $A,B,C$ subsets of $T$ such that $A\indepe[\gamma]{C}B$. Then, $\cone_\gamma(AC)\cap \cone_\gamma(BC)=\cone_\gamma(C)$.
\begin{proof} It is enough to show that if $a\in A$ and $b\in B$ such that $a\in\cone_\gamma(b)$, then $a\in\cone_\gamma(C)$. As $A\indepe[\gamma]{C}B$, either $a\meet b\in\langle C\rangle_\gamma$ or there is $c\in \langle C\rangle_\gamma$ such that $a\meet c\neq b\meet c$. In the first case, $a\meet b\in\langle C\rangle_\gamma$ and $a\in\cone_\gamma(a\meet b)$. In the second case, by \cref{l:three elements}, $\gamma<a\meet b=a\meet b\meet c\leq a\meet c$, so $a\in\cone_\gamma(c)$.
\end{proof}
\end{lem}

Cone{\hyp}independence relations for different base points are related in the following obvious way:
\begin{lem} \label{l:no based cone-independence definition}
Let $T$ be a meet{\hyp}tree and $A,B,C$ subsets of $T$. Then, $A\indepe[\gamma]{C}B$ for some $\gamma\in C$ if and only if $A\indepe[\gamma]{C}B$ for any $\gamma\in C$.
\begin{proof} It is enough to note that, for $\gamma\in C$, we have $\langle AC\rangle_\gamma=\langle AC\rangle_\wedge$, $\langle BC\rangle_\gamma=\langle BC\rangle_\wedge$ and $\langle C\rangle_\gamma=\langle C\rangle_\wedge$ by \cref{l:generated}(\ref{itm:generated 1}). Thus, the \cref{d:cone-independence} of $A\indepe[\gamma]{C}B$ is independent of $\gamma$ when $\gamma\in C$.
\end{proof} 
\end{lem}

\begin{lem} \label{l:cone-independence and generated} Let $T$ be a meet{\hyp}tree and $\gamma$ a point. Then, $\indepe[\gamma]{}$ has the \emph{independent generation property}: 
\[A\indepe[\gamma]{C}B\Rightarrow\langle ABC\rangle_\gamma=\langle AC\rangle_\gamma\cup\langle BC\rangle_\gamma.\]
\begin{proof} Without loss of generality, we can assume that $A=\langle AC\rangle_{\gamma}$, $B=\langle BC\rangle_{\gamma}$ and $C=\langle C\rangle_\gamma$. We want to show that $\langle AB\rangle_\wedge=A\cup B$. By \cref{l:generated}(\ref{itm:generated 3}), we only need to check that $a\meet b\in A\cup B$ for any $a\in A$ and $b\in B$. As $A\indepe[\gamma]{C} B$, either $a\meet b\in C\subseteq A\cup B$, or there is $c\in C$ such that $a\meet c\neq b\meet c$. Now, in the latter case, by \cref{l:three elements}, $a\meet b=a\meet b\meet c=\min\{a\meet c,b\meet c\}\in A\cup B$.
\end{proof}
\end{lem}
\begin{lem} \label{l:cone-independence finite} Let $T$ be a meet{\hyp}tree and $\gamma$ a point. Then, $\indepe[\gamma]{}$ satisfies the following properties:
\begin{enumerate}[wide]
\item[\rm{(Quantifier free definability)}] For any finite variables $x,y,z$, there is a quantifier free formula $\phi(x,y,z)$ in $\lang_\gamma$ such that, for any tuples $a,b,c$ (on the corresponding variables), $a\indepe[\gamma]{c}b$ if and only if $\vDash \phi(a,b,c)$.
\item[\rm{(Quantifier free invariance)}] If $a\indepe[\gamma]{c}b$, then $a'\indepe[\gamma]{c'}b'$  for any $a',b',c'$ such that $\qftp_{\gamma}(a',b',c')=\qftp_{\gamma}(a,b,c)$. 
\item[\rm{(Embedding invariance)}] If $f:\ T\rightarrow T'$ is an embedding of meet{\hyp}trees, then $A\indepe[\gamma]{C}B$ if and only if $f(A)\indepe[f(\gamma)]{f(C)}f(B)$.
\item[\rm{(Monotonicity)}] If $A\indepe[\gamma]{C}B$ and $A_0\subseteq A$ and $B_0\subseteq B$, then $A_0\indepe[\gamma]{C}B_0$.
\item[\rm{(Normality)}] If $A\indepe[\gamma]{C}B$, then $AC\indepe[\gamma]{C}BC$.
\item[\rm{(Left and right existence)}] $A\indepe[\gamma]{A}B$ and $A\indepe[\gamma]{B}B$.
\item[\rm{(Finite character)}] $A\indepe[\gamma]{C}B$ if and only if $A_0\indepe[\gamma]{C}B_0$ for every $A_0\subseteq A$ and $B_0\subseteq B$ finite.
\item[\rm{(Left and right base monotonicity)}] Suppose $A\indepe[\gamma]{C}B$ and $C\subseteq D\subseteq A$, then $A\indepe[\gamma]{D}B$. If $A\indepe[\gamma]{C}B$ and $C\subseteq D\subseteq B$, then $A\indepe[\gamma]{D}B$.
\item[\rm{(Left and right transitivity)}] If $A\indepe[\gamma]{D}B$ and $D\indepe[\gamma]{C}B$ with $C\subseteq D\subseteq A$, then $A\indepe[\gamma]{C}B$. If $A\indepe[\gamma]{D}B$ and $A\indepe[\gamma]{C}D$ with $C\subseteq D\subseteq B$, then $A\indepe[\gamma]{C}B$.
\item[\rm{(Quantifier free stationarity)}] Suppose $a\indepe[\gamma]{C}b$ and $a'\indepe[\gamma]{C}b'$ with $\qftp_{\gamma}(a'/C)$ $=\qftp_{\gamma}(a/C)$ and $\qftp_{\gamma}(b'/C)=\qftp_{\gamma}(b/C)$, then $\qftp_{\gamma}(a'b'/C)=\qftp_{\gamma}(ab/C)$.
\end{enumerate}
\begin{proof} Quantifier free definability is obvious from \cref{d:cone-independence}. Quantifier free invariance and embedding invariance are trivial consequences of quantifier free definability. Monotonicity, normality, existence and finite character are trivial. It remains to check base monotonicity, transitivity and quantifier free stationarity:
\begin{enumerate}[wide]
\item[\rm{(Base monotonicity)}] Let $A\indepe[\gamma]{C}B$ and $C\subseteq D\subseteq B$. We want to show that $A\indepe[\gamma]{C'}B$. Note that, by monotonicity, $A\indepe[\gamma]{C}D$. By  \cref{l:cone-independence and generated}, $\langle AD\rangle_\gamma= \langle AC\rangle_\gamma\cup \langle D\rangle_\gamma$. Then:
\begin{enumerate}[label={\rm{(\roman*)}}, wide]
\item Pick $a\in\langle AD\rangle_\gamma=\langle AC\rangle_\gamma\cup\langle D\rangle_\gamma$ and $b\in\langle B\rangle_\gamma$ such that $b\leq a$. If $a\in \langle D\rangle_\gamma$, then $b\leq a\leq a$ with $a\in\langle D\rangle_\gamma$. Otherwise, $a\in\langle AC\rangle_\gamma$ and, as $A\indepe[\gamma]{C}B$, there is $c\in\langle C\rangle_\gamma\subseteq\langle D\rangle_\gamma$ such that $b\leq c\leq a$.
\item Pick $a\in\langle AD\rangle_\gamma=\langle AC\rangle_\gamma\cup\langle D\rangle_\gamma$ and $b\in\langle B\rangle_\gamma$ such that $a\meet b\notin \langle D\rangle_\gamma$. If $a\in\langle D\rangle_\gamma$, we get that $a\meet a\in\langle D\rangle_\gamma$ and $a\meet b\notin \langle D\rangle_\gamma$, so $a\meet a\neq a\meet b$ with $a\in \langle D\rangle_\gamma$. Otherwise, $a\in\langle AC\rangle_\gamma$ and, as $A\indepe[\gamma]{C}B$, there is $c\in\langle C\rangle_\gamma\subseteq\langle D\rangle_\gamma$ such that $a\meet c\neq b\meet c$.
\end{enumerate}

This proves right base monotonicity, and left base monotonicity follows by the symmetric argument.
\item[\rm{(Transitivity)}] Let $A\indepe[\gamma]{D}B$ and $A\indepe[\gamma]{C}D$ with $C\subseteq D\subseteq B$. We want to show that $A\indepe[\gamma]{C}B$. 
\begin{enumerate}[label={\rm{(\roman*)}}, wide]
\item Pick $a\in\langle AC\rangle_\gamma\subseteq \langle AD\rangle_\gamma$ and $b\in\langle BC\rangle_\gamma\subseteq \langle BD\rangle_\gamma$ with $b\leq a$. As $A\indepe[\gamma]{D}B$, there is $d\in \langle D\rangle_\gamma$ such that $b\leq d\leq a$. Now, since $A\indepe[\gamma]{C}D$, there is $c\in\langle C\rangle_\gamma$ such that $d\leq c\leq a$, concluding $b\leq c\leq a$.
\item Pick $a\in\langle AC\rangle_\gamma\subseteq \langle AD\rangle_\gamma$ and $b\in\langle BC\rangle_\gamma\subseteq \langle BD\rangle_\gamma$ such that $a\meet c= b\meet c$ for all $c\in \langle C\rangle_\gamma$. Either $a\meet d=b\meet d$ for all $d\in\langle D\rangle_\gamma$ or there is $d\in\langle D\rangle_\gamma$ such that $a\meet d\neq b\meet d$. 

In the first case, as $A\indepe[\gamma]{D}B$, we conclude $a\meet b\in\langle D\rangle_\gamma$. Now, for any $c\in\langle C\rangle_\gamma$, as $a\meet c=b\meet c$, using \cref{l:generated}(\ref{itm:generated 2}), we get that $a\meet c=(a\meet c)\meet (b\meet c)=(a\meet b)\meet c$. As $A\indepe[\gamma]{C}D$, we conclude that $a\meet b=a\meet (a\meet b)\in \langle C\rangle_\gamma$.

In the second case, by \cref{l:four elements}, $a\meet c=d\meet c$ for any $c\in\langle C\rangle_\gamma$. Therefore, as $A\indepe[\gamma]{C}D$, we conclude that $a\meet d\in\langle C\rangle_\gamma$. In particular, $a\meet d=a\meet (a\meet d)=b\meet (a\meet d)=a\meet b\meet d$. Since $a\meet d\neq b\meet d$, by \cref{l:three elements}, we conclude that $a\meet b=a\meet b\meet d=a\meet d\in\langle C\rangle_\gamma$.
\end{enumerate}

This proves right transitivity, and left transitivity follows by the symmetric argument.
\item[\rm{(Quantifier free stationarity)}] By embedding invariance it suffices to consider the case of the universal dense meet{\hyp}tree. By quantifier free invariance, it is enough to show that if $a \indepe[\gamma]{C} B$ and $a'\indepe[\gamma]{C}B$ with $\qftp(a/C)= \qftp(a'/C)$, then $\qftp_\gamma(a/CB)=\qftp_\gamma(a'/CB)$. It is enough to do the case where $|a|=1$ by monotonicity, base monotonicity and invariance. 

Without loss of generality, say $C=\langle C\rangle_\gamma$ and $B=\langle BC\rangle_\gamma$. There are two cases: there is $c\in C$ with $a\leq c$ or $e\coloneqq \max\{a\meet c\tq c\in C\}<a$. 

In the first case, take $c\in C$ with $a\leq c$. Since $\qftp(a/C)=\qftp(a'/C)$, $a'\leq c$. By \cref{l:no based cone-independence definition}, $a\indepe[c]{C}B$ and $a'\indepe[c]{C}B$. By point \ref{itm:def cone-independence i} of \cref{d:cone-independence}, $\qftp_{<}(a/\pi(B))=\qftp_{<}(a'/\pi(B))$ where $\pi(B)=\{b\meet c\tq b\in B\}=\{b\in B\tq b\leq c\}$. Then, by \cref{c:semibranch intervals are indiscernible}, $\qftp_\Gamma(a/B)=\qftp_\Gamma(a'/B)$ where $\Gamma$ is the semibranch given by $c$. By \cref{c:quantifier free type tree}(\ref{itm:quantifier free type tree 1},\ref{itm:quantifier free type tree 3}), $\qftp_\gamma(a/B)=\qftp_\gamma(a'/B)$.

In the second case, let $e\coloneqq\max\{a\meet c\tq c\in C\}$ and $e'\coloneqq\max\{a'\meet c\tq c\in C\}$. Since $\qftp_\gamma(a/C)=\qftp_\gamma(a'/C)$, we have $\qftp_\gamma(e/C)=\qftp_\gamma(e'/C)$. Also, by monotonicity, we have $e\indepe[\gamma]{C}B$ and $e'\indepe[\gamma]{C}B$, but now there is $c\in C$ with $e\leq c$. Thus, by the previous case, $\qftp_\gamma(e/B)=\qftp_\gamma(e'/B)$. By invariance, we may assume $e=e'$. By base monotonicity, $a\indepe[\gamma]{Ce}B$ and $a'\indepe[\gamma]{Ce}B$ with $\qftp_\gamma(a/Ce)=\qftp_\gamma(a'/Ce)$. By \cref{l:no based cone-independence definition}, $a\indepe[e]{Ce}B$ and $a'\indepe[e]{Ce}B$. By choice of $e$, $a,a'\notin \cone_e(C)$. By \cref{l:cones and cone independence}, $a,a'\notin\cone_e(B)$. Hence, by \cref{c:cones are indiscernible}, $\qftp_e(a/B)=\qftp_e(a'/B)$. By \cref{c:quantifier free type tree}(\ref{itm:quantifier free type tree 1}), $\qftp_\gamma(a/B)=\qftp_\gamma(a'/B)$.
\qedhere
\end{enumerate}
\end{proof}
\end{lem}

\begin{lem}\label{l:easy independence} Let $T$ be a meet{\hyp}tree, $\gamma$ a point and $A,B,C$ subsets. Assume that, for all $a\in A$ and $b\in B$, we have $a\meet b\leq \gamma$ when $a\geq \gamma$ and $a\meet \gamma<b\meet \gamma$ when $a\meet \gamma<\gamma$. Then, $A\indepe[\gamma]{C}B$.
\begin{proof} Take $b,b'\in B\cup\{\gamma\}$ and $a\in A$. If $a\geq \gamma$, clearly $a\meet (b\meet b')\leq a\meet b\leq\gamma$. If $a\meet \gamma<\gamma$, $a\meet\gamma<\min\{b\meet \gamma,b'\meet\gamma\}=b\meet b'\meet \gamma$. Therefore, we can assume without loss of generality that $B=\langle B\rangle_\gamma$.  

Note that the condition does not mention $C$. Thus, by transitivity (and induction), we can assume that $A=\{a\}$ is a singleton. Also, we can assume without loss of generality that $C=\langle C\rangle_\gamma$. There are two cases, $a\geq \gamma$ or $a\meet \gamma<\gamma$.
\begin{enumerate}[wide]
\item[Case $a\geq \gamma$.] By existence, we can assume $a>\gamma$. 

First, we prove condition \ref{itm:def cone-independence i} of \cref{d:cone-independence}. Let $a'\in \langle a\, C\rangle_\gamma$, and $b'\in\langle BC\rangle_\gamma$ with $b'\leq a'$. If $a'\in C$ or $b'\in C$, there is nothing to prove. By \cref{l:generated}(\ref{itm:generated 3}), we know that $a'=a\meet c_a$ and $b'=b\meet c_b$ with $b\in B$ and $c_a\in C\cup\{a\}$ and $c_b\in C\cup\{b\}$. Thus, $b'\leq a$, $b'\leq b$ and $b'\leq c_a$, so $b'\leq a\meet b\meet c_a\leq \gamma\meet c_a\leq a\meet c_a=a'$ by hypothesis. If $c_a\in C$,  we are done. If $c_a=a$, $\gamma\meet c_a=\gamma\in C$. 

Now, we prove condition \ref{itm:def cone-independence ii} of \cref{d:cone-independence}. Take $a'\in\langle a\, C\rangle_\gamma$ and $b'\in \langle BC\rangle_\gamma$. If $a'\in C$ or $b'\in C$, there is nothing to prove. By \cref{l:generated}(\ref{itm:generated 3}), we know that $a'=a\meet c_a$ and $b'=b\meet c_b$ with $b\in B$ and $c_a\in C\cup\{a\}$ and $c_b\in C\cup\{b\}$. By hypothesis, $a'\meet b'\leq a\meet b\meet c_a\leq \gamma\meet c_a$. If $a'\meet b'=\gamma\meet c_a$, we are done. Otherwise, $a'\meet (\gamma\meet c_a)=(a\meet c_a)\meet(\gamma\meet c_a)=(a\meet \gamma)\meet c_a=\gamma\meet c_a$ and $b'\meet (\gamma\meet c_a)=b'\meet (a'\meet \gamma \meet c_a)=(a'\meet b')\meet(\gamma\meet c_a)=a'\meet b'<\gamma\meet c_a=a'\meet (\gamma\meet c_a)$.
\item[Case $a\meet \gamma<\gamma$.] \

First, we prove condition \ref{itm:def cone-independence i} of \cref{d:cone-independence}. Let $a'\in \langle a\, C\rangle_\gamma$, and $b'\in\langle BC\rangle_\gamma$ with $b'\leq a'$. If $a'\in C$ or $b'\in C$, there is nothing to prove. By \cref{l:generated}(\ref{itm:generated 3}), we know that $a'=a\meet c_a$ and $b'=b\meet c_b$ with $b\in B$ and $c_a\in C\cup\{a\}$ and $c_b\in C\cup\{b\}$. Then, by hypothesis, $\min\{b\meet \gamma,c_b\meet \gamma\}=b'\meet \gamma\leq a\meet \gamma<b\meet \gamma$, so $c_b\meet \gamma<b\meet\gamma$ (in particular, $c_b\in C$). By \cref{l:three elements} applied to $b,c_b,\gamma$, we get $b'=b\meet c_b=c_b\meet \gamma\in C$ and we are done.

Now, we prove condition \ref{itm:def cone-independence ii} of \cref{d:cone-independence}. Take $a'\in \langle a\, C\rangle_\gamma$ and $b'\in\langle BC\rangle_\gamma$ with $a'\meet c=b'\meet c$ for all $c\in C$. If $a'\in C$ or $b'\in C$, there is nothing to prove. By \cref{l:generated}(\ref{itm:generated 3}), we know that $a'=a\meet c_a$ and $b'=b\meet c_b$ with $b\in B$ and $c_a\in C\cup\{a\}$ and $c_b\in C\cup\{b\}$. By assumption, $a\meet \gamma<b\meet \gamma$. As $a'\meet\gamma=b'\meet\gamma$, we get $b'\meet \gamma=a'\meet \gamma\leq a\meet \gamma<b\meet \gamma$. In particular, $c_b\neq b$, so $c_b\in C$. By \cref{l:three elements} applied to $b,c_b,\gamma$, we get that $b'\meet \gamma=b\meet c_b=c_b\meet \gamma\in C$. Therefore, $a'\meet b'=a'\meet (b\meet c_b)=a'\meet (b'\meet \gamma)=(a'\meet \gamma)\meet (b'\meet \gamma)=b'\meet \gamma=c_b\meet \gamma\in C$.\qedhere
\end{enumerate}
\end{proof}
\end{lem}

In this paper, following \cite{li2019automorphism}, we use the following terminology.

\begin{defi}\label{d:stationary weak independence relation} A \emph{stationary weak independence relation} is a ternary relation between finite subsets satisfying invariance, monotonicity, normality, base monotonicity, transitivity, stationarity, existence and extension:
\begin{enumerate}[wide]
\item[\rm{(Left and right extension)}] If $a\indepe{C}b$ and $a\subseteq a'$, then there is $a''\equiv_Ca'$ with $a''\indepe{C}b$. If $a\indepe{C}b$ and $b\subseteq b'$, then there is $b''\equiv_Cb'$ with $a\indepe{C}b''$. 
\end{enumerate}
\end{defi}
\begin{obs}\label{o:left extension implies right extension} Note that right extension follows from left extension, transitivity, normality, existence, monotonicity and invariance. Indeed, let $a\indepe{c}b$ with $b\subseteq b'$. By normality, we have $a\indepe{c}bc$. By existence and monotonicity, we have $\emptyset\indepe{bc}b'$, so there is $a'$ with $a'\indepe{bc}b'$ and $a'\equiv_{bc}a$ by left extension. Take $b''$ such that $a',b'\equiv_{bc} a,b''$. By invariance, $a\indepe{bc}b''$, concluding $a\indepe{c}b''$ by transitivity. 
\end{obs}

\begin{lem} \label{l:cone-independence} $\indepe[\gamma]{}$ in $\Tt_\gamma$ is a stationary weak independence relation satisfying reflexivity and the independent generation property.
\begin{proof} Reflexivity follows from \cref{o:reflexivity}. By \cref{l:cone-independence and generated}, we already know that $\indepe[\gamma]{}$ satisfies the independent generation property. By quantifier elimination and \cref{l:cone-independence finite}, we already have invariance, monotonicity, normality, existence, base monotonicity, transitivity and stationarity. It remains to show extension. By \cref{o:left extension implies right extension}, it suffices to show left extension:
\begin{claim} \label{cl:claim extension} Let $C$ and $B$ be finite subset with $C=\langle C\rangle_\gamma$ and $B=\langle BC\rangle_\gamma$. Let $a$ be a single element such that $\langle a\, C\rangle_\gamma=\langle C\rangle_\gamma\cup\{a\}$. Then, there is $a'\indepe[\gamma]{C}B$ with $a'\equiv^\gamma_Ca$.
\begin{proof} The case where $a\in C$ is trivial, so assume $a\notin C$. Write $C_{<a}\coloneqq\{c\in C\tq c<a\}$ and $C_{>a}\coloneqq\{c\in C\tq c>a\}$. There are three possible cases:
\begin{enumerate}[label={\rm{Case \arabic*:}}, wide]
\item $C_{<a}=\emptyset$. Since we assume that $\langle a\, C\rangle_\gamma=C\cup\{a\}$, we have that $a\meet c\in C$ or $a\leq c$ for all $c\in C$. Thus, if $C_{<a}=\emptyset$, we get that $a< c$ for all $c\in C$, i.e. $C=C_{>a}$. Take $b_0\coloneqq\bigwedge B$ and $a'<b_0$ arbitrary. Since $B=\langle B\rangle_\gamma$, we have that $a'<b_0\leq \gamma$, so $a'\meet \gamma=a'<b_0\leq b\meet \gamma$ for any $b\in B$. Thus, by \cref{l:easy independence}, we conclude that $a'\indepe[\gamma]{C}B$. On the other hand, since $C\subseteq B$, we have that $a'<c$ for all $c\in C$. Thus, $\qftp_{<}(a/C)=\qftp_{<}(a'/C)$. As $a'\meet c=a'$ for all $c\in C=\langle C\rangle_\gamma$, we have that $\langle a'\, C\rangle_\gamma=C\cup \{a'\}$ by \cref{l:generated}(\ref{itm:generated 4}). Thus, we conclude that $a\equiv^\gamma_C a'$ by \cref{l:quantifier free type} and quantifier elimination.
\item $C_{>a}=\emptyset$. Since we assume that $\langle a\, C\rangle_\gamma=C\cup\{a\}$, we have that $e\coloneqq\max\{a\meet c\tq c\in C\}\in C$. Pick any $a'>e$ with $a'\notin \cone_e(B)$, which exists by finiteness of $B$. By \cref{l:easy independence}, $a'\indepe[e]{C}B$. By \cref{l:no based cone-independence definition}, as $\gamma,e\in C$, we conclude that $a'\indepe[\gamma]{C}B$. On the other hand, since $C\subseteq B$, $a'\notin\cone_e(C)$. Also, by choice of $e$, $a\notin\cone_e(C)$. Therefore, by \cref{c:cones are indiscernible} and quantifier elimination, we conclude that $a\equiv^\gamma_Ca'$.  
\item $C_{<a}\neq \emptyset$ and $C_{>a}\neq\emptyset$. Take $c_0\coloneqq\max C_{<a}$ and $c_1\coloneqq \min C_{>a}=\bigwedge C_{>a}$. Now, $c_0<a<c_1$ and there is no $c\in C$ with $c_0<c<c_1$. Take $B_0\coloneqq\{b\in B\tq c_0<b\leq c_1\}$ and $b_0\coloneqq\min B_0$. Pick any $a'$ with $c_0<a'<b_0$. Obviously, $a'\meet c_1=a'=\max\{a'\meet c\tq c\in C\}$, so $\langle a'\, C\rangle_\gamma=C\cup\{a'\}$ by \cref{l:generated}(\ref{itm:generated 4}). Since there is no $c$ with $c_0<c<c_1$, we get that $C_{<a}=\{c\in C\tq c<a'\}$ and $C_{>a}=\{c\in C\tq a'<c\}$. Hence, $\qftp_<(a/C)=\qftp_<(a'/C)$, concluding that $a\equiv^\gamma_Ca'$ by \cref{l:quantifier free type} and quantifier elimination. It remains to check that $a'\indepe[\gamma]{C}B$. 
\begin{enumerate}[label={\rm(\roman*)}, wide]
\item Take $b\in B$ such that $b<a'$. Then, $b<b_0$, so $b\notin B_0$, concluding $b\leq c_0<a'$. 
\item Take $b\in B$. Either $b\meet c_1\leq c_0<a'=a'\meet c_1$ or $a'\meet c_1=a'<b_0\leq b\meet c_1$. In any case, $a'\meet c_1\neq b\meet c_1$ where $c_1\in C$.  \claimqed
\end{enumerate}
\end{enumerate}
\end{proof}
\end{claim}

Now, suppose $A\indepe[\gamma]{C}B$ with $A,B,C$ finite and $a$ an arbitrary point. Without loss of generality, $A=\langle AC\rangle_\gamma$ and $B=\langle BC\rangle_\gamma$. By \cref{l:generated}(\ref{itm:generated 4}), $\langle a\, A\rangle_\gamma=A\cup\{e,a\}$ with $e\coloneqq\max\{a\meet d\tq d\in A\}$; in particular, $\langle e\, A\rangle_\gamma=A\cup\{e\}$. By \cref{cl:claim extension}, there is $e'$ such that $e'\indepe[\gamma]{A}B$ and $e'\equiv^\gamma_A e$. Set $B'\coloneqq\langle B\, e'\rangle_\gamma$. Take $a'$ such that $a'e'\equiv^\gamma_A ae$, so $\langle e'a'A\rangle_\gamma=A\cup \{e',a'\}$. Again, by \cref{cl:claim extension}, there is $a''$ such that $a''\indepe[\gamma]{Ae'}B'$ with $a''\equiv^\gamma_{Ae'}a'$, so $a''\equiv^\gamma_{A} a$. By monotonicity, $a''\indepe[\gamma]{Ae'}B$, by transitivity, $a''e'\indepe[\gamma]{A}B$ and $a''e'A\indepe[\gamma]{C}B$, so $a''A\indepe[\gamma]{C}B$ by monotonicity. Hence, we conclude that $\indepe[\gamma]{}$ satisfies left extension.  
\end{proof}
\end{lem}
Imitating \cite{tent2013isometry}, a \emph{local stationary weak independence relation} is an independence relation between finite subsets, defined only for non{\hyp}empty base sets, that satisfies invariance, monotonicity, normality, base monotonicity, transitivity, stationarity, existence and extension.
For $A,B,C$ with $C\neq\emptyset$, we say that $A$ is \emph{cone{\hyp}independent of $B$ over $C$}, denoted by $A\indepe[\cone]{C}B$, if $A\indepe[\gamma]{C}B$ for some (any) $\gamma\in C$. 
By \cref{l:no based cone-independence definition,l:cone-independence}, we immediately get the following corollary.
\begin{coro} \label{c:cone-independence without base point} $\indepe[\cone]{}$ is a local stationary weak independence relation satisfying reflexivity and the independent generation property in $\Tt$.
\end{coro}
%

\begin{defi}[Semibranch{\hyp}independence] \label{d:semibranch-independence}
Let $T$ be a meet{\hyp}tree, $\Gamma$ a semibranch and $A,B,C$ subsets. We say that $A$ is \emph{$\Gamma${\hyp}semibranch{\hyp}independent} of $B$ over $C$, denoted $A\indepe[\Gamma]{C}B$, if 
\begin{enumerate}[label={\rm{(\roman*)}}, ref={\rm{\roman*}}, wide]
\item \label{itm:semibranch-independence i} for any $a\in\langle AC\rangle_\Gamma$ and $b\in\langle BC\rangle_\Gamma$ with $b\leq a$, there is $c\in \langle C\rangle_\Gamma$ such that $b\leq c\leq a$; and
\item \label{itm:semibranch-independence ii} for any $a\in\langle AC\rangle_\Gamma$ and $b\in \langle BC\rangle_\Gamma$ with $\pi_\Gamma(a)=\pi_\Gamma(b)$ and $a\meet b\notin\langle C\rangle_\Gamma$, there is $c\in \langle C\rangle_\Gamma$ such that $a\meet c\neq b\meet c$. 
\end{enumerate}
\end{defi}
We apply the independence relation to tuples in the usual way.
\begin{obs} Obviously, $A\indepe[\Gamma]{C}B$ if and only if $\langle AC\rangle_\Gamma\indepe[\Gamma]{\langle C\rangle_\Gamma}\langle BC\rangle_{\Gamma}$. Also, note that $A\indepe[\Gamma]{C}B$ imples $\langle AC\rangle_\Gamma\cap\langle BC\rangle_\Gamma=\langle C\rangle_\Gamma$ by point \cref{d:semibranch-independence}(\ref{itm:semibranch-independence i}).
\end{obs}

It is clear from the definitions that semibranch{\hyp}independence and cone{\hyp}independence are very similar. We are going to make this similarity rigorous. First, we consider the case of semibranches that are not branches:
\begin{lem}\label{l:semibranch and cone independence} Let $T$ be a meet{\hyp}tree, $\gamma$ a point, $\Gamma=T_{\leq\gamma}$ and $A,B,C$ subsets. Then, $A\indepe[\gamma]{C}B$ if and only if $A\indepe[\Gamma]{C\gamma}B$. Furthermore, if there is no $a\in\langle AC\rangle_{\wedge}$ such that $a\geq\gamma$, then $A\indepe[\gamma]{C}B$ if an only if $A\indepe[\Gamma]{C}B$.
\begin{proof} By \cref{l:generated}(\ref{itm:generated 4},\ref{itm:generated 7}), note that $\langle X\rangle_{\gamma}=\langle X\gamma\rangle_{\wedge}=\langle X\,\gamma\rangle_{\Gamma}=\langle X\rangle_{\Gamma}\cup \{\gamma\}$ for any $X\subseteq T$. Thus, it obviously follows from the definitions that $A\indepe[\gamma]{C}B$ if and only if $A\indepe[\Gamma]{C\gamma}B$.

Assume now that there is no $a\in\langle AC\rangle_{\wedge}$ with $a\geq\gamma$. We want to show that $A\indepe[\Gamma]{C\gamma}B$ if and only if $A\indepe[\Gamma]{C}B$. We first explain why $A\indepe[\Gamma]{C}B$ implies $A\indepe[\Gamma]{C\gamma}B$ (this does not use the assumption).

\begin{enumerate}[label={\rm{(\roman*)}}, wide]
\item Take $a\in \langle AC\,\gamma\rangle_{\Gamma}$ and $b\in \langle BC\,\gamma\rangle_{\Gamma}$ with $b\leq a$. If $a\neq \gamma$ and $b\neq \gamma$, as $A\indepe[\Gamma]{C}B$, we are done.
If $a=\gamma$ or $b=\gamma$, then $b\leq \gamma\leq a$.
\item Take $a\in \langle AC\,\gamma\rangle_{\Gamma}$ and $b\in \langle BC\,\gamma\rangle_{\Gamma}$ with $\pi_{\Gamma}(a)=\pi_{\Gamma}(b)$ and $a\meet b\notin \langle C\,\gamma\rangle_{\Gamma}$. If $a=\gamma$, then $\pi_\Gamma(b)=\pi_\Gamma(\gamma)=\gamma$, so $b\geq\gamma$ and $a\meet b=\gamma\in\langle C\,\gamma\rangle_\Gamma$. The case when $b=\gamma$ is symmetric. Now, if $a\neq \gamma$ and $b\neq\gamma$, we are done by assumption.
\end{enumerate}

Conversely, suppose $A\indepe[\Gamma]{C\gamma}B$; we want to show that $A\indepe[\Gamma]{C}B$:
\begin{enumerate}[label={\rm{(\roman*)}}, wide]
\item Take $a\in \langle AC\rangle_{\Gamma}$ and $b\in \langle BC\rangle_{\Gamma}$ with $b\leq a$. As $A\indepe[\Gamma]{C\gamma}B$, there is $c\in \langle C\,\gamma\rangle_{\Gamma}$ such that $b\leq c\leq a$. By hypothesis, $a\not\geq \gamma$, so $c\in\langle C\rangle_{\Gamma}$.
\item Take $a\in \langle AC\rangle_{\Gamma}$ and $b\in \langle BC\rangle_{\Gamma}$ with $\pi_{\Gamma}(a)=\pi_{\Gamma}(b)$ and $a\meet b\notin \langle C\rangle_{\Gamma}$. By hypothesis, $a\meet b\neq \gamma$, so $a\meet b\notin\langle C\,\gamma\rangle_{\Gamma}$. Since $A\indepe[\Gamma]{C\gamma}B$, we conclude that there is $c\in \langle C\gamma\rangle_{\Gamma}$ such that $a\meet c\neq b\meet c$. As $a\meet \gamma=\pi_{\Gamma}(a)=\pi_{\Gamma}(b)=b\meet \gamma$, we get that $c\in\langle C\rangle_{\Gamma}$.\qedhere
\end{enumerate}
\end{proof}
\end{lem}

For branches we use the following trivial remark: 
\begin{lem}\label{l:semibranched independence in substructures} Let $(T_1,\Gamma_1)$ be a semibranched meet{\hyp}tree and $(T_2,\Gamma_2)$ a substructure of $(T_1,\Gamma_1)$ as semibranched meet{\hyp}tree. Then, for any $A,B,C\subseteq T_2$, $A\indepe[\Gamma_2]{C}B$ if and only if $A\indepe[\Gamma_1]{C}B$. 
\begin{proof} For $x\in X\subseteq T_2$, $\pi_{\Gamma_2}(x)=\pi_{\Gamma_1}(x)$ and $\langle X\rangle_{\Gamma_2}=\langle X\rangle_{\Gamma_1}$. Therefore, the definitions of $\indepe[\Gamma_1]{}$ and $\indepe[\Gamma_2]{}$ agree on subsets of $T_2$. 
\end{proof}
\end{lem}

We get then the following conclusion:
\begin{coro} \label{c:semibranch and cone independence} Let $T$ be a meet{\hyp}tree, $\Gamma$ a branch and $A,B,C$ subsets of $T$. Let $T^{\infty}$ be the meet{\hyp}tree extending $T$ given by \cref{c:semibranches tree}. It follows that $(T,\Gamma)$ is a substructure of $(T^{\infty},\Gamma^{\infty})$ semibranched meet{\hyp}tree. Then, $A\indepe[\Gamma]{C}B$ in $T$ if and only if $A\indepe[\Gamma]{C}B$ in $T^{\infty}$.
\begin{proof} Consider $T^{\infty}$ with the semibranch $\Gamma^{\infty}\coloneqq T^{\infty}_{\subseteq \Gamma}$ (\cref{c:semibranches tree}). Then, by \cref{l:semibranch and cone independence,l:semibranched independence in substructures}, we conclude that $A\indepe[\Gamma]{C}B$ in $T^{\infty}$ if and only if $A\indepe[\Gamma^{\infty}]{C}B$, if and only if $A\indepe[\Gamma]{C}B$ in $T$.
\end{proof}
\end{coro}

\begin{coro} \label{c:semibranch-independence} For any branch $\Gamma$ of $\Tt$, $\indepe[\Gamma]{}$ is a stationary weak independence relation satisfying reflexivity and the independent generation property in $\Tt_\Gamma$.
\begin{proof} Take a countable elementary extension $\Tt'$ of $\Tt$ with an element $\gamma$ such that $\gamma>a$ for all $a\in \Gamma$, and write $\Gamma'=\Tt'_{\leq \gamma}$. By \cref{l:semibranch and cone independence} and \cref{l:semibranched independence in substructures}, for any $A,B,C\subseteq \Tt$, $A\indepe[\gamma]{C}B$ if and only if $A\indepe[\Gamma']{C}B$, if and only if $A\indepe[\Gamma]{C}B$. By $\omega${\hyp}categoricity and \cref{l:cone-independence}, $\indepe[\gamma]{}$ is a stationary weak independence relation satisfying the independent generation property in $\Tt'_{\gamma}$. Hence, $\indepe[\Gamma]{}$ trivially inherits monotonicity, normality, base monotonicity, transitivity and existence. Reflexivity is given by \cref{d:semibranch-independence}(\ref{itm:semibranch-independence i}). 

By \cref{c:quantifier free type tree}(\ref{itm:quantifier free type tree 3}) and quantifier elimination of $\Tt_{\Gamma}$ and $\Tt'_\gamma$, we have that $\tp_\Gamma(x)=\tp_\Gamma(y)$ in $\Tt_\Gamma$ if and only if $\tp_\gamma(x)=\tp_{\gamma}(y)$ in $\Tt'_\gamma$. Thus, $\indepe[\Gamma]{}$ also inherits invariance and stationarity. 

If $A\indepe[\Gamma]{C}B$ with $A,B,C\subseteq \Tt$, then $A\indepe[\gamma]{C}B$. Thus, by the independent generation property of $\indepe[\gamma]{}$ and \cref{l:generated}(\ref{itm:generated 7}), we conclude that 
%
%
$\langle ACB\rangle_\Gamma=\langle ABC\rangle_\gamma\setminus\{\gamma\}=\langle AC\rangle_\gamma\cup\langle BC\rangle_\gamma\setminus \{\gamma\}=\langle AC\rangle_\Gamma\cup\langle BC\rangle_\Gamma$. As $A,B,C$ are arbitrary, we conclude that $\indepe[\Gamma]{}$ satisfies the independent generation property.

Finally, it remains to show left extension. Suppose $a\indepe[\Gamma]{C}B$ and $a\subseteq a'$ finite tuples in $\Tt$ and $C,B\subseteq \Tt$ finite. Then, by left extension of $\indepe[\gamma]{}$ in $\Tt'_\gamma$, there is $a''$ tuple in $\Tt'$ with $a''\equiv^{\gamma}_{aC}a'$ such that $a''\indepe[\gamma]{C}B$. Now, as $\Tt_\Gamma$ is universal for finite semibranched meet{\hyp}trees, there is $a'''$ tuple in $\Tt_{\Gamma}$ such that $\qftp_{\Gamma'}(a'''/aBC)=\qftp_{\Gamma'}(a''/aBC)$. Hence, by \cref{c:quantifier free type tree}(\ref{itm:quantifier free type tree 3}) and quantifier elimination, $a'''\equiv^\gamma_{aBC}a''$. By invariance, we get that $a'''\indepe[\gamma]{C}B$, so $a'''\indepe[\Gamma]{C}B$. Also, $a'''\equiv^\gamma_{aC}a'$, so $a'''\equiv^{\Gamma}_{aC}a'$ by \cref{c:quantifier free type tree}(\ref{itm:quantifier free type tree 3}) and quantifier elimination. 
\end{proof}
\end{coro}

\section{Generic meet{\hyp}tree expansions} \label{s:section 3}
Let $\lang_0,\lang_1,\lang_2$ be first order languages such that $\lang_0 = \lang_1 \cap \lang_2$ and let $\lang\coloneqq\lang_1\cup \lang_2$. 
%
%
Let $\mathcal{K}_1$ and $\mathcal{K}_2$ be two non{\hyp}empty classes of structures in their respective languages $\lang_1$ and $\lang_2$. The \emph{mix}\footnote{We use the term ``mix'' following \cite{kaplan2019automorphism}. Alternative, one may prefer to call it ``fusion'', following  \cite{kruckman2020interpolative,kruckman2022interpolative} or ``superposition'' following \cite{bodirsky2014ramsey}} of $\mathcal{K}_1$ and $\mathcal{K}_2$ is the class $\mathcal{K}_1\owedge\mathcal{K}_2$ of all 
structures $M$ such that $M_{\mid\lang_i}\in\mathcal{K}_i$ for $i\in\{1,2\}$.

\begin{lem}\label{l:compatible fraisse} In the context of the previous paragraph, suppose $\lang_2\setminus\lang_1$ is relational and $\lang_0=\lang_1\cap\lang_2$ consisting only of constant symbols. Assume $\mathcal{K}_1$ and $\mathcal{K}_2$ are strong Fra\"{\i}ss\'{e} classes of finite structures with infinite Fra\"{\i}ss\'{e} limits $M_1$ and $M_2$ respectively, and assume that all constant terms are in $\lang_0$, i.e. $\langle \emptyset\rangle^{M_1}_1=\langle \emptyset\rangle^{M_1}_0\cong \langle \emptyset\rangle^{M_2}_0=\langle \emptyset\rangle^{M_2}_2$. Then, 
\begin{enumerate}[label={\rm{(\arabic*)}}, wide]
\item $\mathcal{K}_1\owedge\mathcal{K}_2$ is a strong Fra\"{\i}ss\'{e} class with Fra\"{\i}ss\'{e} limit $M$ and
\item $M_{\mid\lang_1}\cong M_1$ and $M_{\mid\lang_2}\cong M_2$.
\end{enumerate}
\begin{proof} 
\begin{enumerate}[label={\rm{(\arabic*)}}, wide]
\item[\hspace{-1.6em}\setcounter{enumi}{1}\theenumi] Given any structure $A_1\in\mathcal{K}_1$ and any structure $A_2\in\mathcal{K}_2$ with $|A_1|=|A_2|$, there is an $\lang${\hyp}structure $A$ with $A_{\mid \lang_1}\cong A_1$ and $A_{\mid\lang_2}\cong A_2$. In particular, $\mathcal{K}_1\owedge\mathcal{K}_2$ is non{\hyp}empty. Furthermore, this also shows that for any $A_1\in\mathcal{K}_1$ and $A_2\in\mathcal{K}_2$, there is $A\in\mathcal{K}_1\owedge\mathcal{K}_2$ such that $A_1\subseteq A_{\mid\lang_1}$ and $A_2\subseteq A_{\mid\lang_2}$. 

The hereditary property {\rm{(HP)}} is obvious as $\mathcal{K}_1$ and $\mathcal{K}_2$ consist of finite structures. We prove the strong amalgamation property {\rm{(SAP)}}. Let $f:\ C\rightarrow A$, $g:\ C\rightarrow B$ embeddings with $A,B,C\in\mathcal{K}_1\owedge\mathcal{K}_2$. By the strong amalgamation properties in $\mathcal{K}_1$ and $\mathcal{K}_2$, there are $\lang_i${\hyp}embeddings $f'_i:\ A_{\mid\lang_i}\rightarrow D_i$ and $g'_i:\ B_{\mid\lang_i}\rightarrow D_i$ with $D_i\in\mathcal{K}_i$ and $f'_i\circ f=g'_i\circ g$ such that $f'_i(A)\cap g'_i(B)=f'_i\circ f(C)$ for $i\in\{1,2\}$. Since $\lang_2$ is relational with constants, we may assume that $|D_2|\leq |D_1|$. As $M_2$ is infinite, we can further assume that $|D_2|=|D_1|$. By the hypotheses, the partial injective maps $f'_2\circ{f'}^{-1}_1$ and $g'_2\circ{g'}_1^{-1}$ agree in their common domain. Thus, there is a bijection $h$ between $D_1$ and $D_2$ extending $f'_2\circ{f'}^{-1}_1$ and ${g'}_2\circ{g'}_1^{-1}$. Take the $\lang${\hyp}expansion of $D_1$ given by making $h$ an $\lang_2${\hyp}isomorphism. Then, $f'\coloneqq f'_1=h^{-1}\circ f'_2:\ A\rightarrow D$ and $g'\coloneqq g'_1=h^{-1}\circ g'_2:\ B\rightarrow D$ are $\lang_1$ and $\lang_2$ embeddings, so $\lang${\hyp}embeddings.

For the strong joint embedding property {\rm{(SJEP)}}, if $\lang_0=\emptyset$, the proof is very similar. Otherwise, this is a special case of SAP because, given $A,B\in \mathcal{K}_1\owedge\mathcal{K}_2$, the assumption on the constant terms implies that $\langle\emptyset\rangle^A\cong\langle\emptyset\rangle^B$ as $\lang${\hyp}structures.

\item This reduces to showing that, for $i\in \{1,2\}$, any $A,B\in\mathcal{K}_i$, $A'\in\mathcal{K}_1\owedge\mathcal{K}_2$ and $\lang_i${\hyp}embeddings $f:\ A\rightarrow A'_{\mid\lang_i}$ and $g:\ A\rightarrow B$, there are an $\lang_i${\hyp}embedding $g':\ B\rightarrow B'_{\mid\lang_i}$ and an $\lang${\hyp}embedding $f':\ A'\rightarrow B'$ with $B'\in\mathcal{K}_1\owedge\mathcal{K}_2$ such that $f'\circ f=g'\circ g$ (because this shows that the restrictions have ``weak ultrahomogeneity'' in the sense of \cite[Lemma 7.1.4]{hodges1993model}). As $\lang_2 \setminus \lang_1$ is relational, the proof for $i=1$ is slightly easier than the proof for $i=2$, so we write only the case for $i=2$.


By the amalgamation property of $\mathcal{K}_2$, we find $B'_0\in\mathcal{K}_2$ and $\lang_2${\hyp}embeddings $f'_0:\ {A'}_{\mid\lang_2}\rightarrow B'_0$ and $g'_0:\ B\rightarrow B'_0$ such that $f'_0\circ f=g'_0\circ g$. Since $M_1$ is infinite, there is ${A'}_{\mid\lang_1}\subseteq B'_1$ with $B'_1\in\mathcal{K}_1$ and $|B'_1|\geq |B'_0|$. As $\lang_2\setminus\lang_1$ is relational, there is $B'_0\subseteq B'_2$ with $B'_2\in\mathcal{K}_2$ and $|B'_1|=|B'_2|$. Take a bijection $h:\ B'_1\rightarrow B'_2$ extending $f'_0$. As $\lang_0$ only contains constants, $h$ is an $\lang_0${\hyp}isomorphism. Let $B'\in\mathcal{K}_1\owedge\mathcal{K}_2$ be the $\lang${\hyp}expansion of $B'_2$ making $\alpha:\ B'_1\rightarrow {B'}_{\mid\lang_1}$ an $\lang_1${\hyp}isomorphism. Then, $g'\coloneqq g'_0:\ B\rightarrow {B'}_{\mid\lang_2}$ is an $\lang_2${\hyp}embedding. On the other hand, $f'\coloneqq f'_0=h_{\mid A'}:\ A'\rightarrow B'$ is an $\lang_2${\hyp}embedding and, also, an $\lang_1${\hyp}embedding, so $f'$ is an $\lang${\hyp}embedding. \qedhere
\end{enumerate}
\end{proof}
\end{lem}
\begin{obs} Since $\lang_2\setminus\lang_1$ is relational, we get that $\langle\tbullet \rangle=\langle\tbullet\rangle_1$ in $M$. Thus, $M$ is uniformly locally finite if $M_1$ is.
\end{obs}

\begin{defi} Let $M$ be an $\lang${\hyp}structure. We say that $M$ is a \emph{generic mix}\footnote{This notion is very close to the notion of interpolative fusion defined in \cite{kruckman2020interpolative,kruckman2022interpolative}. } over $\lang_1,\lang_2$ if for any tuples $a_1,a_2$ of $M$ such that $a_1$ enumerates a finitely generated $\lang_1${\hyp}structure, $a_2$ enumerates a finitely generated $\lang_2${\hyp}structure and $\qftp_0(a_1)=\qftp_0(a_2)$, there is $a$ tuple in $M$ such that $\qftp_1(a)=\qftp_1(a_1)$ and $\qftp_2(a)=\qftp_2(a_2)$. In other words, $M$ is a generic mix over $\lang_1,\lang_2$ if $\Age(M_{\mid\lang_1})\owedge\Age(M_{\mid\lang_2})\subseteq \Age(M)$. 
\end{defi}


The following lemma says that the generic mix property is invariant under adding finite substructures as parameters:
\begin{lem}\label{l:generic mix over parameters} Let $M$ be ultrahomogenous generic mix over $\lang_1,\lang_2$. Let $B$ be a finitely generated $\lang${\hyp}substructure of $M$. Then, the $\lang(B)${\hyp}expansion of $M$ is a generic mix over $\lang_1(B),\lang_2(B)$.
\begin{proof} Set an enumeration $b$ of $B$. Let $a_1$ and $a_2$ be tuples in $M$ such that $a_1$ enumerates a finitely generated $\lang_1(B)${\hyp}structure, $a_2$ enumerates a finitely generated $\lang_2(B)${\hyp}structure and $\qftp_0(a_1/B)=\qftp_0(a_2/B)$. Since $M$ is a generic mix, there are $a'b'$ such that $\qftp_1(a_1b)=\qftp_1(a'b')$ and $\qftp_2(a_2b)=\qftp_2(a'b')$. In particular, $\qftp(b)=\qftp(b')$. So, by ultrahomogeneity, we conclude that there is $a$ such that $ab\equiv a'b'$. Hence, $\qftp_1(a/B)=\qftp_1(a_1/B)$ and $\qftp_2(a/B)=\qftp_2(a_2/B)$, as required.
\end{proof}
\end{lem}

Recall that $\lang_\wedge$ is the language of meet{\hyp}trees, $\lang_{\Gamma}$ is the language of semibranched meet{\hyp}trees and $\lang_\gamma$ is the language of pointed meet{\hyp}trees. Let $\lang_{\Rr}$ be a first{\hyp}order language disjoint to $\lang_\Gamma$ and $\lang_\gamma$. Write $\lang\coloneqq \lang_{\wedge}\cup\lang_{\Rr}$, $\lang_{\Rr\gamma}\coloneqq \lang_{\gamma}\cup\lang_\Rr$ and $\lang_{\Rr\Gamma}\coloneqq \lang_{\Gamma}\cup\lang_\Rr$. Let $M$ be a locally finite Fra\"{\i}ss\'{e} limit over $\lang_{\Rr}$. 

\begin{defi}\label{d:generic meet-tree expansion} We say that $M$ \emph{has a generic meet{\hyp}tree expansion} if the mix of the age of $M$ with the age of $\Tt$ (i.e. the class of finite meet{\hyp}trees) is a Fra\"{\i}ss\'{e} class whose limit is an expansion of $M$ and $\Tt$. The \emph{generic meet{\hyp}tree expansion} of $M$ is the corresponding Fra\"{\i}ss\'{e} limit. Similarly, we define the \emph{generic branched meet{\hyp}tree expansion}.

Let $p\in M$ be a point and $M_p$ be the $(\lang_\Rr\cup\{\gamma\})${\hyp}expansion of $M$ given by interpreting $\gamma$ as $p$. We say that $M$ \emph{has a generic pointed meet{\hyp}tree expansion at $p$} if the mix of the age of $M_p$ with the age of $\Tt_\gamma$ (i.e. the class of pointed finite meet{\hyp}trees) is a Fra\"{\i}ss\'{e} class whose limit is an expansion of $M_p$ and $\Tt_\gamma$ (note that here $\lang_0=\{\gamma\}$). The \emph{generic pointed meet{\hyp}tree expansion at $p$} of $M$ is the corresponding Fra\"{\i}ss\'{e} limit. We say that $M$ \emph{has generic pointed meet{\hyp}tree expansions} if it has a generic pointed meet{\hyp}tree expansion at $p$ for every $p\in M$. Abusing notation, we denote $\gamma\coloneqq p$.
\end{defi}

\begin{lem}\label{l:generic meet-tree expansions} Assume $M$ is a strong Fra\"{\i}ss\'{e} limit and $\lang_\Rr$ is relational. Then:
\begin{enumerate}[label={\rm{(\arabic*)}}, ref={\rm{\arabic*}}, wide]
\item \label{itm:generic meet-tree expansions 1} $M$ has a generic meet{\hyp}tree expansion; denote it by $\Tt^{M}$.
\item \label{itm:generic meet-tree expansions 2} $M$ has generic pointed meet{\hyp}tree expansions and they are isomorphic to the expansions of $\Tt^{M}$ given by adding a constant for a point $\gamma$; denote it by $\Tt^M_\gamma$.
\item \label{itm:generic meet-tree expansions 3} $M$ has a generic branched meet{\hyp}tree expansion and it is isomorphic to the expansion of $\Tt^{M}$ given by adding a predicate for a branch $\Gamma$ and a function for the corresponding branch{\hyp}projection; denote it by $\Tt^{M}_\Gamma$. In particular, $\Tt^M_\Gamma$ is independent up to isomorphism of the choice of the branch\footnote{Note that, contrary to the branch case, $\Tt^M_\gamma$ is not in general independent of the choice of $\gamma$ unless there is a unique $1${\hyp}type in $M$.}.
\end{enumerate}
\begin{proof}
\begin{enumerate}[label={\rm{(\arabic*)}}, wide]
\item[\hspace{-1.6em}\setcounter{enumi}{1}\theenumi] It is straightforward to check that the class of branched/pointed meet{\hyp}trees is a strong amalgamation class. The fact that $M$ has generic meet{\hyp}tree, generic branched meet{\hyp}tree and generic pointed meet{\hyp}tree expansions follows from \cref{l:compatible fraisse}. This gives (1) and also the first parts of (2) and (3). 
\item Let $\gamma\in\Tt^M$ which we consider as an expansion of $M$, $(\Tt^M)'$ be the generic pointed meet{\hyp}tree expansion at $\gamma$ and $\Tt^M_\gamma$ the $\lang\cup \{\gamma\}${\hyp}expansion of $\Tt^M$ adding a constant for $\gamma$. We need to show that $\Tt^M_\gamma\cong (\Tt^M)'$. Since both are ultrahomogeneous, it is enough to show that both have the same age. Obviously, $(\Tt^M_\gamma)_{\mid\lang_\Rr\cup \{\gamma\}}= M_\gamma$ and $(\Tt^M_\gamma)_{\mid \lang_\gamma}\cong \Tt_\gamma$. By \cref{l:generic mix over parameters}, $\Tt^M_\gamma$ is a generic mix over $\lang_\Rr\cup\{\gamma\}$ and $\lang_\gamma$, so $\Age(M_\gamma)\owedge\Age(\Tt_\gamma)\subseteq\Age(\Tt^M_\gamma)$. Furthermore, as $\Tt^M_\gamma$ is locally finite, it follows that $\Age((\Tt^M)')=\Age(M_\gamma)\owedge\Age(\Tt_\gamma)=\Age(\Tt^M_\gamma)$.
\item Let $\Tt^M_\Gamma$ be the $\lang_{\Rr\Gamma}${\hyp}expansion of $\Tt^M$ given by adding a predicate for a branch and a function for the branch projection. It suffices to show that $\Tt^M_\Gamma$ is ultrahomogeneous and $\Age(\Tt^M_\Gamma)=\Age(M)\owedge\Age(\Tt_\Gamma)$. Note that $\Tt^M_\Gamma$ is locally finite as $\Tt_\Gamma$ is, so $\Age(\Tt^M_\Gamma)\subseteq\Age(M)\owedge\Age(\Tt_\Gamma)$.

We prove the following statement:
\begin{enumerate}[label={\rm{($\ast$)}}, wide]
\item Let $A\subseteq \Tt^M_\Gamma$ be a finite substructure (or empty) and $f:\ A\rightarrow B$ be an $\lang_{\Rr\Gamma}${\hyp}embedding with  $B\in\Age(M)\owedge\Age(\Tt_\Gamma)$. Then, there is some $g:\ B\rightarrow \Tt^M_\Gamma$ such that $g \circ f=\id_A$.
\end{enumerate}
Note that when $A=\emptyset$ this gives $\Age(\Tt^M_\Gamma)=\Age(M)\owedge\Age(\Tt_\Gamma)$, and for $A\neq\emptyset$ it gives ultrahomogeneity. 

First we reduce to the case that 
\begin{enumerate}[label={\rm{($\ast \ast$)}}, wide]
\item $A$ is non empty, the semibranches of $A$ and $B$ are branches and the maximum of the branch in $A$ is mapped to the the maximum of the branch in $B$. 
\end{enumerate}

Let $A'=\langle A\,\gamma\rangle_\Gamma$ with $\gamma\in\Gamma\subseteq \Tt^M$ and $\gamma>\max A\cap \Gamma$ and let $B'_1=\langle B\,\gamma'\rangle_\Gamma$ be a semibranched meet{\hyp}tree extending $B_{\mid\lang_\Gamma}$ with $\gamma'$ greater than the maximum of the semibranch in $B$. Note that by \cref{l:generated}(\ref{itm:generated 7}) $A'=A\cup\{\gamma\}$ and $B'_1=B\cup\{\gamma'\}$. Consider $f':\ A'_{\mid\lang_\Gamma}\rightarrow B'_1$ extending $f$ by mapping $\gamma$ to $\gamma'$. It is easy to see that $f'$ is an embedding of semibranched meet{\hyp}trees. Take the $\lang_{\Rr\Gamma}${\hyp}expansion $f'(A')$ of $f'(A'_{\mid\lang_\Gamma})$ making $f':\ A'\rightarrow f'(A')$ an $\lang_{\Rr\Gamma}${\hyp}isomorphism. Let $B'_2$ be a strong amalgam of ${f'(A')}_{\mid\lang_\Rr}$ and $B_{\mid\lang_\Rr}$ over ${f(A)}_{\mid \lang_\Rr}$. As $\lang_\Rr$ is relational, we may assume that $B'_2$ is the union of $f'(A')$ and $B'_1$ as sets, i.e. its universe is the same as $B'_1$. Let $B'$ be the unique common $\lang_{\Rr\Gamma}${\hyp}expansion of $B'_1$ and $B'_2$. Then, $f':\ A'\rightarrow B'$ is an $\lang_{\Rr\Gamma}${\hyp}embedding extending $f$ with $B'\in\Age(M)\owedge\Age(\Tt_\Gamma)$, $A'\neq\emptyset$, the semibranches of $A'$ and $B'$ are branches and the maximum of the branch in $A'$ is mapped to the the maximum of the branch in $B'$. 

So, assume ($\ast\ast$). By (1), there is $g:\ {B}_{\mid\lang}\rightarrow \Tt^M$ such that $g\circ f=\id_{A_{\mid\lang}}$. In particular, $g(\gamma_B)=\gamma_A\in \Gamma$ where $\gamma_B$ is the maximum of the branch in $B$ and $\gamma_A$ is the maximum of the branch in $A$. Thus, the image of the branch of $B$ is precisely $g(B)\cap \Gamma$. Furthermore, for any $b\in B$, $\pi(b)=b\meet\gamma_B$, so $g(\pi(b))=g(b)\meet\gamma_A=\pi(g(b))$. Therefore, $g:\ B\rightarrow \Tt^M_\Gamma$ is also an $\lang_{\Rr\Gamma}${\hyp}embedding.  
\qedhere
\end{enumerate}
\end{proof}
\end{lem}
\begin{nota} One might find more natural at first sight to denote the generic meet{\hyp}tree expansion of $M$ by $M^{\Tt}$ (or some variation like $M^{\mathrm{tree}}$) rather than $\Tt^M$. This notation would suggest that $M$ is considered the main structure while the tree is the extra one. For this paper, however, it is more convenient to consider the tree as the main structure and $M$ as the extra one, so we prefer the notation $\Tt^M$. Nevertheless, it is important to note that, contrary to what these notations suggest, the operation is symmetric.
\end{nota}

\begin{ejem}
\begin{enumerate}[label={\rm{(\arabic*)}}, wide] 
\item[\hspace{-1.4em}\setcounter{enumi}{1}\theenumi] The universal dense (branched, pointed) meet{\hyp}tree is the generic (branched, pointed) meet{\hyp}tree expansion of the countable infinite set over the empty language (i.e. just equality).
\item The \emph{Rado (branched, pointed) meet{\hyp}tree} is the generic (branched, pointed) meet{\hyp}tree expansion of the Rado graph (i.e. the random graph). We denote the Rado meet{\hyp}tree by $\Tt^\Rr$.
\item The \emph{generically ordered universal dense (branched, pointed) meet{\hyp}tree} is the generic (branched, pointed) meet{\hyp}tree expansion go the countable dense linear order. We denote the generically ordered universal dense meet{\hyp}tree by $\Tt^{<}$.
\end{enumerate}
\end{ejem}



In the following lemma we show how to mix stationary weak independence relations (see \cref{d:stationary weak independence relation}). Recall that that the independent generation property of an independence relation $\indepe{}$ says that if $A\indepe{C}B$, then $\langle ABC\rangle=\langle AC\rangle\cup\langle BC\rangle$.

\begin{lem}\label{l:mix independence} Let $M$ be a locally finite ultrahomogeneous generic mix over $\lang_1,\lang_2$. Let $\indepe[1]{}$ and $\indepe[2]{}$ be stationary weak independence relations on finite subsets satisfying the independent generation property on $M_{\mid\lang_1}$ and $M_{\mid\lang_2}$ respectively. Let $\indepe{}$ be the ternary relation on finite subsets defined by $A\indepe{C}B$ if and only if $\langle AC\rangle\indepe[1]{\langle C\rangle}\langle BC\rangle$ and $\langle AC\rangle\indepe[2]{\langle C\rangle}\langle BC\rangle$. Assume that $\indepe{}$ is stationary in $M_{\mid\lang_0}$. Then, $\indepe{}$ is a stationary weak independence relation satisfying the independent generation property on $M$.
\begin{proof} Invariance, monotonicity, normality, transitivity, existence and the independent generation property are obvious. 
\begin{enumerate}[ , wide]
\item[\rm{(Base monotonicity)}] Suppose $A\indepe{C}B$ and $C\subseteq C'\subseteq B$. First of all, note that by right monotonicity and independent generation, we get that $\langle AC'\rangle=\langle AC\rangle\cup \langle C'\rangle$. Then, $\langle AC\rangle\indepe[1]{\langle C\rangle}\langle B\rangle$ and $\langle AC\rangle\indepe[2]{\langle C\rangle}\langle B\rangle$ where $\langle C\rangle\subseteq \langle C'\rangle\subseteq \langle B\rangle$. Thus, by right base monotonicity and left normality of $\indepe[1]{}$ and $\indepe[2]{}$, we get $\langle AC\rangle\cup \langle C'\rangle\indepe[1]{\langle C'\rangle}\langle B\rangle$ and $\langle AC\rangle\cup\langle C'\rangle\indepe[1]{\langle C'\rangle}\langle B\rangle$. Hence, we have concluded that $\langle AC'\rangle\indepe[1]{\langle C'\rangle}\langle B\rangle$ and $\langle AC'\rangle\indepe[2]{\langle C'\rangle}\langle B\rangle$, i.e. $A\indepe{C'}B$, so $\indepe{}$ satisfies right base monotonicity. Left base monotonicity is proved symmetrically.
\item[\rm{(Extension)}] Take $a,B,C$ with $a\indepe{C}B$ and $a\subseteq a'$. We may assume without loss of generality that $B=\langle BC\rangle$, $C=\langle C\rangle$, and $a$ and $a'$ enumerate substructures containing $C$. By extension of $\indepe[1]{}$ and $\indepe[2]{}$, there are $a'_1$ and $a'_2$ such that $a'_1\indepe[1]{C}B$, $a'_2\indepe[2]{C}B$, $a'_1\equiv^{1}_Ca'$ and $a'_2\equiv^{2}_Ca'$. Since $a'_1 \equiv^{0}_Ca' \equiv^{0}_C a'_2$, by stationarity of $\indepe{}$ in $M_{\mid\lang_0}$, we get that $a'_1\equiv^{0}_{B}a'_2$. As $M$ is a generic mix over $\lang_1,\lang_2$, by \cref{l:generic mix over parameters}, $M$ is a generic mix over $\lang_1(B),\lang_2(B)$, so there is $a''$ such that $a''\equiv^{1}_{B}a'_1$ and $a''\equiv^{2}_{B}a'_2$. In particular, by invariance of $\indepe[1]{}$  and $\indepe[2]{}$, we have $a''\indepe{C}B$ and $a''\equiv^{1}_Ca'$ and $a''\equiv^{2}_Ca'$. Since $M$ is ultrahomogeneous, we conclude that $a''\equiv_Ca'$. In sum, $a''\equiv_Ca'$ and $a''\indepe{C}B$, so $\indepe{}$ satisfies left extension. Right extension follows from \cref{o:left extension implies right extension}.
\item[\rm{(Stationarity)}] Suppose $a\indepe{C}B$, $a'\indepe{C}B$ and $a\equiv_Ca'$. We may assume without loss of generality that $B=\langle BC\rangle$ and $C=\langle C\rangle$, and $a$ and $a'$ enumerate substructures $A$ and $A'$ containing $C$. Then, $a\indepe[1]{B}B$, $a'\indepe[1]{C}B$, $a\indepe[2]{C}B$, $a'\indepe[2]{C}B$, $a\equiv^{1}_Ca'$, $a\equiv^{2}_Ca'$. Then, by stationarity of $\indepe[1]{}$ and $\indepe[2]{}$, we conclude that $a\equiv^{1}_{B}a'$ and $a\equiv^{2}_{B}a'$. Now, by the independent generation property, $\langle AB\rangle=A\cup B$ and $\langle A' B\rangle=A'\cup B$. As $M$ is ultrahomogeneous, we conclude that $a\equiv_{B}a'$, so $\indepe{}$ satisfies stationarity.
\qedhere 
\end{enumerate}
\end{proof}
\end{lem}

\begin{coro} Let $M$ be an infinite strong Fra\"{\i}ss\'{e} limit over a relational language $\lang_\Rr$ disjoint to $\lang_{\Gamma}$ and $\lang_\gamma$. Assume that $M$ has a stationary weak independence relation. Let $\Tt^M$ be the generic meet{\hyp}tree expansion of $M$, $\Gamma$ a branch and $\gamma$ a point. Then, $\Tt^M_\gamma$ and $\Tt^{M}_\Gamma$ have stationary weak independence relations satisfying reflexivity and the independent generation property.
\end{coro}

We also recall the following properties (see \cite{conant2016axiomatic}):
\begin{defi} \label{d:free weak independence} A \emph{free weak independence relation} $\indepe{}$ is a stationary weak independence relation satisfying the independent generation property, reflexivity and freeness, where:
\begin{enumerate}[wide]
\item[\rm{(Reflexivity)}] If $a\indepe{C}a$, then $a\in \langle C\rangle$.
\item[\rm{(Freeness)}] If $A\indepe{C}B$ and $\langle C\rangle \cap\langle AB\rangle\subseteq \langle D\rangle$ where $D\subseteq C$, then $A\indepe{D}B$.
\end{enumerate}
\end{defi}
\begin{obs} Typically, reflexivity, independent generation and freeness are stated with $\acl$ in place of $\langle\tbullet\rangle$ --- and independent generation is called independent closure. Here we take this variant noting that, for the structures studied in this paper, $\acl(C)=\langle C\rangle$. 
\end{obs} 
\begin{ejem} The countable universal homogeneous tournament (i.e. the Fra\"{\i}ss\'{e} limit of asymmetric complete binary relations) has the free weak independence relation given by $A\indepe{C}B$ if every element of $A\setminus C$ wins against every element of $B\setminus C$.
\end{ejem}

The following is a slight generalisation of \cite[Lemma 3.1]{li2020automorphism}.
\begin{lem} \label{l:mix free independence} Let $M$, $\indepe[1]{}$, $\indepe[2]{}$ and $\indepe{}$ be as in \cref{l:mix independence}. Assume that $\indepe[1]{}$ satisfies reflexivity, $\indepe[2]{}$ satisfies freeness and $\lang_2\setminus\lang_1$ is relational. Let $A,B,C,D$ be finite subsets of $M$ with $\langle AC\rangle\cap\langle BC\rangle=\langle C\rangle$ such that $A\indepe{BC}D$ and $B\indepe[1]{C}D$. Then, $A\indepe{C}D$.
\begin{proof} We need to show that $\langle AC\rangle\indepe[1]{\langle C\rangle}\langle DC\rangle$ and $\langle AC\rangle\indepe[2]{\langle C\rangle}\langle DC\rangle$. Since $A\indepe{BC}D$, in particular, $\langle ABC\rangle\indepe[1]{\langle BC\rangle}\langle DBC\rangle$. Since $B\indepe[1]{C}D$, by normality, base monotonicity and extension of $\indepe[1]{}$, we have $\langle BC\rangle_1\indepe[1]{\langle C\rangle_1}\langle DC\rangle_1$. Since $\lang_2\setminus\lang_1$ is relational, $\langle\tbullet\rangle_1=\langle\tbullet \rangle$. By transitivity and monotonicity of $\indepe[1]{}$, we have that $\langle AC\rangle \indepe[1]{\langle C\rangle}\langle DC\rangle$. By the independent generation property of $\indepe[1]{}$, we have $\langle ADC\rangle=\langle AC\rangle\cup\langle DC\rangle$.

On the other hand, since $A\indepe{BC}D$, in particular, $\langle AC\rangle \indepe[2]{\langle BC\rangle }\langle DC\rangle$. Now, $\langle BC\rangle\cap \langle ADC\rangle= (\langle BC\rangle\cap \langle AC\rangle)\cup(\langle BC\rangle\cap \langle DC\rangle)$. By hypothesis $\langle BC\rangle\cap\langle AC\rangle=\langle C\rangle$. Also, $\langle BC\rangle\indepe[1]{\langle C\rangle}\langle DC\rangle$, so $\langle BC\rangle\cap\langle DC\rangle=\langle C\rangle$ by reflexivity of $\indepe[1]{}$. Hence, $\langle BC\rangle\cap \langle ADC\rangle= \langle C\rangle$, concluding $\langle AC\rangle\indepe[2]{\langle C\rangle}\langle DC\rangle$ by freeness of $\indepe[2]{}$. 
\end{proof} 
\end{lem}



\section{Back{\hyp}and{\hyp}forth arguments} \label{s:section 4}
Through this section, unless otherwise stated, we work in the generic meet{\hyp}tree expansion $\Tt^M$ of an infinite strong Fra\"{\i}ss\'{e} limit $M$ over a finite relational language $\lang_\Rr$ disjoint to $\lang_{\Gamma}$ and $\lang_\gamma$. As prototypical cases, one may consider the generic meet{\hyp}tree expansion of the (countable) dense linear order, which we denote by $\Tt^{<}$, or the generic meet{\hyp}tree expansion of the random graph, which we denote by $\Tt^\Rr$. \medskip

We start this section with the key result (\cref{l:dicotomia}), which is an improvement of \cite[Proposition 4.19]{kaplan2019automorphism}. 
\begin{defi} \label{d:fan} Let $T$ be a meet{\hyp}tree, $\gamma$ a point and $\alpha$ an automorphism of $T$. We say that $\alpha$ is a \emph{fan above $\gamma$} if $a\wedge \alpha(a)=\gamma$ for any $a\geq \gamma$. We say that $\alpha$ is an \emph{infinite fan above $\gamma$} if $\alpha^n$ is a \emph{fan above $\gamma$} for any $n\in\N_{>0}$. We say that $\alpha$ is an \emph{(infinite) fan} if it is so at some point $\gamma$.
\end{defi}
\begin{obs} If $\alpha$ is a fan above $\gamma$, then $\alpha(\gamma)=\gamma$. Indeed, $\alpha(\gamma)\geq \gamma\wedge\alpha(\gamma)=\gamma$, so $\alpha(\gamma)=\alpha(\gamma\wedge\alpha(\gamma))=\alpha(\gamma)\wedge\alpha^2(\gamma)=\gamma$. Consequently, $\alpha$ is a fan above $\gamma$ if and only if $\alpha$ fixes $\gamma$ and there is no $\alpha${\hyp}invariant cone above $\gamma$. Similarly, $\alpha$ is an infinite fan above $\gamma$ if and only $\alpha$ fixes $\gamma$ and there is no $\alpha${\hyp}invariant finite set of cones above $\gamma$.
\end{obs}
\begin{lem} \label{l:dicotomia} Let $T$ be a meet{\hyp}tree and $\alpha$ an automorphism of $T$. Then, either $\alpha$ setwise fixes a branch or $\alpha$ is a fan.
\begin{proof} By Zorn's Lemma, there is a maximal setwise fixed chain $\Gamma$. 
\begin{enumerate}[wide]
\item[($\ast$)] Note that if $a\leq b$ for some $b\in \Gamma$, then $a \in \Gamma$. Indeed, $a\in \Gamma'\coloneqq \bigcup_{b\in\Gamma} T_{\leq b}$ and $\Gamma'$ is a setwise fixed chain containing $\Gamma$, so we have that $a\in \Gamma'=\Gamma$.
\item[($\ast\ast$)] Now, note that if $a\geq b$ for every $b\in\Gamma$ then $c\coloneqq a \meet \alpha(a)=\max\Gamma$. Indeed, $\alpha(a)\geq b$ for every $b\in\Gamma$ too, so $c\geq b$ for any $b\in\Gamma$. Since $c\leq \alpha(a)$ and $\alpha(c)=\alpha(a)\wedge \alpha^2(a)\leq\alpha(a)$, we have that $c$ and $\alpha(c)$ are comparable and above all the elements of $\Gamma$. Thus, $\Gamma''\coloneqq\Gamma\cup\{\alpha^i(c)\tq i\in\Z\}$ is a setwise fixed chain containing $\Gamma$, so $c\in \Gamma''=\Gamma$. Since $c\geq b$ for all $b\in\Gamma$, we have that $c=\max\Gamma$. 
\end{enumerate}

Finally, if $\Gamma$ is not a branch, then by ($\ast$) there must be some $a$ as in ($\ast\ast$), and thus $\Gamma$ must have a maximal element $c$, and $\alpha$ is a fan above $c$: if $a \geq c$ then ($\ast\ast$) applies and we are done.
\end{proof}
\end{lem}

Thus, in particular, any automorphism $\alpha$ of $\Tt^M$ setwise fixes  a branch or is a fan. Our aim for the rest of the section is to show that, up to a finite product of conjugates, we can make $\alpha$ strictly increasing (inside the corresponding semibranch that it fixes setwise) and, in the fan case, an infinite fan. The full argument is a long sequence of back{\hyp}and{\hyp}forth manipulations very similar to the ones developed in \cite{li2020automorphism}. All the arguments fundamentally rely on the following two basic facts:  


\begin{lem} \label{l:saturation of cones} 
Let $\gamma$ be a point, $a,b$ single elements above $\gamma$ and $\bar{c},\bar{d}$ finite tuples with $\bar{c}\equiv^{\Rr\gamma} \bar{d}$. Suppose $a\notin\cone_\gamma(\bar{c})$ and $b\notin\cone_\gamma(\bar{d})$. Then, there is  $a'\in \cone_\gamma(b)$ such that $a' \bar{d}\equiv^{\Rr\gamma}a \bar{c}$. Furthermore, we can take $a'>b$.
\begin{proof} Without loss of generality, we can assume that $\bar{c}=\bar{d}$ and write $C$ for the set enumerated by $\bar{c}$. By \cref{l:generated}(\ref{itm:generated 8}), assume without loss of generality that $C=\langle C\rangle_\gamma$. Take $a''>b$ arbitrary. By \cref{c:cones are indiscernible}, $a\equiv^\gamma_C a''$. Now, by \cref{l:generated}(\ref{itm:generated 4}), $D\coloneqq\langle b\, C\rangle_{\gamma}=C\cup \{b\}$. Also, by \cref{l:generated}(\ref{itm:generated 4}), $\langle a\, D\rangle_\gamma=D\cup\{a\}$ and $\langle a''D\rangle_{\gamma}=D\cup \{a''\}$. Thus, by the generic mix property over $D$ (see \cref{l:generic mix over parameters}), there is $a'$ such that $a'\equiv^{\gamma}_D a''$ and $a'\equiv^{\Rr}_{D}a$. In particular, $a'\in\cone_\gamma(b)$, $a'\equiv^{\gamma}_C a''\equiv^{\gamma}_C a$ and $a'\equiv^\Rr_C a$. Since $\langle a\, C\rangle_{\gamma}=C\cup \{a\}$ by \cref{l:generated}(\ref{itm:generated 4}), we conclude that $a'\equiv^{\Rr\gamma}_C a$ by ultrahomogeneity.
\end{proof}
\end{lem}

\begin{lem}\label{l:saturation of branch intervals} 
\begin{enumerate}[label={\rm{(\arabic*)}}, ref={\rm{\arabic*}}, wide]
\item[\hspace{-1.4em}\setcounter{enumi}{1}\rm{(\theenumi)}] 
Let $\gamma$ be a point, $a,b_1,b_2<\gamma$ with $b_1<b_2$ and $\bar{c}=(c_i)_{i<n}$, $\bar{d}=(d_i)_{i<n}$ finite tuples with $\bar{c}\equiv^{\Rr\gamma}\bar{d}$. Assume $a>c_i\meet\gamma$ for all $i<n$ such that $b_1\geq d_i\meet \gamma$, and assume $a<c_i\meet \gamma$ for all $i<n$ such that $b_2\leq d_i\meet\gamma$. Then, there is $a'$ with $b_1<a'<b_2$ such that $a' \bar{d}\equiv^{\Rr\gamma}a \bar{c}$.
\item 
Let $\Gamma$ be a branch, $a,b_1,b_2\in\Gamma$ with $b_1<b_2$ and $\bar{c}=(c_i)_{i<n},\bar{d}=(d_i)_{i<n}$ finite tuples with $\bar{c}\equiv^{\Rr\Gamma}\bar{d}$. Assume $a>\pi_\Gamma(c_i)$ for all $i<n$ such that $b_1\geq \pi_\Gamma(d_i)$, and assume $a<\pi_\Gamma(c_i)$ for all $i<n$ such that $b_2\leq \pi_\Gamma(d_i)$. Then, there is $a'$ with $b_1<a'<b_2$ such that $a' \bar{d}\equiv^{\Rr\Gamma}a \bar{c}$.
\end{enumerate}
\begin{proof} Both statements are very similar, so we just prove {\rm{(1)}}. Without loss of generality, we can assume the following:
\begin{enumerate}[label={$\bullet$}, noitemsep, wide]
\item $\bar{c}=\bar{d}$  and let $C$ for the set enumerated by $\bar{c}$.
\item There is $c\in C$ with $c<a,b_1,b_2,\gamma$.
\item $C=\langle C\rangle_{\gamma}$ by \cref{l:generated}(\ref{itm:generated 2}).
\end{enumerate} 
Write $\Gamma\coloneqq\{x\tq x\leq\gamma\}$ and $\pi(C)\coloneqq\{c\meet\gamma\tq c\in C\cup\{\gamma\}\}$. 

First of all, note that, if $a\in C$, then $b_1<a<b_2$ by hypothesis, so the lemma holds trivially with $a'=a$. So let's assume $a\notin C$. Additionally, we can assume that $a\neq b_1$ and $a\neq b_2$. Indeed, if $a=b_i$ with $i\in\{1,2\}$, then $b_i \notin C$ and, by density of the tree and finiteness of $C$, we can replace $b_i$ by $b_i' \neq a$ with $b_1<b_i'<b_2$ such that $\qftp_<(b_i'/\pi(C))=\qftp_<(b_i/\pi(C))$.

Let $c_1=\max\{c\in\pi(C)\tq c<a\}$ and $c_2=\min\{c\in \pi(C)\tq c>a\}$. By hypothesis, $c_1<b_2$ and $c_2>b_1$. Take $d_1=\max\{b_1,c_1\}$ and $d_2=\min\{b_2,c_2\}$, so $d_1<d_2$. Now, for any $a''$ with $d_1<a''<d_2$, we have that $\qftp_{<}(a''/\pi(C))=\qftp_{<}(a/\pi(C))$, so $\qftp_{\Gamma}(a''/C\,\gamma)=\qftp_{\Gamma}(a/C\,\gamma)$ by \cref{c:semibranch intervals are indiscernible}. In particular, we have $\qftp_{\gamma}(a''/C)=\qftp_{\gamma}(a/C)$, so $a''\equiv^{\gamma}_C a$ by quantifier elimination of $\Tt_\gamma$. 

By \cref{l:generated}(\ref{itm:generated 4}), we have that $D\coloneqq\langle d_1d_2 C\rangle_{\gamma}=C\cup \{d_1,d_2\}$. Similarly, by \cref{l:generated}(\ref{itm:generated 4}), we have $\langle a\, D\rangle_\gamma=D\cup\{a\}$ and $\langle a'' D\rangle_{\gamma}=D\cup \{a''\}$ with $a,a''\notin D$. Thus, by the generic mix property over $D$ (see \cref{l:generic mix over parameters}), there is $a'$ such that $a'\equiv^{\gamma}_D a''$ and $a'\equiv^\Rr_D a$. In particular, $d_1<a'<d_2$, $a'\equiv^{\gamma}_C a''\equiv^{\gamma}_C a$ and $a''\equiv^\Rr_C a$. Since $\langle a\, C\rangle_{\gamma}=C\cup \{a\}$ by \cref{l:generated}(\ref{itm:generated 7}), we conclude that $a'\equiv^{\Rr\gamma}_C a$ by ultrahomogeneity.
%
%
\end{proof}
\end{lem}

Recall that, for a semibranch $\Gamma$, we write $\Gamma^\circ\coloneqq \{a\in \Gamma\tq \mathrm{there\ is\ }b\in \Gamma\mathrm{\ with\ }a<b\}$.  By \cref{l:semibranches}, we have that $\Gamma$ is not a branch in $\Tt$ if and only if $\Gamma\neq \Gamma^{\circ}$ in which case $\Gamma\setminus \Gamma^\circ=\{\gamma\}$ with $\gamma=\max\Gamma$.

\begin{defi} \label{d:unbounded}
Let $\Gamma$ be a semibranch and $\alpha$ an automorphism fixing $\Gamma$ setwise.
\begin{enumerate}[label={$\bullet$}, wide]
\item $\alpha$ is \emph{strictly increasing at a point $x$} if $\alpha(x)>x$. Similarly, $\alpha$ is \emph{strictly decreasing at $x$} if $\alpha(x)<x$.
\item A subset $L$ of $\Gamma^{\circ}$ is \emph{upwards unbounded (downwards unbounded)} if for any $x\in \Gamma^{\circ}$ there is $y\in L$ with $y>x$ (with $y<x$); a subset is \emph{unbounded} in $\Gamma^{\circ}$ if it is upwards and downwards unbounded in $\Gamma^{\circ}$.   
\item $\alpha$ is \emph{unbounded (upwards unbounded/downwards unbounded) increasing in $\Gamma^{\circ}$} if the set of points of $\Gamma^{\circ}$ where $\alpha$ is strictly increasing is unbounded (upwards unbounded/downwards unbounded) in $\Gamma^{\circ}$. 
\item $\alpha$ is \emph{unbounded (upwards unbounded/downwards unbounded) decreasing in $\Gamma^{\circ}$} if the set of points of $\Gamma^{\circ}$ where $\alpha$ is strictly decreasing is unbounded (upwards unbounded/downwards unbounded) in $\Gamma^{\circ}$.
\item $\alpha$ is \emph{unbounded (upwards unbounded/downwards unbounded) in $\Gamma^{\circ}$} if  the set of points of $\Gamma^{\circ}$ not fixed by $\alpha$ is unbounded (upwards unbounded/downwards unbounded) in $\Gamma^{\circ}$.
\end{enumerate}

When $\Gamma\neq\Gamma^{\circ}$ and $\gamma\coloneqq\max\Gamma$, we may say \emph{below $\gamma$} in place of \emph{in $\Gamma^{\circ}$}.
\end{defi}

\begin{lem} \label{l:saturation of increasing points} Let $\Gamma$ be a semibranch, $a,b\in\Gamma^\circ$ and $\bar{c}=(c_i)_{i<n}, \bar{d}=(d_i)_{i<n}$ finite tuples with $\bar{c}\equiv^{\Rr\Gamma}\bar{d}$. Let $C_0\coloneqq\{c_i\tq \pi_\Gamma(c_i)\in\Gamma^\circ\}$ and $\alpha$ an automorphism of $\Tt^M$ fixing $\Gamma$ setwise. 
\begin{enumerate}[label={\rm{(\arabic*)}}, ref={\rm{\arabic*}}, wide]
\item \label{itm:saturation of increasing points 1} Assume that $\alpha$ is upwards unboundedly increasing in $\Gamma$ and $a>\pi_{\Gamma}(c)$ for all $c\in C_0$. Then, there is $a'>b$ with $a'\bar{d}\equiv^{\Rr\Gamma}a\bar{c}$ and $\alpha(a')>a'$.
\item \label{itm:saturation of increasing points 2} Assume that $\alpha$ is downwards unboundedly increasing in $\Gamma$ and $a<\pi_\Gamma(c)$ for all $c\in C_0$. Then, there is $a'<b$ with $a'\bar{d}\equiv^{\Rr\Gamma}a\bar{c}$ and $\alpha(a')>a'$.
\end{enumerate}
\begin{proof} \begin{enumerate}[label={\rm{(\arabic*)}}, wide]
\item[\hspace{-1.4em}\setcounter{enumi}{1}\theenumi] Write $D_0\coloneqq\{d_i\tq \pi_{\Gamma}(d_i)\in\Gamma^{\circ}\}$. Pick $e\in\Gamma^\circ$ with $e>b$ and $e>\pi_\Gamma(d)$ for each $d\in D_0$ such that $\alpha(e)>e$. Then, by \cref{l:saturation of branch intervals}, there is $a'\bar{d}\equiv^{\Rr\Gamma}a\bar{c}$ with $e<a'<\alpha(e)$. Hence, $a'>b$ and $\alpha(a')>\alpha(e)>a'$.
\item Write $D_0\coloneqq\{d_i\tq \pi_{\Gamma}(d_i)\in\Gamma^{\circ}\}$. Pick $e\in\Gamma^\circ$ with $e<\alpha^{-1}(b)$ and $e<\alpha^{-1}(\pi_\Gamma(d))$ and $e<\pi_\Gamma(d)$ for each $d\in D_0$ such that $\alpha(e)>e$. Then, by \cref{l:saturation of branch intervals}, there is $a'\bar{d}\equiv^{\Rr\Gamma}a\bar{c}$ with $e<a'<\alpha(e)$. Hence, $a'<b$ and $\alpha(a')>\alpha(e)>a'$. \qedhere
\end{enumerate}
\end{proof}
\end{lem}

Recall that $\sim_\gamma$ denotes the cone equivalence relation above $\gamma$. Note that if $\Gamma$ is a branch and $\gamma$ a point, then if $\gamma \notin \Gamma$ then $\cone_\gamma(a)$ is disjoint to $\Gamma$ for all $a>\gamma$. Otherwise, there is a unique cone intersecting $\Gamma$ (in fact, containing $\Gamma_{>\gamma}$).
\begin{lem} \label{l:trivial obs} Let $\Gamma$ be a branch, $\gamma$ a point and $\mathcal{C}$ the set of all cones above $\gamma$ disjoint to $\Gamma$. Let $\sigma$ be a permutation of $\mathcal{C}$. Then, there is an automorphism $\alpha$ of $\Tt^M_\Gamma$ satisfying the following properties: 
\begin{enumerate}[label={\rm{(\arabic*)}}, wide]
\item The only points fixed by $\alpha$ are $\gamma$ and $\pi_{\Gamma}(\gamma)$.
\item The induced action on $\mathcal{C}$ by $\alpha$ is $\sigma$, i.e. $\alpha(C)=\sigma(C)$ for any $C\in \mathcal{C}$.
\end{enumerate}
\begin{proof} 
First of all, note that: 
\begin{enumerate}[label={($\ast$)}, ref={($\ast$)}, wide]
\item \label{itm:neuman} For any isomorphism $\alpha:\ A\rightarrow B$ between finite substructures of $\Tt^M_\Gamma$ and any element $a$, there is an isomorphism $\alpha':\ A'\rightarrow B'$ between finite substructures of $\Tt^M_\Gamma$ extending $\alpha$ with $a\in A'$ such that all its fixed points are contained in $A$. 
\end{enumerate}

Indeed, by ultrahomogeneity, there is $\alpha'_0:\ A'\rightarrow B'_0$ extending $\alpha$ with $a\in A'$ finite. Now, by strong amalgamation, $\acl$ is trivial, i.e. $\acl(X)=X$ for any set $X$. By Neumann's Lemma \cite[Lemma 4.1]{Neumann}, there is an automorphism $\sigma$ fixing $B$ with $A'\cap\sigma(B'_0)\subseteq B$. Then, $\alpha'\coloneqq\sigma\circ\alpha'_0:\ A'\rightarrow B'\coloneqq \sigma(B'_0)$ extends $\alpha$ and all its fixed points are contained in $A$.

Let $\{w_n\}_{n\in\N}$ be an enumeration of $\Tt^M$ with $w_0=\gamma$. By a back{\hyp}and{\hyp}forth recursion, we define a chain of isomorphisms $\alpha_n:\ A_n\rightarrow B_n$, between finite substructures of $\Tt^M_\Gamma$ with $w_n\in A_n$ and $w_n\in B_n$, satisfying the properties {\rm{(1,2)}} restricted to their domains and images (and {\rm{(2)}} also for $\alpha_n^{-1}$ and $\sigma^{-1}$). 
\begin{enumerate}[wide]
\item[Step $n=0$:] Simply take $\alpha_0:\ \gamma\mapsto \gamma;\ \pi_{\Gamma}(\gamma)\mapsto \pi_{\Gamma}(\gamma)$ and $A_0=B_0=\{\gamma,\pi_{\Gamma}(\gamma)\}=\langle \gamma\rangle_{\Gamma}$. This trivially satisfies properties {\rm{(1,2)}}.
\item[Recursive step:] Suppose $\alpha_n$ has been constructed. There are two cases: 
\begin{enumerate}[label={\rm{Case {(\alph*)}:}}, topsep=0pt, wide]
\item Suppose $w_{n+1}\notin C$ for any $C \in\mathcal{C}$. By \ref{itm:neuman}, we can extend $\alpha_n$ to an isomorphism $\alpha_{n+1}:\ A_{n+1}\rightarrow B_{n+1}$ such that every fixed point of $\alpha_{n+1}$ is contained in $A_n$. Furthermore, we can assume $A_{n+1}=\langle A_n\, w_{n+1} u_{n+1}\rangle_{\Gamma}$ and $B_{n+1}=\langle B_n\, v_{n+1}w_{n+1}\rangle_{\Gamma}$ where $v_{n+1}\coloneqq\alpha_{n+1}(w_{n+1})$ and $u_{n+1}\coloneqq \alpha^{-1}_{n+1}(w_{n+1})$. 

As $\alpha_{n+1}$ has the same fixed points as $\alpha_n$, {\rm{(1)}} follows by recursion hypothesis. On the other hand, note that $\bigcup\mathcal{C}$ is definable over $\gamma$ in $\Tt^M_\Gamma$. Therefore, $w_{n+1},u_{n+1},v_{n+1}\notin \bigcup \mathcal{C}$, so $\cone_\gamma(A_{n+1})\cap\bigcup\mathcal{C}=\cone_\gamma(A_n)\cap\bigcup\mathcal{C}$ and $\cone_\gamma(B_{n+1})\cap \bigcup \mathcal{C}=\cone_\gamma(B_n)\cap \bigcup \mathcal{C}$ by \cref{l:generated}(\ref{itm:generated 8}). Thus, $\alpha_{n+1}$ satisfies {\rm{(2)}} as $\alpha_n$ does.
\item Suppose $w_{n+1}\in C$ for some $C\in\mathcal{C}$. 
\begin{enumerate}[topsep=0pt, wide]
\item[First we add $w_{n+1}$ to the domain.] There are two cases:
\begin{enumerate}[topsep=0pt, wide]
\item[Case $w_{n+1}\in \cone_\gamma(A_n)$.] By \ref{itm:neuman}, we extend $\alpha_n$ to an isomorphism $\alpha'_{n+1}:\ A'_{n+1}\rightarrow B'_{n+1}$ such that $A'_{n+1}=\langle A_n\, w_{n+1}\rangle_{\Gamma}$ and $B'_{n+1}=\langle B_n\, v_{n+1}\rangle_{\Gamma}$, where $\alpha'_{n+1}:\ w_{n+1}\mapsto v_{n+1}$, such that every fixed point of $\alpha_{n+1}$ is contained in $A_n$.

As $\alpha'_{n+1}$ has the same fixed points as $\alpha_n$, {\rm{(1)}} follows by recursion hypothesis. On the other hand, since $w_{n+1}\in\cone_\gamma(A_n)$, we get $\cone_\gamma(A'_{n+1})=\cone_\gamma(A_n)$ by \cref{l:generated}(\ref{itm:generated 8}). Thus, $\alpha'_{n+1}$ satisfies {\rm{(2)}} as $\alpha_n$ does.
\item[Case $w_{n+1}\notin\cone_\gamma(A_n)$.] Then, there is no $b\in B_n$ such that $b\in \sigma(C)$. By \cref{l:saturation of cones}, there is $v_{n+1}\in \sigma(C)$ such that $w_{n+1}\bar{a}\equiv^{\Rr\gamma}v_{n+1}\bar{b}$, where $\bar{a}$ enumerates $A_n$ and $\bar{b}$ is the enumeration of $B_n$ given by $\alpha_n:\ \bar{a}\mapsto \bar{b}$. Furthermore, note that we can take $v_{n+1}\neq w_{n+1}$, as according to \cref{l:saturation of cones}, we can take $v_{n+1}>w_{n+1}$ in case $\sigma(C)=C$. Take the extension of $\alpha_n$ to $A'_{n+1}=A_n\cup \{w_{n+1}\}$ and $B'_{n+1}=B_n\cup \{v_{n+1}\}$ setting $\alpha'_{n+1}:\ w_{n+1}\mapsto v_{n+1}$. 

Note that $A'_{n+1}=\langle A'_{n+1}\rangle_{\Gamma}$ and $B'_{n+1}=\langle B'_{n+1}\rangle_{\Gamma}$. Indeed, for any $a\in A_n$, $a\meet w_{n+1}\leq \gamma$. Thus, by \cref{l:generated}(\ref{itm:generated 2}), $w_{n+1}\meet a=w_{n+1}\meet a\meet \gamma=\min\{a\meet \gamma,\gamma\}\in A_n$ for all $a\in A_n$. On the other hand, since $w_{n+1}$ is in a cone disjoint to $\Gamma$, we get $\pi_\Gamma(w_{n+1})=\pi_\Gamma(w_{n+1})\meet w_{n+1}\leq \gamma\leq w_{n+1}$. Therefore, $\pi_{\Gamma}(w_{n+1})=\pi_\Gamma(\gamma)\in A_n$. Hence, by \cref{l:generated}(\ref{itm:generated 4},\ref{itm:generated 6}), we get $A'_{n+1}=\langle A'_{n+1}\rangle_{\Gamma}$. Similarly, we conclude $B'_{n+1}=\langle B'_{n+1}\rangle_{\Gamma}$.

Since $\Gamma\cap A'_{n+1}=\Gamma\cap A_n$ and $\Gamma\cap B'_{n+1}=\Gamma\cap B_n$, we have that $\alpha'_{n+1}$ is an isomorphism between $A'_{n+1}$ and $B'_{n+1}$ as substructures of $\Tt^M_\Gamma$. As $\gamma\in A_n$ and $\alpha'_{n+1}(\gamma)=\alpha_n(\gamma)=\gamma$ by recursion hypothesis, we conclude that $\alpha'_{n+1}$ satisfies condition {\rm{(1)}}. On the other hand, by choice of $v_{n+1}$ and recursion hypothesis, we get that $\alpha'_{n+1}$ satisfies property {\rm{(2)}}.
\end{enumerate}
\item[] In a similar manner, we extend $\alpha'_{n+1}$ to an isomorphism $\alpha_{n+1}$ satisfying properties {\rm{(1,2)}} by adding $w_{n+1}$ to the image. 
\end{enumerate}
\end{enumerate} 
\end{enumerate}

Finally, taking $\alpha=\bigcup \alpha_n$, we get an automorphism of $\Tt^M_\Gamma$ fixing $\gamma$ and satisfying conditions {\rm{(1,2)}}. 
\end{proof}
\end{lem} 
\begin{coro} \label{c:trivial obs} Let $\gamma$ be a point and $\sigma$ a permutation of the set of cones above $\gamma$. There is an automorphism of $\Tt^M$ fixing $\gamma$ whose induced action on the cones above $\gamma$ is $\sigma$. 
\begin{proof} By \cref{l:trivial obs} taking any branch $\Gamma$ that does not contain $\gamma$. 
\end{proof}
\end{coro}

\begin{lem} \label{l:infinite fan} Let $\gamma$ be a point and $\alpha$ an automorphism of $\Tt^M$ fixing $\gamma$. Assume that the set of cones above $\gamma$ which $\alpha$ moves is infinite (i.e. the induced permutation on the cones above $\gamma$ has infinite support). Then, there are $8$ conjugates of $\alpha$ or $\alpha^{-1}$ whose product is an infinite fan above $\gamma$ (see \cref{d:fan}).
\begin{proof} Let $\mathcal{C}=\{C_i\}_{i\in\Z}$ be the set of cones above $\gamma$ and consider the permutation $\sigma$ on $\mathcal{C}$ induced by $\alpha$. By \cite[Corollary 2.10]{tent2013isometry}\footnote{To apply the corollary we use the trivial independence relation on an infinite set $A \indepe{B}C$ iff $A \cap C \subseteq B$. For any permutation $\sigma$ with infinite support, any finite set $X$ and any element $b\notin X$, there is some $b'\notin X$ such that $\sigma(b') \notin X \cup \sigma(X) \cup \{b'\}$. Thus, by definition and \cite[Lemma 2.6]{tent2013isometry}, $\sigma$ moves maximally.} there are $\sigma_1,\ldots,\sigma_8$ permutations of $\mathcal{C}$ and signs $\epsilon_1,\ldots,\epsilon_8\in\{-1,+1\}$ such that $(\sigma^{\epsilon_1})^{\sigma_1}\cdots(\sigma^{\epsilon_8})^{\sigma_8}$ is the permutation $\mathcal{C}_i\mapsto \mathcal{C}_{i+1}$ for each $i\in \Z$. For each $j\in\{1,\ldots,8\}$, take any $\beta_j\in\Aut(\Tt^M_\gamma)$ whose induced action on $\mathcal{C}$ is $\sigma_j$, which exists by \cref{c:trivial obs}. Then, $(\alpha^{\epsilon_1})^{\beta_1}\cdots (\alpha^{\epsilon_8})^{\beta_8}$ acts on $\mathcal{C}$ by $\mathcal{C}_i\mapsto \mathcal{C}_{i+1}$ for each $i\in\Z$, so it is an infinite fan above $\gamma$. 
\end{proof}
\end{lem}

\begin{lem}\label{l:construction 1} Let $\Gamma$ be a semibranch and $\alpha$ an automorphism of $\Tt^M$ setwise fixing $\Gamma$. Suppose that $\alpha$ is unboundedly increasing in $\Gamma^{\circ}$ (see \cref{d:unbounded}). Then, there is an automorphism $\beta$ of $\Tt^M$ setwise fixing $\Gamma$ such that:
\begin{enumerate}[label={\rm{(\alph*)}}, wide]
\item $\alpha\beta^{-1}\alpha\beta$ is strictly increasing in $\Gamma^{\circ}$. 
\item If $\Gamma$ is not a branch, $\beta$ setwise fixes every cone above $\max\Gamma$. 
\end{enumerate}
\begin{proof} Let $W\coloneqq\{w \tq \pi_\Gamma(w)\in \Gamma^{\circ}\}$ and $V\coloneqq\{v\tq \pi_\Gamma(v)\notin \Gamma^{\circ}\}$. Let $(w_n)_{n\in\N}$ be an enumeration of $W$ and $(v_n)_{n\in\N}$ an enumeration of $V$.

By back{\hyp}and{\hyp}forth, we recursively construct a chain of partial finite isomorphisms $(\beta_n)_{n\in\N}$ of $\Tt^M$ setwise fixing $\Gamma$ and an increasing sequence $(a_n)_{n\in\Z}$ of elements in $\Gamma^\circ$ with the following properties:
\begin{enumerate}[label={(\roman*)}, wide]
\item $\beta_n:\ X_n\rightarrow Y_n$ with $X_n=\langle X_n\rangle_{\Gamma}$ and $Y_n=\langle Y_n\rangle_{\Gamma}$, and $w_k,v_k\in X_n$ and $w_k,v_k\in Y_n$ for $k<n$.
\item For each $-n<k<n$, we have $\alpha\beta_n(a_{k-1})=\beta_n\alpha^{-1}(a_k)$. In particular, $a_k\in X_n$ for $-n-1<k<n-1$, and $\alpha^{-1}(a_k)\in X_n$ for $-n<k<n$.
\item We have $\alpha(a_{k})>a_{k}$ for $-n-1<k<n$, and $\alpha(c_n)>c_n$ where $c_n\coloneqq\beta_n^{-1}\alpha\beta_n(a_{n-1})$. In particular, $a_{n-1}\in X_n$ and $\alpha\beta_n(a_{n-1})\in Y_n$.
\item For $n\geq 1$, $a_{-n}=\min X_n\cap \Gamma^{\circ}$ and $c_n=\max X_n\cap \Gamma^\circ$.
\item (In case $\Gamma\neq \Gamma^{\circ}$) $\beta_n(v)$ and $v$ are in the same cone above $\gamma\coloneqq \max\Gamma$ for any $v\in X_n$ with $v>\gamma$. Also, $\gamma\in X_n$, $\gamma\in Y_n$ and $\beta_n(\gamma)=\gamma$.
\end{enumerate}

In stage $n$ we construct $\beta_n$ and $(a_k)_{-n\leq k<n}$. 
\begin{enumerate}[wide]
\item[Case $n=0$:] If $\Gamma$ is not a branch, we take $\beta_{0}:\ \gamma\mapsto \gamma$, with $X_{0}=Y_{0}=\{\gamma\}$, where $\gamma\coloneqq\max\Gamma$. If $\Gamma$ is a branch, we take $\beta_0=\emptyset$.
\item[Case $n=1$:] Take any $a_0\in\Gamma^{\circ}$ with $\alpha(a_0)>a_0$. Applying sequentially \cref{l:saturation of increasing points}(\ref{itm:saturation of increasing points 1}) twice and then \cref{l:saturation of increasing points}(\ref{itm:saturation of increasing points 2}) twice, we find $b,c_1,b',a_{-1}\in\Gamma^{\circ}$ with $b'<\pi_\Gamma(w_0)<b$ and $a_{-1}<\pi_\Gamma(w_0)<c_1$ such that $\alpha(b)>b$, $\alpha(c_1)>c_1$, $\alpha(b')>b'$, $\alpha(a_{-1})>a_{-1}$ and $a_0c_1\alpha^{-1}(a_0)a_{-1}\equiv^{\Rr\Gamma}b\alpha(b)b'\alpha^{-1}(b')$. By homogeneity, we extend $\beta_0$ to $\beta'_1:\ X'_1\rightarrow Y'_1$ with $\beta'_1:\ a_0c_1\alpha^{-1}(a_0)a_{-1}\mapsto b\alpha(b)b'\alpha^{-1}(b')$ and $w_0\in X'_1\cap Y'_1$.  Further, we restrict $\beta'_1$ so that $X'_1=\langle X_0\, a_0c_1\alpha^{-1}(a_0)a_{-1}w_0{\beta'}^{-1}_1(w_0)\rangle_\Gamma$. If $\Gamma=\Gamma^{\circ}$, we set $\beta_1=\beta_1'$. If $\Gamma\neq\Gamma^\circ$, let $\bar{x}'$ enumerate $X'_1$ and $\bar{y}'$ be the enumeration of $Y'_1$ given by $\beta_1':\ \bar{x}'\mapsto \bar{y}'$. By \cref{l:saturation of cones} and homogeneity, find $u,u'\in \cone_\gamma(v_0)$ such that $v_0u'\bar{x}'\equiv^{\Rr\Gamma}uv_0\bar{y}'$. By homogeneity, we extend $\beta'_1$ to $\beta_1:\ X_1\rightarrow Y_1$ where $X_1=\langle X'_1\, v_0 u'\rangle_\Gamma$ and $v_0u'\mapsto uv_0$. 

Obviously, $\beta_1$ satisfies {\rm{(i)}}. For $-1<k<1$, we have $k=0$, and $\beta_1\alpha^{-1}(a_0)=b'=\alpha\beta_1(a_{-1})$, so {\rm{(ii)}} holds. By construction, $c_1=\beta_1^{-1}\alpha\beta_1(a_0)$ and we have $\alpha(a_k)>a_k$ for $-2<k<1$ and $\alpha(c_1)>c_1$, so we get {\rm{(iii)}}. By construction, $\alpha^{-1}(b')<b'<\pi_\Gamma(w_0)<b<\alpha(b)$, so applying $\beta_1^{-1}$ we get $a_{-1}<\alpha^{-1}(a_0)<\pi_\Gamma(\beta_1^{-1}(w_0))<a_0<c_1$. On the other hand, $a_{-1}<\pi_\Gamma(w_0)<c_1$ too. Hence, by \cref{l:generated}(\ref{itm:generated 6}), we get that $c_1=\max X_1\cap \Gamma^{\circ}$ and $a_{-1}=\min X_1\cap \Gamma^{\circ}$, i.e. {\rm{(iv)}} holds. Finally, in case $\Gamma\neq \Gamma^{\circ}$, by \cref{l:generated}(\ref{itm:generated 8}), we have $\cone(X_1)=\cone_\gamma(v_0\, u')$ and $\cone_\gamma(v_0)=\cone_\gamma(u)=\cone_\gamma(u')$, so {\rm{(v)}} holds.    
\item[Recursion case (case $n+1$):] Suppose $\beta_n:\ X_n\rightarrow Y_n$ and $(a_k)_{-n\leq k<n}$ are already given for $n\geq 1$. Let $\bar{x}$ enumerate $X_n$ and $\bar{y}$ be the enumeration of $Y_n$ given by $\beta_n:\ \bar{x}\mapsto \bar{y}$ .
\begin{enumerate}[label={\rm{Step {\arabic*}:}}, wide]
\item By recursion hypothesis {\rm{(iv)}}, we know that $c_n=\max X_n\cap \Gamma^\circ$ where $c_n\coloneqq\beta_n^{-1}\alpha\beta_n(a_{n-1})$. Set $a_n\coloneqq\alpha(c_n)$. By recursion hypothesis {\rm{(iii)}}, $a_n=\alpha(c_n)>c_n$. Hence, by \cref{l:saturation of increasing points}(\ref{itm:saturation of increasing points 1}), we can find $b\in\Gamma^\circ$ with $\alpha(b)>b$ and $b>\pi_{\Gamma}(w_n)$ such that $a_n\bar{x}\equiv^{\Rr\Gamma} b\bar{y}$. 
\item Applying \cref{l:saturation of increasing points}(\ref{itm:saturation of increasing points 1}) again, we can find $c_{n+1}\in\Gamma^{\circ}$ with $\alpha(c_{n+1})>c_{n+1}$ and $c_{n+1}>\pi_{\Gamma}(w_{n})$ such that $c_{n+1}a_n\bar{x}\equiv^{\Rr\Gamma}\alpha(b)b\bar{y}$. 
\begin{figure}[ht]
\begin{center}
\begin{tikzpicture}[scale=0.4]
\node (xL) at (-5*3+1.7,-4) {};  \node (yL) at (-5*3+1.7,4) {}; \node (xR) at (5*3,-4) {}; \node (yR) at (5*3,4) {};


\node (topL) [label={[xshift=5, yshift=-10]:$\Gamma$}] at (-2*3,3) {}; \node (botL) at (-2*3,-3) {}; \node (piwnL) [draw, shape=circle, fill, scale=0.15, label={right:$\pi_\Gamma(w_n)$}] at (-2*3,0) {}; \node (rL) at (-3*3,0.75) {}; \node (rmeetwnL) [draw, shape=coordinate] at (-2.4*3,0.3) {}; \node (wnL) [draw, shape=circle, fill, scale=0.15, label=above:$w_n$] at (-2.6*3,2) {};

\node (topR) [label={[xshift=5, yshift=-10]:$\Gamma$}] at (2*3,3) {}; \node (botR) at (2*3,-3) {}; \node (rR) at (1*3,0.75) {}; \node (rmeetwnR) [draw, shape=coordinate] at (1.6*3,0.3) {}; \node (wnR) [draw, shape=circle, fill, scale=0.15, label=above:$w_n$] at (1.4*3,2) {}; \node (piwnR) [draw, shape=circle, fill, scale=0.15, label={right:$\pi_\Gamma(w_n)$}] at (2*3,0) {};

\node (cn) [draw, shape=circle, fill, scale=0.15, label=right:$c_n$] at (-2*3,-2) {}; 
\node (an) [draw, shape=circle, fill, scale=0.15, label={right:$a_n=\alpha(c_n)$}] at (-2*3,-1) {}; 
\node (b) [draw, shape=circle, fill, scale=0.15, label={right:$b$}] at (2*3,1) {}; 
\node (alphab) [draw, shape=circle, fill, scale=0.15, label=right:$\alpha(b)$] at (2*3,2) {};
\node (cn1) [draw, shape=circle, fill, scale=0.15, label={right:$c_{n+1}$}] at (-2*3,1) {}; 
\node (alphacn1) [draw, shape=circle, fill, scale=0.15, label={right:$\alpha(c_{n+1})$}] at (-2*3,2) {};

\draw (botL) -- (topL); \draw (piwnL) -- (rL); \draw (rmeetwnL) -- (wnL);
\draw (botR) -- (topR); \draw (piwnR) -- (rR); \draw (rmeetwnR) -- (wnR);
\draw[->, shorten <=20, shorten >=20] (topL) to [bend left=20] node[auto]{$\begin{array}{c}a_n\mapsto b\\ c_{n+1}\mapsto \alpha(b)\end{array}$} (topR) ;
\draw[->, shorten <=1.5, shorten >=1.2] (cn) to [out=180-45, in=180+45] node[auto]{$\alpha$} (an) ;
\draw[->, shorten <=1.5, shorten >=1.2] (cn1) to [out=180-45, in=180+45] node[auto]{$\alpha$} (alphacn1);
\draw[->, shorten <=1.5, shorten >=1.2] (b) to [out=180-45, in=180+45] node[auto]{$\alpha$} (alphab);
\end{tikzpicture}
\caption{Illustration of steps 1 and 2.}
\subcaption*{(In general, the position of $w_n$ respect to $a_n$ and $c_n$ is unknown --- the position drawn is just one possibility)}
\end{center}
\end{figure}
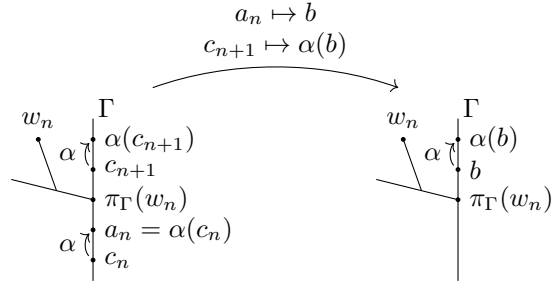
\item By recursion hypothesis {\rm{(iii,iv)}}, $\alpha^{-1}(a_{-n})<a_{-n}=\min X_n\cap\Gamma^{\circ}$. Hence, by \cref{l:saturation of increasing points}(\ref{itm:saturation of increasing points 2}), we can find $b'\in \Gamma^{\circ}$ such that $\alpha^{-1}(a_{-n})c_{n+1}a_n\bar{x}\equiv^{\Rr\Gamma} b'\alpha(b)b\bar{y}$ with $\alpha(b')>b'$ and $b'<\pi_{\Gamma}(w_n)$. 
\item Then, $\alpha^{-1}(b')<b'$, so, applying \cref{l:saturation of increasing points}(\ref{itm:saturation of increasing points 2}) again, we can find $a_{-n-1}\in \Gamma^{\circ}$ such that $a_{-n-1}\alpha^{-1}(a_{-n})c_{n+1}a_n\bar{x}\equiv^{\Rr\Gamma}\alpha^{-1}(b')b'\alpha(b)b\bar{y}$ with $\alpha(a_{-n-1})>a_{-n-1}$ and $a_{-n-1}<\pi_{\Gamma}(w_{n})$. 
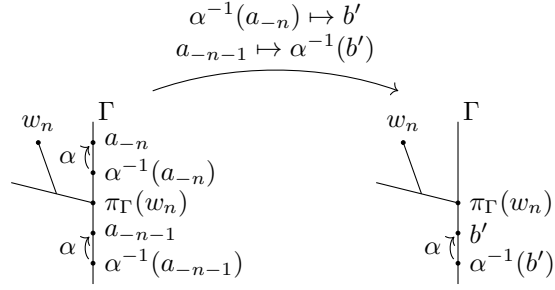
\begin{figure}[ht]
\begin{center}
\begin{tikzpicture}[scale=0.4]
\node (xL) at (-5*3+1.7,-4) {};  \node (yL) at (-5*3+1.7,4) {}; \node (xR) at (5*3,-4) {}; \node (yR) at (5*3,4) {};


\node (topL) [label={[xshift=5, yshift=-10]:$\Gamma$}] at (-2*3,3) {}; \node (botL) at (-2*3,-3) {}; \node (piwnL) [draw, shape=circle, fill, scale=0.15, label={right:$\pi_\Gamma(w_n)$}] at (-2*3,0) {}; \node (rL) at (-3*3,0.75) {}; \node (rmeetwnL) [draw, shape=coordinate] at (-2.4*3,0.3) {}; \node (wnL) [draw, shape=circle, fill, scale=0.15, label=above:$w_n$] at (-2.6*3,2) {};

\node (topR) [label={[xshift=5, yshift=-10]:$\Gamma$}] at (2*3,3) {}; \node (botR) at (2*3,-3) {}; \node (rR) at (1*3,0.75) {}; \node (rmeetwnR) [draw, shape=coordinate] at (1.6*3,0.3) {}; \node (wnR) [draw, shape=circle, fill, scale=0.15, label=above:$w_n$] at (1.4*3,2) {}; \node (piwnR) [draw, shape=circle, fill, scale=0.15, label={right:$\pi_\Gamma(w_n)$}] at (2*3,0) {};

\node (alphaan1) [draw, shape=circle, fill, scale=0.15, label=right:$\alpha^{-1}(a_{-n-1})$] at (-2*3,-2) {}; 
\node (an1) [draw, shape=circle, fill, scale=0.15, label={right:$a_{-n-1}$}] at (-2*3,-1) {}; 
\node (alphab) [draw, shape=circle, fill, scale=0.15, label={right:$\alpha^{-1}(b')$}] at (2*3,-2) {}; 
\node (b) [draw, shape=circle, fill, scale=0.15, label={right:$b'$}] at (2*3,-1) {};
\node (alphaan) [draw, shape=circle, fill, scale=0.15, label={right:$\alpha^{-1}(a_{-n})$}] at (-2*3,1) {}; 
\node (an) [draw, shape=circle, fill, scale=0.15, label={right:$a_{-n}$}] at (-2*3,2) {};

\draw (botL) -- (topL); \draw (piwnL) -- (rL); \draw (rmeetwnL) -- (wnL);
\draw (botR) -- (topR); \draw (piwnR) -- (rR); \draw (rmeetwnR) -- (wnR);
\draw[->, shorten <=20, shorten >=20] (topL) to [bend left=20] node[auto]{$\begin{array}{c} \alpha^{-1}(a_{-n})\mapsto b' \\ a_{-n-1}\mapsto \alpha^{-1}(b') \end{array}$} (topR) ;
\draw[->, shorten <=1.5, shorten >=1.2] (alphaan1) to [out=180-45, in=180+45] node[auto]{$\alpha$} (an1) ;
\draw[->, shorten <=1.5, shorten >=1.2] (alphaan) to [out=180-45, in=180+45] node[auto]{$\alpha$} (an);
\draw[->, shorten <=1.5, shorten >=1.2] (alphab) to [out=180-45, in=180+45] node[auto]{$\alpha$} (b);
\end{tikzpicture}
\caption{Illustration of steps 3 and 4.} 
\subcaption*{(In general, the position of $w_n$ respect to $a_{-n}$ and $\alpha^{-1}(a_{-n})$ is unknown --- the position drawn is just one possibility)}
\end{center}
\end{figure}
\item By homogeneity, we extend $\beta_n$ to a partial isomorphism $\beta'_{n+1}:\ X'_{n+1}\rightarrow Y'_{n+1}$ with $\beta'_{n+1}:\ a_{-n-1}\alpha^{-1}(a_{-n})c_{n+1}a_n\mapsto \alpha^{-1}(b')b'\alpha(b)b$ and $w_n\in X'_{n+1}\cap Y'_{n+1}$. Further, we restrict $\beta'_{n+1}$ so that $X'_{n+1}=\langle X_n\, a_{-n-1}\alpha^{-1}(a_{-n})$ $c_{n+1}a_nw_n{\beta'}^{-1}_{n+1}(w_n)\rangle_\Gamma$. In case $\Gamma=\Gamma^{\circ}$, we conclude setting $\beta_{n+1}=\beta'_{n+1}$. Otherwise, continue to step 6.  
\item (In case $\Gamma\neq \Gamma^{\circ}$) Let $\bar{x}'$ enumerate $X'_{n+1}$ and $\bar{y}'$ be the enumeration of $Y'_{n+1}$ given by $\beta'_{n+1}:\ \bar{x}'\mapsto \bar{y}'$. By \cref{l:generated}(\ref{itm:generated 8}) and recursion hypothesis {\rm{(v)}}, we have $C\coloneqq\cone_\gamma(X'_{n+1})=\cone_\gamma(X_n)=\cone_\gamma(Y_n)=\cone_\gamma(Y'_{n+1})$. If $v_n\in C$, by homogeneity, extend $\beta'_{n+1}$ to a partial isomorphism $\beta_{n+1}$ by adding $v_n$ to the domain and the image. If $v_n\notin C$, applying \cref{l:saturation of cones}, we can find $u$ with $v_n \bar{x}'\equiv^{\Rr\Gamma}u \bar{y}'$ and $u$ in the same cone as $v_n$. By homogeneity, find $u'$ with $u'v_n\bar{x}'\equiv^{\Rr\Gamma}v_nu\bar{y}'$. Extend $\beta'_{n+1}$ to $\beta_{n+1}:\ X_{n+1}\rightarrow Y_{n+1}$ by mapping $v_nu'\mapsto uv_n$, where $X_{n+1}=\langle X'_{n+1}\, u' v_n\rangle_\Gamma$ and $Y_{n+1}=\langle Y'_{n+1}\, v_n u\rangle_{\Gamma}$.
\end{enumerate} 

Now, we check that $\beta_{n+1}$ and $(a_k)_{-n-1\leq k< n+1}$ just constructed satisfy the recursion hypotheses:
\begin{enumerate}[label={(\roman*)}, wide]
\item By construction, $\beta_{n+1}:\ X_{n+1}\rightarrow X_{n+1}$ is a partial isomorphism between substructures extending $\beta_n$ with $w_n,v_n\in X_{n+1}$ and $w_n,v_n\in Y_{n+1}$. 
\item By construction, $a_{-n-1},\alpha^{-1}(a_{-n}),a_n\in X_{n+1}$, $a_{-n}\in X_n\subseteq X_{n+1}$ and $\alpha^{-1}(a_n)=c_n\in X_{n}\subseteq X_{n+1}$. Also, by construction, $\beta_{n+1}\alpha^{-1}(a_n)=\beta_{n+1}(c_n)=\alpha\beta_n(a_{n-1})$. Finally, we get $\beta_{n+1}(a_{-n-1})=\alpha^{-1}(b')$, concluding $\alpha\beta_{n+1}(a_{-n-1})=\beta_{n+1}\alpha^{-1}(a_{-n})$.
\item By construction, $a_n=\alpha(c_n)$ where $c_n=\beta_n^{-1}\alpha\beta_n(a_{n-1})$. Thus, by recursion hypothesis {\rm{(iii)}}, $a_n=\alpha(c_n)>c_n$, so $\alpha(a_n)>\alpha(c_n)=a_n$. Also, by construction, $\alpha(a_{-n-1})>a_{-n-1}$. Finally, $\beta_{n+1}^{-1}\alpha\beta_{n+1}(a_n)=c_{n+1}$ and $\alpha(c_{n+1})>c_{n+1}$. 
\item By recursion hypotheses {\rm{(iii,iv)}}, $a_n=\alpha(c_n)>c_n=\max X_n\cap \Gamma^{\circ}$ and $\alpha^{-1}(a_{-n})<a_{-n}=\min X_n\cap \Gamma^{\circ}$. Since $\beta_{n+1}$ is a partial isomorphism, $b>\max Y_n\cap \Gamma^{\circ}$ and $b'<\min Y_n\cap \Gamma^{\circ}$. By construction, $\alpha(b)>b$ and $\alpha^{-1}(b')<b'$. As $\beta_{n+1}$ is a partial isomorphism, $c_{n+1}>a_n$ and $a_{-n-1}<\alpha^{-1}(a_{-n})$. Also, by construction, $a_{-n-1}<\pi_{\Gamma}(w_n)<c_{n+1}$ and $b'<\pi_{\Gamma}(w_n)<b$. Hence, by \cref{l:generated}(\ref{itm:generated 6}), $a_{-n-1}=\min X_{n+1}\cap \Gamma^{\circ}$ and $c_{n+1}=\max X_{n+1}\cap \Gamma^\circ$.
\item (In case $\Gamma\neq \Gamma^{\circ}$) If $v_n\sim_\gamma v\in X_n$, then $\beta_{n+1}(v_n)\sim_\gamma \beta_{n+1}(v)\sim_\gamma v\sim_\gamma v_n$ by recursion hypothesis {\rm{(v)}}. On the other hand, if $v_n\notin \cone(X_n)$, by construction, $\beta_{n+1}(v_n)$ is in the same cone as $v_n$.  By \cref{l:generated}(\ref{itm:generated 8}), (v) holds.
\end{enumerate}
\end{enumerate}

Set $\beta=\bigcup \beta_n$, so $\beta\in\Aut(\Tt^{M}_\Gamma)$ by property {\rm{(i)}}. 
\begin{enumerate}[label={\rm{(\alph*)}}, wide]
\item For any $x\in \Gamma^{\circ}$, $x=w_n$ for some $n\in\N$. Then, $a_{-n-1}=\min X_{n+1}\leq x\leq \max X_{n+1}=\beta^{-1}\alpha\beta(a_n)<\alpha\beta^{-1}\alpha\beta(a_n)=a_{n+1}$ by properties {\rm{(i,ii,iii,iv)}}. Take $m\coloneqq\min\{n\in\Z\tq x<a_n\}$. Then, $a_{m-1}\leq x< a_m$, so $x<a_m=\alpha\beta^{-1}\alpha\beta(a_{m-1})\leq \alpha\beta^{-1}\alpha\beta(x)$ by property {\rm{(ii)}}. As $x$ is arbitrary, we conclude that $\alpha\beta^{-1}\alpha\beta$ is strictly increasing in $\Gamma^{\circ}$. 
\item (In case $\Gamma\neq \Gamma^{\circ}$) For any $v>\gamma\coloneqq\max \Gamma$, $v=v_n$ for some $n\in\N$. Then, $\beta(v)$ is in the same cone as $v$ by property {\rm{(v)}}. As $v$ is arbitrary, we get that $\beta$ leaves invariant every cone above $\gamma$. \qedhere
\end{enumerate}
\end{proof}
\end{lem}

The proof of \cref{l:construction 1} is a blueprint for the proofs of the following lemmas, which are quite similar.

\begin{lem} \label{l:construction 2} Let $\Gamma$ be a semibranch and $\alpha$ an automorphism of $\Tt^M$ fixing $\Gamma$ setwise. Suppose that $\alpha$ is upwards unboundedly increasing and downwards unboundedly decreasing in $\Gamma^{\circ}$. Then, there are an automorphism $\beta$ of $\Tt^M$ fixing $\Gamma$ setwise and $a,a'\in\Gamma^{\circ}$ such that:
\begin{enumerate}[label={\rm{(\alph*)}}]
\item $((\alpha\beta^{-1}\alpha\beta)^n(a))_{n\in\N}$ is strictly increasing and upwards unbounded in $\Gamma^{\circ}$.
\item $((\alpha\beta^{-1}\alpha\beta)^n(a'))_{n\in\N}$ is strictly decreasing and downwards unbounded. 
\item If $\Gamma$ is not a branch, $\beta$ fixes setwise every cone above $\max\Gamma$. 
\end{enumerate}
\begin{proof} Let $W\coloneqq\{w\tq \pi_\Gamma(w)\in \Gamma^{\circ}\}$ and $V\coloneqq\{v\tq \pi_\Gamma(v)\notin \Gamma^{\circ}\}$. Let $(w_n)_{n\in\N}$ be an enumeration of $W$ and $(v_n)_{n\in\N}$ an enumeration of $V$.

By back{\hyp}and{\hyp}forth, we recursively construct a chain of partial finite isomorphisms $(\beta_n)_{n\in\N}$ of $\Tt^M$ setwise fixing $\Gamma$, an increasing sequence $(a_n)_{n\in\N}$ in $\Gamma^\circ$ and a decreasing sequence $(a'_n)_{n\in\N}$ in $\Gamma^\circ$ with the following properties:
\begin{enumerate}[label={(\roman*)}, wide]
\item $\beta_n:\ X_n\rightarrow Y_n$ with $X_n=\langle X_n\rangle_{\Gamma}$ and $Y_n=\langle Y_n\rangle_{\Gamma}$, and $w_k,v_k\in X_n$ and $w_k,v_k\in Y_n$ for $k<n$.
\item For each $0<k<n$, we have $\alpha\beta_n(a_{k-1})=\beta_n\alpha^{-1}(a_k)$ and $\alpha\beta_n(a'_{k-1})=\beta_n\alpha^{-1}(a'_k)$. In particular, $a_k,a'_k\in X_n $ for $k<n-1$, and $\alpha^{-1}(a_k),\alpha^{-1}(a'_k) \in X_n$ for $0<k<n$. 
\item We have $\alpha(a_{k})>a_{k}$ and $\alpha(a'_k)<a'_k$ for $k<n$, and $\alpha(c_n)>c_n$ and $\alpha(c'_n)<c'_n$ where $c_n\coloneqq \beta_n^{-1}\alpha\beta_n(a_{n-1})$ and $c'_n\coloneqq \beta_n^{-1}\alpha\beta_n(a'_{n-1})$. In particular, $a_{n-1},a'_{n-1}\in X_n$ and $\alpha\beta_n(a_{n-1}),\alpha\beta_n(a'_{n-1})\in Y_n$.
\item For $n\geq 1$, $c'_n=\min X_n\cap \Gamma^{\circ}$ and $c_n=\max X_n\cap \Gamma^\circ$.
\item (In case $\Gamma\neq\Gamma^{\circ}$) $\beta_n(v)$ and $v$ are in the same cone above $\gamma\coloneqq \max\Gamma$ for any $v\in X_n$ with $v>\gamma$. Also, $\gamma\in X_n$, $\gamma\in Y_n$ and $\beta_n(\gamma)=\gamma$.
\end{enumerate}

In stage $n$ we construct $\beta_n$, $(a_k)_{k<n}$ and $(a'_k)_{k<n}$.
\begin{enumerate}[wide]
\item[Case $n=0$:] If $\Gamma$ is not a branch, we take $\beta_{0}:\ \gamma\mapsto \gamma$, with $X_{0}=Y_{0}=\{\gamma\}$, where $\gamma\coloneqq\max\Gamma$. If $\Gamma$ is a branch, we take $\beta_0=\emptyset$.
\item[Case $n=1$:] Take any $a_0,a'_0\in\Gamma^{\circ}$ with $a_0>a_0'$ such that $\alpha(a_0)>a_0$ and $\alpha(a'_0)<a'_0$. Applying sequentially \cref{l:saturation of increasing points}(\ref{itm:saturation of increasing points 1}) twice and then \cref{l:saturation of increasing points}(\ref{itm:saturation of increasing points 2}) twice (with $\alpha^{-1}$), we find $b,c_1,b',c_1'\in\Gamma^{\circ}$ with $b'<\pi_\Gamma(w_0)<b$ and $c'_1<\pi_\Gamma(w_0)<c_1$ such that $\alpha(b)>b$, $\alpha(c_1)>c_1$, $b'>\alpha(b')$, $c'_1>\alpha(c'_1)$ and $a_0c_1a'_0c'_1\equiv^{\Rr\Gamma}b\alpha(b)b'\alpha(b')$. By homogeneity, we extend $\beta_0$ to $\beta'_1:\ X'_1\rightarrow Y'_1$ with $\beta_1':\ a_0c_1a'_0c'_1\mapsto b\alpha(b)b'\alpha(b')$ and $w_0\in X_1'\cap Y'_1$. Further, we restrict $\beta'_1$ so that $X'_1=\langle X_0\, a_0c_1a'_0c'_1w_0{\beta'}^{-1}_1(w_0)\rangle_\Gamma$. If $\Gamma=\Gamma^{\circ}$, we set $\beta_1=\beta_1'$. If $\Gamma\neq\Gamma^\circ$, by \cref{l:saturation of cones} and homogeneity, we extend $\beta'_1$ to $\beta_1:\ X_1\rightarrow Y_1$ by adding $v_0$ to the domain and the image such that $v_0$ and $\beta_1(v_0)$ are in the same cone. As above, we restrict so that $X_1=\langle X_0\, v_0 \beta_1^{-1}(v_0)\rangle_\Gamma$.

Obviously, $\beta_1$ satisfies {\rm{(i)}}. For {\rm{(ii)}} there is nothing to check. By construction, $c_1=\beta_1^{-1}\alpha\beta(a_0)$ and $c'_1=\beta_1^{-1}\alpha\beta_1(a'_0)$, and we have $\alpha(a_0)>a_0$, $\alpha(a'_0)<a'_0$, $\alpha(c_1)>c_1$ and $\alpha(c'_1)<c'_1$, so we get {\rm{(iii)}}. By construction, $\alpha(b')<b'<\pi_\Gamma(w_0)<b<\alpha(b)$, so applying $\beta^{-1}_1$ we get $c'_1<a'_0<\pi_\Gamma(\beta^{-1}(w_0))<a_0<c_1$. On the other hand, $c'_1<\pi_\Gamma(w_0)<c_1$ too. Hence, by \cref{l:generated}(\ref{itm:generated 6}), we get that $c_1=\max X_1\cap \Gamma^{\circ}$ and $c'_1=\min X_1\cap \Gamma^{\circ}$, so {\rm{(iv)}} holds. Finally, in case $\Gamma\neq\Gamma^{\circ}$, by \cref{l:generated}(\ref{itm:generated 8}), we have $\cone(X_1)=\cone_\gamma(v_0)$, so {\rm{(v)}} holds.
\item[Recursion case (case $n+1$):] Suppose $\beta_n:\ X_n\rightarrow Y_n$, $(a_k)_{k<n}$ and $(a'_k)_{k<n}$ are already given for $n\geq 1$. Let $\bar{x}$ enumerate $X_n$ and $\bar{y}$ be the enumeration of $Y_n$ given by $\beta_n:\ \bar{x}\mapsto \bar{y}$.
\begin{enumerate}[label={\rm{Step {\arabic*}:}}, wide]
\item By recursion hypothesis {\rm{(iv)}}, we know that $c_n=\max X_n\cap \Gamma^\circ$ where $c_n\coloneqq\beta_n^{-1}\alpha\beta_n(a_{n-1})$. Set $a_n\coloneqq\alpha(c_n)$. By recursion hypothesis {\rm{(iii)}}, $a_n=\alpha(c_n)>c_n$. Hence, by \cref{l:saturation of increasing points}(\ref{itm:saturation of increasing points 1}), we can find $b\in\Gamma^\circ$ with $\alpha(b)>b$ and $b>\pi_{\Gamma}(w_n)$ such that $a_n\bar{x}\equiv^{\Rr\Gamma} b\bar{y}$.
\item Applying \cref{l:saturation of increasing points}(\ref{itm:saturation of increasing points 1}) again, we can find $c_{n+1}\in\Gamma^{\circ}$ with $\alpha(c_{n+1})>c_{n+1}$ and $c_{n+1}>\pi_{\Gamma}(w_{n})$ such that $c_{n+1}a_n\bar{x}\equiv^{\Rr\Gamma}\alpha(b)b\bar{y}$.
\begin{figure}[ht]
\begin{center}
\begin{tikzpicture}[scale=0.4]
\node (xL) at (-5*3+1.7,-4) {};  \node (yL) at (-5*3+1.7,4) {}; \node (xR) at (5*3,-4) {}; \node (yR) at (5*3,4) {};


\node (topL) [label={[xshift=5, yshift=-10]:$\Gamma$}] at (-2*3,3) {}; \node (botL) at (-2*3,-3) {}; \node (piwnL) [draw, shape=circle, fill, scale=0.15, label={right:$\pi_\Gamma(w_n)$}] at (-2*3,0) {}; \node (rL) at (-3*3,0.75) {}; \node (rmeetwnL) [draw, shape=coordinate] at (-2.4*3,0.3) {}; \node (wnL) [draw, shape=circle, fill, scale=0.15, label=above:$w_n$] at (-2.6*3,2) {};

\node (topR) [label={[xshift=5, yshift=-10]:$\Gamma$}] at (2*3,3) {}; \node (botR) at (2*3,-3) {}; \node (rR) at (1*3,0.75) {}; \node (rmeetwnR) [draw, shape=coordinate] at (1.6*3,0.3) {}; \node (wnR) [draw, shape=circle, fill, scale=0.15, label=above:$w_n$] at (1.4*3,2) {}; \node (piwnR) [draw, shape=circle, fill, scale=0.15, label={right:$\pi_\Gamma(w_n)$}] at (2*3,0) {};

\node (cn) [draw, shape=circle, fill, scale=0.15, label=right:$c_n$] at (-2*3,-2) {}; 
\node (an) [draw, shape=circle, fill, scale=0.15, label={right:$a_n=\alpha(c_n)$}] at (-2*3,-1) {}; 
\node (b) [draw, shape=circle, fill, scale=0.15, label={right:$b$}] at (2*3,1) {}; 
\node (alphab) [draw, shape=circle, fill, scale=0.15, label=right:$\alpha(b)$] at (2*3,2) {};
\node (cn1) [draw, shape=circle, fill, scale=0.15, label={right:$c_{n+1}$}] at (-2*3,1) {}; 
\node (alphacn1) [draw, shape=circle, fill, scale=0.15, label={right:$\alpha(c_{n+1})$}] at (-2*3,2) {};

\draw (botL) -- (topL); \draw (piwnL) -- (rL); \draw (rmeetwnL) -- (wnL);
\draw (botR) -- (topR); \draw (piwnR) -- (rR); \draw (rmeetwnR) -- (wnR);
\draw[->, shorten <=20, shorten >=20] (topL) to [bend left=20] node[auto]{$\begin{array}{c}a_n\mapsto b\\ c_{n+1}\mapsto \alpha(b)\end{array}$} (topR) ;
\draw[->, shorten <=1.5, shorten >=1.2] (cn) to [out=180-45, in=180+45] node[auto]{$\alpha$} (an) ;
\draw[->, shorten <=1.5, shorten >=1.2] (cn1) to [out=180-45, in=180+45] node[auto]{$\alpha$} (alphacn1);
\draw[->, shorten <=1.5, shorten >=1.2] (b) to [out=180-45, in=180+45] node[auto]{$\alpha$} (alphab);
\end{tikzpicture}
\caption{Illustration of steps 1 and 2.}
\subcaption*{(In general, the position of $w_n$ respect to $a_n$ and $c_n$ is unknown --- the position drawn is just one possibility)}
\end{center}
\end{figure}
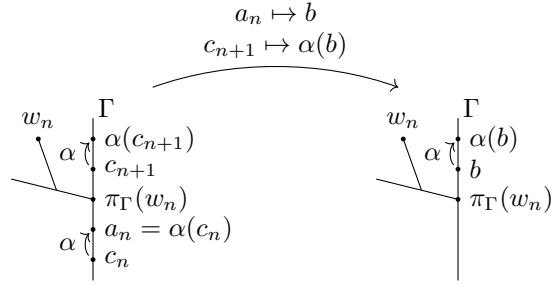
\item By recursion hypothesis {\rm{(iv)}}, we know that $c'_n=\min X_n\cap\Gamma^{\circ}$ where $c'_n=\beta_n^{-1}\alpha\beta_n(a'_{n-1})$. Set $a'_n=\alpha(c'_n)$. By recursion hypothesis {\rm{(iii)}}, $a'_n=\alpha(c'_n)<c'_n$. Hence, by \cref{l:saturation of increasing points}(\ref{itm:saturation of increasing points 2}) applied to $\alpha^{-1}$, we can find $b'\in\Gamma^\circ$ such that $a'_nc_{n+1}a_n\bar{x}\equiv^{\Rr\Gamma} b'\alpha(b)b\bar{y}$ with $\alpha(b')<b'$ and $b'<\pi_{\Gamma}(w_n)$.
\item Applying \cref{l:saturation of increasing points}(\ref{itm:saturation of increasing points 2}) to $\alpha^{-1}$ again, we can find $c'_{n+1}\in\Gamma^{\circ}$ with $\alpha(c'_{n+1})<c'_{n+1}$ and $c'_{n+1}<\pi_{\Gamma}(w_{n})$ such that $c'_{n+1}a'_nc_{n+1}a_n\bar{x}\equiv^{\Rr\Gamma}\alpha(b')b'\alpha(b)b\bar{y}$.
\begin{figure}[ht]
\begin{center}
\begin{tikzpicture}[scale=0.4]
\node (xL) at (-5*3+1.7,-4) {};  \node (yL) at (-5*3+1.7,4) {}; \node (xR) at (5*3,-4) {}; \node (yR) at (5*3,4) {};


\node (topL) [label={[xshift=5, yshift=-10]:$\Gamma$}] at (-2*3,3) {}; \node (botL) at (-2*3,-3) {}; \node (piwnL) [draw, shape=circle, fill, scale=0.15, label={right:$\pi_\Gamma(w_n)$}] at (-2*3,0) {}; \node (rL) at (-3*3,0.75) {}; \node (rmeetwnL) [draw, shape=coordinate] at (-2.4*3,0.3) {}; \node (wnL) [draw, shape=circle, fill, scale=0.15, label=above:$w_n$] at (-2.6*3,2) {};

\node (topR) [label={[xshift=5, yshift=-10]:$\Gamma$}] at (2*3,3) {}; \node (botR) at (2*3,-3) {}; \node (rR) at (1*3,0.75) {}; \node (rmeetwnR) [draw, shape=coordinate] at (1.6*3,0.3) {}; \node (wnR) [draw, shape=circle, fill, scale=0.15, label=above:$w_n$] at (1.4*3,2) {}; \node (piwnR) [draw, shape=circle, fill, scale=0.15, label={right:$\pi_\Gamma(w_n)$}] at (2*3,0) {};

\node (alphacn1) [draw, shape=circle, fill, scale=0.15, label=right:$\alpha(c'_{n+1})$] at (-2*3,-2) {}; 
\node (cn1) [draw, shape=circle, fill, scale=0.15, label={right:$c'_{n+1}$}] at (-2*3,-1) {}; 
\node (alphab) [draw, shape=circle, fill, scale=0.15, label={right:$\alpha(b')$}] at (2*3,-2) {}; 
\node (b) [draw, shape=circle, fill, scale=0.15, label={right:$b'$}] at (2*3,-1) {};
\node (an) [draw, shape=circle, fill, scale=0.15, label={right:$a'_n=\alpha(c'_n)$}] at (-2*3,1) {}; 
\node (cn) [draw, shape=circle, fill, scale=0.15, label={right:$c'_n$}] at (-2*3,2) {};

\draw (botL) -- (topL); \draw (piwnL) -- (rL); \draw (rmeetwnL) -- (wnL);
\draw (botR) -- (topR); \draw (piwnR) -- (rR); \draw (rmeetwnR) -- (wnR);
\draw[->, shorten <=20, shorten >=20] (topL) to [bend left=20] node[auto]{$\begin{array}{c} a'_n\mapsto b' \\ c'_{n+1}\mapsto \alpha(b') \end{array}$} (topR) ;
\draw[<-, shorten <=1.5, shorten >=1.2] (alphacn1) to [out=180-45, in=180+45] node[auto]{$\alpha$} (cn1) ;
\draw[<-, shorten <=1.5, shorten >=1.2] (an) to [out=180-45, in=180+45] node[auto]{$\alpha$} (cn);
\draw[<-, shorten <=1.5, shorten >=1.2] (alphab) to [out=180-45, in=180+45] node[auto]{$\alpha$} (b);
\end{tikzpicture}
\caption{Illustration of steps 3 and 4.}
\subcaption*{(In general, the position of $w_n$ respect to $a'_n$ and $c'_n$ is unknown --- the position drawn is just an illustration)}
\end{center}
\end{figure}
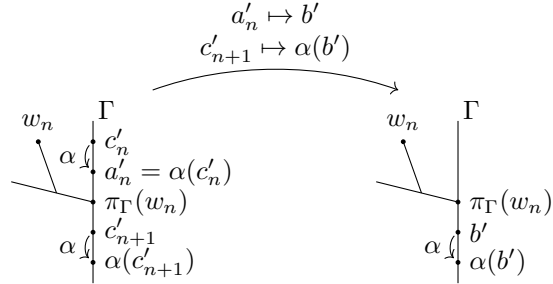
\item By homogeneity, we extend $\beta_n$ to a partial isomorphism $\beta'_{n+1}:\ X'_{n+1}\rightarrow Y'_{n+1}$ with $\beta'_{n+1}:\ c'_{n+1}a'_nc_{n+1}a_n\mapsto \alpha(b')b'\alpha(b)b$ and $w_n\in X'_{n+1}\cap Y'_{n+1}$. Further, we restrict $\beta'_{n+1}$ so that $X'_{n+1}=\langle X_n\, c'_{n+1}a'_nc_{n+1}a_nw_n{\beta'}^{-1}_{n+1}(w_n)\rangle_\Gamma$. In case $\Gamma=\Gamma^{\circ}$, we conclude setting $\beta_{n+1}=\beta'_{n+1}$. Otherwise, continue to step 6.  
\item (In case $\Gamma\neq \Gamma^{\circ}$) By recursion hypothesis {\rm{(v)}}, \cref{l:saturation of cones} and homogeneity, we can extend $\beta''_{n+1}$ to a partial elementary map $\beta_{n+1}$ by adding $v_n$ to the domain and the image such that $v_n$ and $\beta_{n+1}(v_n)$ are in the same cone. As above, we restrict so that $X_{n+1}=\langle X'_n\, v_n \beta_{n+1}^{-1}(v_n)\rangle_\Gamma$.
\end{enumerate}

Now, we check that $\beta_{n+1}$, $(a_k)_{k<n+1}$ and $(a'_k)_{k<n+1}$ just constructed satisfy the recursion hypotheses:
\begin{enumerate}[label={(\roman*)}, wide]
\item By construction, $\beta_{n+1}:\ X_{n+1}\rightarrow Y_{n+1}$ is a partial isomorphism extending $\beta_n$ with $w_n,v_n\in X_{n+1}$ and $w_n,v_n\in Y_{n+1}$. 
\item By construction, $c'_{n+1},a'_n,c_{n+1},a_n\in X_{n+1}$. Also, we have $\alpha^{-1}(a_n)=c_n\in X_{n}\subseteq X_{n+1}$ and $\alpha^{-1}(a'_n)=c'_n\in X_n\subseteq X_{n+1}$, so $\beta_{n+1}\alpha^{-1}(a_n)=\beta_{n+1}(c_n)=\alpha\beta_n(a_{n-1})$ and $\beta_{n+1}\alpha^{-1}(a'_n)=\beta_{n+1}(c'_n)=\alpha\beta_n(a'_{n-1})$ by recursion hypothesis {\rm{(iii)}}.
\item  By construction, $a_n=\alpha(c_n)$ and $a'_n=\alpha(c'_n)$ where $c_n=\beta_n^{-1}\alpha\beta_n(a_{n-1})$ and $c'_n=\beta_n^{-1}\alpha\beta_n(a'_{n-1})$. Thus, by recursion hypothesis {\rm{(iii)}}, $a_n=\alpha(c_n)>c_n$ and $a'_n=\alpha(c'_n)<c'_n$, so $\alpha(a_n)>\alpha(c_n)=a_n$ and $\alpha(a'_n)<\alpha(c'_n)=a'_n$. On the other hand, $\beta_{n+1}^{-1}\alpha\beta_{n+1}(a_n)=c_{n+1}$ and $\beta_{n+1}^{-1}\alpha\beta_{n+1}(a'_n)=c'_{n+1}$ where $\alpha(c_{n+1})>c_{n+1}$ and $\alpha(c'_{n+1})<c'_{n+1}$. 
\item By recursion hypotheses {\rm{(iii,iv)}}, $a_n=\alpha(c_n)>c_n=\max X_n\cap \Gamma^{\circ}$ and $a'_n=\alpha(c'_n)<c'_n=\min X_n\cap \Gamma^{\circ}$. By construction, $\alpha(b)>b$ and $\alpha(b')<b'$. Applying $\beta^{-1}_{n+1}$, we get $c_{n+1}>a_n$ and $c'_{n+1}<a'_n$. Also, by construction, $c'_{n+1}<\pi_{\Gamma}(w_n)<c_{n+1}$ and $b'<\pi_{\Gamma}(w_n)<b$. Hence, by \cref{l:generated}(\ref{itm:generated 6}), $c'_{n+1}=\min X_{n+1}\cap \Gamma^{\circ}$ and $c_{n+1}=\max X_{n+1}\cap \Gamma^\circ$.
\item (In case $\Gamma\neq \Gamma^{\circ}$) By construction applying \cref{l:generated}(\ref{itm:generated 8}).
\end{enumerate}
\end{enumerate}

Set $\beta=\bigcup \beta_n$, so $\beta\in\Aut(\Tt^{M}_\Gamma)$ by property {\rm{(i)}}. Set $a=a_0$ and $a'=a'_0$. For any $x\in \Gamma^{\circ}$, $x=w_n$ for some $n\in\N$. Then, $a'_{n+1}=\alpha(c'_{n+1})<c'_{n+1}=\min X_{n+1}\cap \Gamma^{\circ}\leq x\leq \max X_{n+1}\cap\Gamma^{\circ}=c_{n+1}< \alpha(c_{n+1})=a_{n+1}$ by properties {\rm{(i,ii,iii,iv)}}. The same is true for $x=a_n$ or $x=a'_n$. 
\begin{enumerate}[label={\rm{(\alph*)}}, wide]
\item By property {\rm{(ii)}}, $a_k=(\alpha\beta^{-1}\alpha\beta)^k(a)$. The sequence $(a_k)_{k\in\N}$ is strictly increasing and, for any $x\in \Gamma^\circ$, $x\in X_n$ for some $n$, so $x<\max X_n\cap \Gamma^\circ=c_n<\alpha(c_n)=a_n$ by properties {\rm{(iii,iv)}}, concluding that $(a_k)_{k\in\N}$ is upwards unbounded. 
\item By property {\rm{(ii)}}, $a'_k=(\alpha\beta^{-1}\alpha\beta)^k(a')$. The sequence $(a'_k)_{k\in\N}$ is strictly decreasing and, for any $x\in \Gamma^\circ$, $x\in X_n$ for some $n$, so $x>\min X_n\cap \Gamma^\circ=c'_n>\alpha(c'_n)=a'_n$ by properties {\rm{(iii,iv)}}, concluding that $(a'_k)_{k\in\N}$ is downwards unbounded.
\item (In case $\Gamma\neq\Gamma^{\circ}$) For any $v>\gamma\coloneqq\max\Gamma$, $v=v_n$ for some $n\in\N$. Then, $\beta(v)$ is in the same cone as $v$ by property {\rm{(v)}}. As $v$ is arbitrary, we get that $\beta$ setwise fixes every cone above $\gamma$. \qedhere
\end{enumerate}
\end{proof}
\end{lem}

\begin{lem} \label{l:construction 3} Let $\Gamma$ be a semibranch and $\alpha$ an automorphism of $\Tt^M$ fixing $\Gamma$ setwise. Suppose that there are $a,a'\in \Gamma^\circ$ with $(\alpha^n(a))_{n\in\N}$ strictly increasing and upwards unbounded in $\Gamma^{\circ}$ and $(\alpha^n(a'))_{n\in\N}$ strictly decreasing and downwards unbounded in $\Gamma^\circ$. Then, there is an automorphism $\beta$ of $\Tt^M$ fixing $\Gamma$ setwise such that:
\begin{enumerate}[label={\rm{(\alph*)}}, wide]
\item $\alpha^{-1}\beta^{-1}\alpha\beta$ is strictly increasing in $\Gamma^{\circ}$. 
\item If $\Gamma$ is not a branch, $\beta$ setwise fixes every cone above $\max\Gamma$. 
\end{enumerate}
\begin{proof} First of all, note that it follows from the assumptions that $\alpha$ is strictly decreasing in $\Gamma^{\circ}_{<a'}$ and strictly increasing in $\Gamma^{\circ}_{>a}$. In particular, $a>a'$.

Let $W\coloneqq\{w\tq \pi_\Gamma(w)\in \Gamma^{\circ}\}$ and $V\coloneqq\{v\tq \pi_\Gamma(v)\notin \Gamma^{\circ}\}$. Let $(w_n)_{n\in\N}$ be an enumeration of $W$ and $(v_n)_{n\in\N}$ an enumeration of $V$. Without loss of generality, assume $\pi_\Gamma(w_0)>a$.  Let 
\[l(n)\coloneqq\min\{k\in\N\tq \pi_\Gamma(w_k)\leq \alpha^{\lfloor (n-1)/2\rfloor}(a')\mathrm{\ or\ }\pi_\Gamma(w_k)\geq \alpha^{\lfloor (n-1)/2\rfloor}(a)\}.\]
Note that, since $(\alpha^n(a))_{n\in\N}$ is upwards unbounded and $(\alpha^n(a'))_{n\in\N}$ is downwards unbounded, we get that $l(n)\rightarrow\infty$.

By back{\hyp}and{\hyp}forth, we recursively construct a chain of partial finite isomorphisms $(\beta_n)_{n\in\N}$ of $\Tt^M$ setwise fixing $\Gamma$, and an increasing sequence $(a_n)_{n\in\Z}$ of elements in $\Gamma^\circ$ with the following properties:
\begin{enumerate}[label={(\roman*)}, wide]
\item $\beta_n:\ X_n\rightarrow Y_n$ with $X_n=\langle X_n\rangle_{\Gamma}$ and $Y_n=\langle Y_n\rangle_{\Gamma}$, and $w_k,v_k\in X_n$ and $w_k,v_k\in Y_n$ for every $k<l(n)$.
\item For each $-n<k<n$, we have $\alpha\beta_n(a_{k-1})=\beta_n\alpha(a_k)$. In particular, $a_k\in X_n$ for $-n-1<k<n-1$, and $\alpha(a_k)\in X_n$ for $-n<k<n$.  
\item For $n\geq 1$, we have $a_1>a>a'>a_0$ and $\beta_n(a_0)>a>a'>\beta_n(a_{-1})$. Furthermore, we have $a_k>\alpha^{k-1}(a)$ for $0<k<n$; $\beta_n(a_k)>\alpha^{\lfloor k/2\rfloor}(a)$ for $0\leq k<n$; $a_k<\alpha^{\lfloor -k/2\rfloor}(a')$ for $-n\leq k\leq 0$, and $\beta_n(a_k)<\alpha^{-k-1}(a')$ for $-n\leq k<0$.
\item For $n\geq 1$, we have $\alpha(a_{-n+1})=\min X_n\cap \Gamma^{\circ}$ and $c_{n}=\max X_n\cap \Gamma^\circ$ with $c_n>\alpha^2(a_{n-1})$ and $c_n>\alpha(a)$ where $c_n\coloneqq\beta_n^{-1}\alpha\beta_n(a_{n-1})$. In particular, $a_{n-1}\in X_n$ and $\alpha\beta_n(a_{n-1})\in Y_n$. 
\item (In case $\Gamma\neq \Gamma^{\circ}$) $\beta_n(v)$ and $v$ are in the same cone above $\gamma\coloneqq \max\Gamma$ for any $v\in X_n$ with $v>\gamma$. Also, $\gamma\in X_n$, $\gamma\in Y_n$ and $\beta_n(\gamma)=\gamma$.
\end{enumerate}

In stage $n$ we construct $\beta_n$ and $(a_k)_{-n\leq k<n}$.
\begin{enumerate}[wide]
\item[Case $n=0$:] If $\Gamma$ is not a branch, we take $\beta_{0}:\ \gamma\mapsto \gamma$, with $X_{0}=Y_{0}=\{\gamma\}$, where $\gamma\coloneqq\max\Gamma$. If $\Gamma$ is a branch, we take $\beta_0=\emptyset$.
\item[Case $n=1$:] Take $a_0<a'$ arbitrary. Applying \cref{l:saturation of branch intervals} twice, we find $b,b'\in \Gamma^{\circ}$ with $b'<\alpha(a')<a'<a<b$ such that $a_0\alpha(a_0)\equiv^{\Rr\Gamma}bb'$. Applying \cref{l:saturation of branch intervals} again, we find $c_1>\alpha(a)$ such that $a_0\alpha(a_0)c_1\equiv^{\Rr\Gamma}bb'\alpha(b)$. By homogeneity, we find $a_{-1}$ with $a_0\alpha(a_0)c_1a_{-1}\equiv^{\Rr\Gamma}bb'\alpha(b)\alpha^{-1}(b')$. By homogeneity, we extend $\beta_0$ to $\beta'_1:\ X'_1\rightarrow Y'_1$ with $\beta'_1:\ a_0\alpha(a_0)c_1a_{-1}\mapsto bb'\alpha(b)\alpha^{-1}(b')$. Further, we restrict $\beta'_1$ so that $X_1'=\langle X_0\, a_0\alpha(a_0)c_1a_{-1}\rangle_\Gamma$. If $\Gamma=\Gamma^{\circ}$, we set $\beta_1=\beta_1'$. If $\Gamma\neq\Gamma^\circ$, let $\bar{x}'$ enumerate $X'_1$ and $\bar{y}'$ be the enumeration of $Y'_1$ given by $\beta_1':\ \bar{x}'\mapsto \bar{y}'$. By \cref{l:saturation of cones} and homogeneity, we extend $\beta'_1$ to $\beta_1:\ X_1\rightarrow Y_1$ by adding $v_0$ to the domain and the image such that $v_0$ and $\beta_1(v_0)$ are in the same cone. As above, we restrict so that $X_1=\langle X_0\, v_0\beta^{-1}_1(v_0)\rangle_\Gamma$.

Since $\pi_\Gamma(w_0)>a$, we get $l(1)=0$, so $\beta_1$ satisfies {\rm{(i)}} vacuously. If $-1<k<1$, then $k=0$ and $\alpha\beta_1(a_{-1})=b'=\beta_1\alpha(a_0)$, so {\rm{(ii)}} holds. Since $b'<\alpha(a')$, we get $\beta_1(a_{-1})=\alpha^{-1}(b')<a'$, while obviously $a_0<a'<a<b=\beta_1(a_0)$. In particular, $\beta_1(a_{-1})<a'<\beta_1(a_0)$, so $a_{-1}<a_0$ too. Thus, we get {\rm{(iii)}}. By construction, $c_1=\beta_1^{-1}\alpha\beta(a_0)$, and $c_1>\alpha(a)>a>a_0>\alpha^2(a_0)$. Also, $b'<\alpha^{-1}(b')<a'<a<b<\alpha(b)$, getting $\alpha(a_0)<a_{-1}<a_0<c_1$ by applying $\beta^{-1}_1$. Hence, by \cref{l:generated}(\ref{itm:generated 6}), we get that $c_1=\max X_1\cap \Gamma^{\circ}$ and $\alpha(a_0)=\min X_1\cap \Gamma^{\circ}$, so {\rm{(iv)}} holds. Finally, by \cref{l:generated}(\ref{itm:generated 8}), we have $\cone(X_1)=\cone_\gamma(v_0)$, so {\rm{(v)}} holds.
\item[Recursion case (case $n+1$):] Suppose $\beta_n:\ X_n\rightarrow Y_n$ and $(a_k)_{-n\leq k<n}$ are already given where $n\geq 1$. Let $\bar{x}$ enumerate $X_n$ and $\bar{y}$ be the enumeration of $Y_n$ given by $\beta_n:\ \bar{x}\mapsto \bar{y}$. 
\begin{enumerate}[label={\rm{Step {\arabic*}:}}, wide]
\item Set $a_n\coloneqq\alpha^{-1}(c_n)$ where $c_n\coloneqq\beta_n^{-1}\alpha\beta_n(a_{n-1})$. By recursion hypothesis, $a_n\in \Gamma^{\circ}$. By homogeneity, find $b\in\Gamma^\circ$ such that $a_n\bar{x}\equiv^{\Rr\Gamma} b\bar{y}$. 
\item When $n>1$, by recursion hypotheses {\rm{(iii, iv)}}, $a_n=\alpha^{-1}(c_n)>\alpha(a_{n-1})>a_{n-1}>a_0$. When $n=1$, $a_1=\alpha^{-1}(c_1)>a>a_0$ by recursion hypothesis {\rm{(iii,iv)}}. 

Hence, $b>\beta_n(a_{n-1})\geq \beta_n(a_0)>a$ by recursion hypothesis {\rm{(iii)}}, concluding that $\alpha(b)>b$. By recursion hypothesis {\rm{(iv)}}, we also get $\alpha(b)>\alpha\beta_n(a_{n-1})=\beta_n(c_n)=\max Y_n\cap \Gamma^{\circ}$. Thus, by \cref{l:saturation of branch intervals}, find $c_{n+1}\in\Gamma^{\circ}$ such that $c_{n+1}a_n\bar{x}\equiv^{\Rr\Gamma} \alpha(b)b\bar{y}$ with $c_{n+1}>\alpha^2(a_n)$. Note that, since $\alpha(b)>\max Y_n\cap \Gamma^{\circ}$ and $\alpha(b)>b$, we get $c_{n+1}>\max X_n\cap \Gamma^{\circ}$ and $c_{n+1}>a_n$.
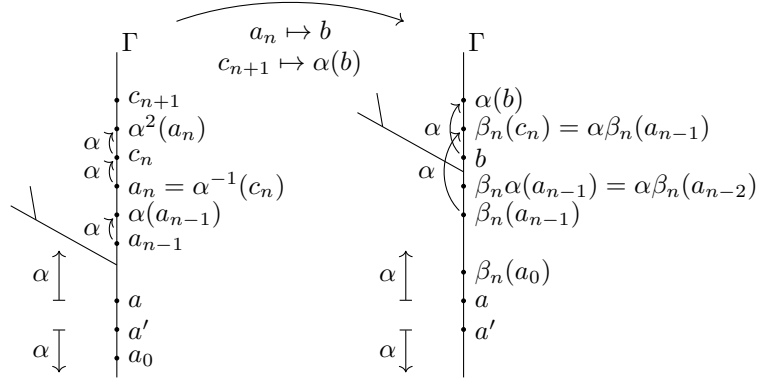
\begin{figure}[ht]
\begin{center}
\begin{tikzpicture}[scale=0.38]
\node (xL) at (-15,-6) {}; \node (yL) at (-15,6) {}; \node (xR) at (15,-6) {}; \node (yR) at (15,6) {};


\node (topL) [label={[xshift=5, yshift=-10]:$\Gamma$}] at (-6,6) {}; \node (botL) at (-6,-6) {}; \node (piwkL) [draw, shape=coordinate] at (-6,-1.75) {}; \node (rL) at (-10,0.5) {}; \node (rmeetwkL) [draw, shape=coordinate] at (-8.8,-0.175) {}; \node (wkL) [draw, shape=coordinate] at (-9,1) {};

\node (topR) [label={[xshift=5, yshift=-10]:$\Gamma$}] at (6,6) {}; \node (botR) at (6,-6) {}; \node (piwkR) [draw, shape=coordinate] at (6,1.5) {}; \node (rR) at (12-10,0.5+3.25) {}; \node (rmeetwkR) [draw, shape=coordinate] at (12-8.8,-0.175+3.25) {}; \node (wkR) [draw, shape=coordinate] at (12-9,1+3.25) {}; 

\node (a0) [draw, shape=circle, fill, scale=0.15, label=right:$a_0$] at (-6,-5) {};
\node (a') [draw, shape=circle, fill, scale=0.15, label=right:$a'$] at (-6,-4) {};
\node (a) [draw, shape=circle, fill, scale=0.15, label=right:$a$] at (-6,-3) {};
\node (an1) [draw, shape=circle, fill, scale=0.15, label=right:$a_{n-1}$] at (-6,-1) {};
\node (alphaan1) [draw, shape=circle, fill, scale=0.15, label=right:$\alpha(a_{n-1})$] at (-6,0) {};
\node (an) [draw, shape=circle, fill, scale=0.15, label={right:$a_n=\alpha^{-1}(c_n)$}] at (-6,1) {};
\node (cn) [draw, shape=circle, fill, scale=0.15, label=right:$c_n$] at (-6,2) {};
\node (alpha2an) [draw, shape=circle, fill, scale=0.15, label=right:$\alpha^2(a_n)$] at (-6,3) {};
\node (cn1) [draw, shape=circle, fill, scale=0.15, label=right:$c_{n+1}$] at (-6,4) {};

\node (a') [draw, shape=circle, fill, scale=0.15, label=right:$a'$] at (6,-4) {};
\node (a) [draw, shape=circle, fill, scale=0.15, label=right:$a$] at (6,-3) {};
\node (betaa0) [draw, shape=circle, fill, scale=0.15, label=right:$\beta_n(a_0)$] at (6,-2) {};
\node (betaan1) [draw, shape=circle, fill, scale=0.15, label={right:$\beta_n(a_{n-1})$}] at (6,0) {};
\node (betaalphaan1) [draw, shape=circle, fill, scale=0.15, label={right:$\beta_n\alpha(a_{n-1})=\alpha\beta_n(a_{n-2})$}] at (6,1) {};
\node (b) [draw, shape=circle, fill, scale=0.15, label=right:$b$] at (6,2) {};
\node (betacn) [draw, shape=circle, fill, scale=0.15, label={right:$\beta_n(c_n)=\alpha\beta_n(a_{n-1})$}] at (6,3) {};
\node (alphab) [draw, shape=circle, fill, scale=0.15, label=right:$\alpha(b)$] at (6,4) {};

\draw (botL) -- (topL); \draw (piwkL) -- (rL); \draw (rmeetwkL) -- (wkL);
\draw (botR) -- (topR); \draw (piwkR) -- (rR); \draw (rmeetwkR) -- (wkR);
\draw[->, shorten <=20, shorten >=20] (topL) to [bend left=20] node[below]{$\begin{array}{c}a_n\mapsto b\\ c_{n+1}\mapsto \alpha(b)\end{array}$} (topR);

\draw[|->] (-8,-3) to node[auto]{$\alpha$} (-8,-1.25);
\draw[|->] (-8,-4) to node[left]{$\alpha$} (-8,-5.5);
\draw[|->] (4,-3) to node[auto]{$\alpha$} (4,-1.25);
\draw[|->] (4,-4) to node[left]{$\alpha$} (4,-5.5);

\draw[->, shorten <=1.5, shorten >=1.2] (an1) to [out=180-45, in=180+45] node[auto]{$\alpha$} (alphaan1) ;
\draw[->, shorten <=1.5, shorten >=1.2] (an) to [out=180-45, in=180+45] node[auto]{$\alpha$} (cn) ;
\draw[->, shorten <=1.5, shorten >=1.2] (cn) to [out=180-45, in=180+45] node[auto]{$\alpha$} (alpha2an) ;

\draw[->, shorten <=1.5, shorten >=1.2] (betaan1) to [out=180-45, in=180+45] node[auto]{$\alpha$} (betacn) ;
\draw[->, shorten <=1.5, shorten >=1.2] (b) to [out=180-45, in=180+45] node[auto]{$\alpha$} (alphab) ;
\end{tikzpicture}
\caption{Illustration of steps 1 and 2 when $n>1$.}
\end{center}
\end{figure}

\item Since the sequence is increasing, $a_{-n}<a_{-n+1}\leq a_0$ and by recursion hypothesis {\rm{(iii)}}, $a_0<a'$. Thus, $\alpha(a_{-n})<\alpha(a_{-n+1})=\min X_n\cap \Gamma^{\circ}$ by recursion hypothesis {\rm{(iv)}}. By recursion hypothesis {\rm{(ii)}}, $a_{-n}\in X_n$, so $\alpha(a_{-n+1})<a_{-n}$. Hence, by \cref{l:saturation of branch intervals} again, find $b'\in\Gamma^{\circ}$ such that $\alpha(a_{-n})c_{n+1}a_n\bar{x}\equiv^{\Rr\Gamma}b'\alpha(b)b\bar{y}$ with $b'<\alpha\beta_n\alpha(a_{-n+1})$. Since $\alpha(a_{-n+1})<a_{-n}\leq a_{-1}$, we get that $\beta_n\alpha(a_{-n+1})<\beta_n(a_{-n})\leq \beta_n(a_{-1})$. By recursion hypothesis {\rm{(iii)}}, $\beta_n(a_{-1})<a'$, concluding that $\alpha\beta_n\alpha(a_{-n+1})<\beta_n\alpha(a_{-n+1})<a'$.
\item By homogeneity, we find $a_{-n-1}\in\Gamma^{\circ}$ such that $a_{-n-1}\alpha(a_{-n})c_{n+1}a_n\bar{x}$ $\equiv^{\Rr\Gamma}\alpha^{-1}(b')b'\alpha(b)b\bar{y}$. 
\begin{figure}[ht]
\begin{center}
\begin{tikzpicture}[scale=0.38]
\node (xL) at (-15+0.4,-6) {}; \node (yL) at (-15+0.4,6) {}; \node (xR) at (15,-6) {}; \node (yR) at (15,6) {};


\node (topL) [label={[xshift=5, yshift=-10]:$\Gamma$}] at (-6,6) {}; \node (botL) at (-6,-6) {}; 
\node (piwkL) [draw, shape=coordinate] at (-6,-2.5) {}; \node (rL) at (-10,-0.5) {}; \node (rmeetwkL) [draw, shape=coordinate] at (-8.8,-1.1) {}; \node (wkL) [draw, shape=coordinate] at (-9,0) {};

\node (topR) [label={[xshift=5, yshift=-10]:$\Gamma$}] at (6,6) {}; \node (botR) at (6,-6) {}; 
\node (piwkR) [draw, shape=coordinate] at (6,0.5) {}; \node (rR) at (12-10,-0.5+3) {}; \node (rmeetwkR) [draw, shape=coordinate] at (12-8.8,-1.1+3) {}; \node (wkR) [draw, shape=coordinate] at (12-9,3) {}; 

\node (a) [draw, shape=circle, fill, scale=0.15, label=right:$a$] at (-6,4) {};
\node (a') [draw, shape=circle, fill, scale=0.15, label=right:$a'$] at (-6,3) {};
\node (a0) [draw, shape=circle, fill, scale=0.15, label=right:$a_0$] at (-6,2) {};
\node (a1n) [draw, shape=circle, fill, scale=0.15, label=right:$a_{-n+1}$] at (-6,0) {};
\node (an) [draw, shape=circle, fill, scale=0.15, label=right:$a_{-n}$] at (-6,-1) {};
\node (alphaa1n) [draw, shape=circle, fill, scale=0.15, label={right:$\alpha(a_{-n+1})$}] at (-6,-2) {};
\node (an1) [draw, shape=circle, fill, scale=0.15, label=right:$a_{-n-1}$] at (-6,-3) {};
\node (alphaan) [draw, shape=circle, fill, scale=0.15, label=right:$\alpha(a_{-n})$] at (-6,-4) {};

\node (a) [draw, shape=circle, fill, scale=0.15, label=right:$a$] at (6,4) {};
\node (a') [draw, shape=circle, fill, scale=0.15, label=right:$a'$] at (6,3) {};
\node (betaalphaa1n) [draw, shape=circle, fill, scale=0.15, label={right:$\beta_n\alpha(a_{-n+1})=\alpha\beta_{n}(a_{-n})$}] at (6,0) {};
\node (alphab) [draw, shape=circle, fill, scale=0.15, label={right:$\alpha^{-1}(b')$}] at (6,-1) {};
\node (alphabetaalphaa1n) [draw, shape=circle, fill, scale=0.15, label=right:$\alpha\beta_n\alpha(a_{-n+1})$] at (6,-2) {};
\node (b) [draw, shape=circle, fill, scale=0.15, label={right:$b'$}] at (6,-3) {};

\draw (botL) -- (topL); 
\draw (piwkL) -- (rL); \draw (rmeetwkL) -- (wkL);
\draw (botR) -- (topR); 
\draw (piwkR) -- (rR); \draw (rmeetwkR) -- (wkR);
\draw[->, shorten <=20, shorten >=20] (topL) to [bend left=20] node[below]{$\begin{array}{c}\alpha(a_{-n})\mapsto b'\\ a_{-n-1}\mapsto \alpha^{-1}(b')\end{array}$} (topR);

\draw[|->] (-7,4) to node[auto]{$\alpha$} (-7,5.5);
\draw[|->] (-7,3) to node[left]{$\alpha$} (-7,1.25);
\draw[|->] (8,4) to node[right]{$\alpha$} (8,5.5);
\draw[|->] (8,3) to node[auto]{$\alpha$} (8,1.25);

\draw[->, shorten <=1.5, shorten >=1.2] (a1n) to [out=180+45, in=180-45] node[left]{$\alpha$} (alphaa1n) ;
\draw[->, shorten <=1.5, shorten >=1.2] (an) to [out=180+45, in=180-45] node[left]{$\alpha$} (alphaan) ;

\draw[->, shorten <=1.5, shorten >=1.2] (betaalphaa1n) to [out=180+45, in=180-45] node[left]{$\alpha$} (alphabetaalphaa1n) ;
\draw[->, shorten <=1.5, shorten >=1.2] (alphab) to [out=180+45, in=180-45] node[left]{$\alpha$} (b) ;
\end{tikzpicture}
\caption{Illustration of steps 3 and 4.}
\subcaption*{(The position of $\alpha^{-1}(b')$ respect to $\alpha\beta_n\alpha(a_{-n+1})$ is not precise --- there is no reason for it to be greater than $\alpha \beta_n \alpha (a_{-n+1})$)}
\end{center}
\end{figure}
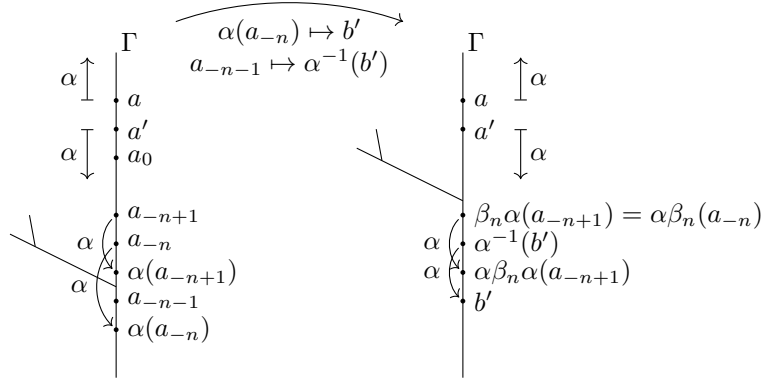
\item By homogeneity, we extend $\beta_n$ to a partial isomorphism $\beta'_{n+1}:\ X'_{n+1}\rightarrow Y'_{n+1}$ with $\beta'_{n+1}:\ a_{-n-1}\alpha(a_{-n})c_{n+1}a_n\mapsto\alpha^{-1}(b')b'\alpha(b)b$ and $w_k\in X'_{n+1}\cap Y'_{n+1}$ for $k<l(n+1)$. Further, we restrict $\beta'_{n+1}$ so that $X'_{n+1}=\langle X_n\, a_{-n-1}\alpha(a_{-n})c_{n+1}a_n\,$ $\{w_k{\beta'}^{-1}_{n+1}(w_k)\}_{k<l(n+1)}\rangle_\Gamma$. In case $\Gamma=\Gamma^{\circ}$, we conclude setting $\beta_{n+1}=\beta'_{n+1}$. Otherwise, continue to step 6.
\item (In case $\Gamma\neq \Gamma^{\circ}$) By recursion hypothesis {\rm{(vi)}}, \cref{l:saturation of cones} and homogeneity, we can extend $\beta'_{n+1}$ to a partial isomorphism $\beta_{n+1}$ by adding $v_k$, for $k<l(n+1)$, to the domain and the image such that $v_k$ and $\beta_{n+1}(v_k)$ are in the same cone for each $k<l(n+1)$.  As above, we restrict so that $X_{n+1}=\langle X'_{n+1}\, \{v_k \beta_{n+1}^{-1}(v_k)\}_{k<l(n+1)}\rangle_\Gamma$.
\end{enumerate}

Now, we check that $\beta_{n+1}$ and $(a_k)_{-n-1\leq k<n+1}$ just constructed satisfy the recursion hypotheses:
\begin{enumerate}[label={(\roman*)}, wide]
\item By construction, $\beta_{n+1}:\ X_{n+1}\rightarrow Y_{n+1}$ is a partial isomorphism extending $\beta_n$ with $w_k,v_k\in X_{n+1}$ and $w_k,v_k\in Y_{n+1}$ for $k<l(n+1)$.
\item By construction, $a_{-n-1},\alpha(a_{-n}),a_n\in X_{n+1}$, $a_{-n}\in X_n\subseteq X_{n+1}$ and $\alpha(a_n)=c_n\in X_{n}\subseteq X_{n+1}$. Also, by construction, $\beta_{n+1}\alpha(a_n)=\beta_{n+1}(c_n)=\alpha\beta_n(a_{n-1})$. Finally, we get that $\beta_{n+1}(a_{-n-1})=\alpha^{-1}(b')$, so $\alpha\beta_{n+1}(a_{-n-1})=\beta_{n+1}\alpha(a_{-n})$.
\item By recursion hypothesis {\rm{(iii)}}, we have that $a'>a_0$ and $\beta_n(a_0)>a>a'>\beta_n(a_{-1})$. 

When $n=1$, by recursion hypothesis {\rm{(iv)}}, we have that $a_1=\alpha^{-1}(c_1)>a>a_0$. For $n>1$, by construction and recursion hypotheses {\rm{(iii,iv)}}, $a_n=\alpha^{-1}(c_n)>\alpha(a_{n-1})>a_{n-1}$. On the other hand, by construction, $\beta_{n+1}(a_{-n-1})=\alpha^{-1}(b')<\beta_n\alpha(a_{-n+1})$, concluding that $a_{-n-1}<\alpha(a_{-n+1})$. In particular, as $\alpha(a_{-n+1})=\min X_n\cap\Gamma^{\circ}\leq a_{-n}$ by recursion hypothesis {\rm(iv)}, we get $a_{-n-1}<a_{-n}$. Thus, $(a_k)_{-n-1\leq k<n+1}$ is strictly increasing with $a_1>a>a'>a_0$.

As we have already noted, $a_1>a$ and $a_n>\alpha(a_{n-1})$ for $n>1$, concluding $a_n>\alpha^{n-1}(a)$ by recursion hypothesis {\rm{(iii)}}. Similarly, $a_{-1}<a_0<a'$ and $a_{-n-1}<\alpha(a_{-n+1})$, concluding $a_{-n-1}<\alpha^{\lfloor (n+1)/2\rfloor}(a')$ by recursion hypothesis {\rm{(iii)}}.

For $n=1$, $\beta_2(a_1)>\beta_2(a_0)>a$. For $n>1$, by recursion hypothesis {\rm{(iv)}}, $a_n=\alpha^{-1}(c_n)>\alpha(a_{n-1})$, so $\beta_{n+1}(a_n)>\beta_{n+1}\alpha(a_{n-1})$. Thus, by recursion hypotheses {\rm{(ii,iii)}}, $\beta_{n+1}(a_n)>\beta_n\alpha(a_{n-1})=\alpha\beta_n(a_{n-2})>\alpha\circ\alpha^{\lfloor (n-2)/2\rfloor}(a)=\alpha^{\lfloor n/2\rfloor}(a)$. 

Finally, by construction, $b'<\alpha\beta_n\alpha(a_{-n+1})$ where $\beta_{n+1}(a_{-n-1})=\alpha^{-1}(b')$. Thus, $\beta_{n+1}(a_{-n-1})<\beta_n\alpha(a_{-n+1})$, getting $\beta_{n+1}(a_{-n-1})<\alpha\beta_n(a_{-n})$ by recursion hypothesis {\rm{(ii)}}. Hence, $\beta_{n+1}(a_{-n-1})<\alpha\circ\alpha^{n-1}(a')=\alpha^{n}(a')$ by recursion hypothesis {\rm{(iii)}}.
\item By construction, $c_{n+1}>\max X_n\cap \Gamma^{\circ}$ and $c_{n+1}>a_n$, and $\alpha(a_{-n})<\min X_n\cap\Gamma^{\circ}$. Also, $b'<a'$, so $b'<\alpha^{-1}(b')$, getting $\alpha(a_{-n})<a_{-n-1}$ by applying $\beta_{n+1}$. Now, by property {\rm{(iii)}}, $c_{n+1}>a_n>\alpha^{\lfloor n/2\rfloor}(a)>\pi_{\Gamma}(w_k)>\alpha^{\lfloor n/2 \rfloor}(a')>a_{-n-1}>\alpha(a_{-n})$ and $\beta_{n+1}(c_{n+1})>\beta_{n+1}(a_n)>\alpha^{\lfloor n/2\rfloor}(a)>\pi_{\Gamma}(w_k)>\alpha^{\lfloor n/2 \rfloor}(a')>\beta_{n+1}(a_{-n-1})>\alpha(a_{-n})$  for any $k<l(n+1)$. Hence, $c_{n+1}=\max X_{n+1}\cap \Gamma^{\circ}$ and $\alpha(a_{-n})=\min X_{n+1}\cap \Gamma^{\circ}$.
\item (In case $\Gamma\neq \Gamma^{\circ}$) By construction.
\end{enumerate}
\end{enumerate}

Set $\beta=\bigcup \beta_n$, so $\beta\in\Aut(\Tt^{M}_\Gamma)$ by {\rm{(i)}} noting that $l(n)\rightarrow \infty$. 
\begin{enumerate}[label={\rm{(\alph*)}}, wide]
\item For any $x\in \Gamma^{\circ}$, $x=w_k$ for some $k\in\N$ and $k<l(n)$ for some $n$. Thus, $a_{-n}<x< a_n$ by property {\rm{(iii)}}. Take $m\coloneqq\min\{n\in\Z\tq x<a_n\}$. Then, $a_{m-1}\leq x< a_m$, so $x<a_m=\alpha^{-1}\beta^{-1}\alpha\beta(a_{k-1})\leq \alpha^{-1}\beta^{-1}\alpha\beta(x)$. As $x$ is arbitrary, we conclude that $\alpha^{-1}\beta^{-1}\alpha\beta$ is strictly increasing in $\Gamma^{\circ}$. 
\item (In case $\Gamma\neq \Gamma^{\circ}$) For any $v>\gamma\coloneqq\max \Gamma$, $v=v_n$ for some $n\in\N$. Then, $\beta(v)$ is in the same cone as $v$. As $v$ is arbitrary, we get that $\beta$ leaves invariant every cone above $\gamma$.\qedhere
\end{enumerate}
\end{proof}
\end{lem}

\begin{coro}\label{c:semibranch strictly increasing} Let $\Gamma$ be a semibranch and $\alpha$ an automorphism of $\Tt^M$ setwise fixing $\Gamma$. Suppose that $\alpha$ is unbounded in $\Gamma^{\circ}$. Then, there are $4$ conjugates of $\alpha$ or $\alpha^{-1}$ by automorphisms $\beta_1,\ldots,\beta_4$ setwise fixing $\Gamma$ whose product is strictly increasing in $\Gamma^{\circ}$. Furthermore, if $\Gamma$ is not a branch, we can take $\beta_1,\ldots,\beta_4$ setwise fixing every cone above $\gamma=\max\Gamma$. 
\begin{proof} Up to replacing $\alpha$ by $\alpha^{-1}$, we can assume that $\alpha$ is upwards unboundedly increasing. Either $\alpha$ is unboundedly increasing or it is downwards unboundedly decreasing. In the first case, by \cref{l:construction 1}, there is $\beta$ such that $\alpha\cdot \alpha^{\beta}$ is strictly increasing in $\Gamma^{\circ}$. Then, $\alpha\cdot \alpha^{\beta}\cdot \alpha\cdot \alpha^{\beta}$ is strictly increasing in $\Gamma^{\circ}$, so we conclude by taking $\beta_1=\beta_3=\id$ and $\beta_2=\beta_4=\beta$. In the second case, by \cref{l:construction 2,l:construction 3}, there are $\beta,\beta'$ such that $(\alpha\cdot \alpha^{\beta})^{-1}\cdot (\alpha\cdot \alpha^{\beta})^{\beta'}$ is strictly increasing in $\Gamma^{\circ}$. Then, we conclude by taking $\beta_1=\beta$, $\beta_2=\id$, $\beta_3=\beta'$ and $\beta_4=\beta\beta'$.
Furthermore, if $\Gamma$ is not a branch, according to \cref{l:construction 1,l:construction 2,l:construction 3}, we can take $\beta_1,\ldots,\beta_4$ setwise fixing every cone above $\gamma=\max\Gamma$. 
\end{proof}
\end{coro}
\begin{coro} \label{c:fan strictly increasing} Let $\alpha$ be an automorphism of $\Tt^{M}$ permuting infinitely many cones above $\gamma$ and unbounded below $\gamma$. Then, there are $32$ conjugates of $\alpha$ or $\alpha^{-1}$ whose product is an infinite fan above $\gamma$ and strictly increasing below $\gamma$.
\begin{proof} By \cref{l:infinite fan}, there are $8$ conjugates of $\alpha$ or $\alpha^{-1}$ whose product $\alpha'$ is an infinite fan above $\gamma$. By \cref{c:semibranch strictly increasing}, there are $4$ conjugates of $\alpha'$ or ${\alpha'}^{-1}$ by elements setwise fixing every cone above $\gamma$ whose product $\alpha''$ is strictly increasing below $\gamma$. Thus, $\alpha''$ is the product of $32$ conjugates of $\alpha$ or $\alpha^{-1}$. Also, $\alpha''$ is strictly increasing below $\gamma$ and permutes the cones above $\gamma$ as ${\alpha'}^4$, so it is an infinite fan above $\gamma$. 
\end{proof}
\end{coro}


\section{Main results} \label{s:section 5}
Throughout this section, unless otherwise stated, we work in the generic meet{\hyp}tree expansion $\Tt^M$ of an infinite strong Fra\"{\i}ss\'{e} limit $M$ over a finite relational language $\lang_\Rr$ disjoint to $\lang_{\Gamma}$ and $\lang_\gamma$ that has a free weak independence relation $\indepe[\Rr]{}$ (see \cref{d:free weak independence}). As a prototypical example, one may consider the generic meet{\hyp}tree expansion of the random graph, which we call the Rado meet{\hyp}tree.\medskip

In this section we prove our main theorem. The proof is fundamentally based on \cite[Corollary 2.16]{li2019automorphism}. The first step is the following \cref{l:alpha ordering}. Recall that $\sim_\gamma$ denotes the cone equivalence relation above $\gamma$.
\begin{lem} \label{l:alpha ordering} Let $U$ be the universe of an elementary substructure of $\Tt^{M}$. Let $\gamma$ be a point and $a',b,c,d$ finite tuples in $U$ with $c\subseteq d$. Let $\alpha$ be an automorphism of $\Tt^M$. Assume that $\alpha$ is an infinite fan above $\gamma$ and $\alpha(U)=U$.
\begin{enumerate}[label={\rm{(\arabic*)}}, ref={\rm{\arabic*}}, wide]
\item \label{itm:alpha ordering 1} If $\alpha$ is strictly increasing below $\gamma$, then there is a tuple  $a$ in $U$ with $a\equiv^{\Rr\gamma}_ca'$ such that $a\indepe[\Rr\gamma]{c}b$ and $a\indepe[\gamma]{d}\alpha(a)$.
\item \label{itm:alpha ordering 2} If $\alpha$ is strictly decreasing below $\gamma$, then there is a tuple $a$ in $U$ with $a\equiv^{\Rr\gamma}_ca'$ such that $b\indepe[\Rr\gamma]{c}a$ and $\alpha(a)\indepe[\gamma]{d}a$.
\end{enumerate}
\begin{proof} We start with the following claim:

\begin{claim} \label{cl:alpha ordering} For any $x,y$ finite tuples of $U$ and any $x'$ tuple in $\Tt^M$ such that $x'\equiv^{\gamma} x$, there is $x''$ tuple in $U$ such that $x''\equiv^{\Rr\gamma}_yx'$.
\begin{proof} Set $\Gamma\coloneqq U_{\leq \gamma}$. Suppose first that $\Gamma$ is not a branch in $U$. We claim that there is $w$ in $U$ with $w>\gamma$. Indeed, as $\Gamma$ is not a branch, there is $w$ in $U$ with $w> v$ for all $v\in \Gamma$. Suppose aiming a contradiction that $w\meet \gamma<\gamma$. Then, in case {\rm{(1)}}, $w\meet \gamma<\alpha(w)\meet \gamma$. By \cref{l:three elements}, we get $v\coloneqq w\meet \gamma= w\meet \alpha(w)\in\Gamma$. However, $\alpha(v)\in \Gamma$ and $\alpha(v)>v$, concluding $w\not>\alpha(v)$
, a contradiction with the choice of $w$. In case {\rm{(2)}}, $w\meet\gamma<\alpha^{-1}(w)\meet \gamma$. By \cref{l:three elements}, we get $v\coloneqq w\meet \gamma=w\meet\alpha^{-1}(w)\in\Gamma$. However, $\alpha^{-1}(v)\in\Gamma$ and $w\not>\alpha^{-1}(v)$, getting a contradiction. 

Take $w$ in $U$ with $w>\gamma$. Then, $w,\alpha(w)$ are in $U$ and, as $\alpha$ is a fan, $\gamma=w\meet\alpha(w)$, concluding that $\gamma$ is in $U$. Then, by $\omega${\hyp}categoricity, $(U,\gamma)\cong\Tt^M_\gamma$ and, thus, by saturation, we can find $x''$ in $U$ realising the type of $x'$ over $y$.\smallskip

Suppose now that $\Gamma$ is a branch of $U$. By $\omega${\hyp}categoricity of $\Tt^M$, we get that $(U,\Gamma)$ is isomorphic to $(\Tt^M,\tilde{\Gamma})$ with $\tilde{\Gamma}$ a branch of $\Tt^M$. Thus, by \cref{l:generic meet-tree expansions}, $(U,\Gamma)$ is the generic branched meet{\hyp}tree expansion of $U_{\mid\lang_\Rr}\cong M$.  Also, consider the semibranch $\Gamma'\coloneqq\{v\in \Tt^M\tq v\leq \gamma\}$ of $\Tt^M$ and note that $(U,\Gamma)$ is a substructure of $(\Tt^M,\Gamma')$. Indeed, obviously $\Gamma=\Gamma' \cap U$. On the other hand, in both cases {\rm{(1)}} and {\rm{(2)}}, for any $w$ in $U$, as $w\meet \gamma<\gamma$, we have $w\meet \gamma\neq \alpha(w)\meet \gamma$. Hence, by \cref{l:three elements}, $\min\{w\meet \gamma,\alpha(w)\meet \gamma\}=w\meet \alpha(w)\in U$, concluding that $w\meet \gamma\in U$. Thus, $\pi_\Gamma(w)=\pi_{\Gamma'}(w)=w\meet \gamma$. 

Let $T=\langle x',y\rangle_{\Gamma'}$ be the generated substructure in $(\Tt^M,\Gamma')$. Pick a tuple $\bar{x}'$ extending $x'$ and enumerating $T$ and a tuple $\bar{y}$ extending $y$ and enumerating $\langle y\rangle_{\Gamma'}$. By universality of $(U,\Gamma)$, we find an embedding $f:\ T\rightarrow U$ pointwise fixing $\langle y\rangle_{\Gamma}=\langle y\rangle_{\Gamma'}$. Write $\bar{x}''\coloneqq f(\bar{x}')$. Then, as $f$ and $\id:\ (U,\Gamma)\rightarrow (\Tt^M,\Gamma')$ are embeddings, we have that $\qftp_{\Rr}(\bar{x}''/\bar{y})=\qftp_{\Rr}(\bar{x}'/\bar{y})$ and $\qftp_{\Gamma'}(\bar{x}''/\bar{y})=\qftp_{\Gamma'}(\bar{x}'/\bar{y})$. Since $x'\equiv^{\gamma}x$, we get that $\pi_{\Gamma'}(x'_i)=x'_i\meet\gamma<\gamma$ for each $i<n$ where $\bar{x}'=(x'_i)_{i<n}$. Also, as $\bar{x}''$ is in $U$, $\pi_{\Gamma'}(x''_i)=\pi_\Gamma(x''_i)<\gamma$ for each $i<n$ where $\bar{x}''=(x''_i)_{i<n}$. Therefore, by \cref{c:quantifier free type tree}(\ref{itm:quantifier free type tree 3}), we get that $\qftp_\gamma(\bar{x}''/\bar{y})=\qftp_\gamma(\bar{x}'/\bar{y})$. Together with $\qftp_\Rr(\bar{x}''/\bar{y})=\qftp_\Rr(\bar{x}'/\bar{y})$, we get $\qftp_{\Rr\gamma}(\bar{x}''/\bar{y})=\qftp_{\Rr\gamma}(\bar{x}'/\bar{y})$, and so $x''\equiv^{\Rr\gamma}_yx'$ by quantifier elimination of $\Tt^M_\gamma$.\claimqed
\end{proof}
\end{claim}

\begin{enumerate}[label={\rm{(\arabic*)}}, wide]  
\item Let $x=(x_i)_{i\leq n}$ be a finite tuple without $\gamma$ (i.e. $x_i\neq \gamma$ for $i\leq n$) and $n_0\coloneqq \min\{i\leq n\tq x_i>\gamma\}$. We say that $x$ is \emph{$\alpha${\hyp}ordered} if 
\begin{enumerate}[label={\rm{(\roman*)}}, wide]
\item $x_i\meet \gamma\leq x_{i+1}\meet \gamma$ for each $i<n$,
\item if $x_i\sim_\gamma x_j$ with $n_0\leq i<j\leq n$, then $x_{i+1}\sim_\gamma x_j$,
\item $x_i\meet \alpha(x_j)\leq \gamma$ for each $i\leq j\leq n$ with $j\geq n_0$.
\end{enumerate} 
Let's remark three properties of $\alpha${\hyp}orderings:
\begin{enumerate}[label={$\bullet$}, wide]
\item By definition, any subtuple of an $\alpha${\hyp}ordered tuple is $\alpha${\hyp}ordered. 
\item Since $\alpha$ is an infinite fan, up to a permutation, every finite tuple (without $\gamma$) can be $\alpha${\hyp}ordered. Indeed, we can easily permute $x$ to guarantee that $x_i\meet\gamma\leq x_{i+1}\meet \gamma$ for each $i<n$. Once we get the first condition, for $i\geq n_0$, this condition imposes nothing, so we can freely permute the elements $x_{n_0},\ldots,x_n$ without losing it. In particular, we can permute $x_{n_0}\ldots x_n$ to be sure that any set of indexes of elements belonging to the same cone is convex; this is the second condition. Finally, it remains to get the third condition. Since $\alpha$ is an infinite fan, no finite set of cones above $\gamma$ is $\alpha${\hyp}invariant. For instance, the set of cones $\mathcal{C}$  of elements of $x$ is not $\alpha${\hyp}invariant, meaning that there is a cone $C\in\mathcal{C}$ such that $\alpha(C)\notin\mathcal{C}$. Up to permuting the consecutive blocks of elements of the same cones, we can assume that $C=\cone(x_n)$. Thus, continuing this way, we can order the cones so that $\alpha(\cone_\gamma(x_j))\notin \{\cone_{\gamma}(x_i)\}_{n_0\leq i\leq j}$ for each $n_0\leq j\leq n$; this is the third condition. 
\item Given an $\alpha${\hyp}ordered tuple $x$, up to adding finitely many elements to $x$, we can assume that $\langle x_{\leq i}\rangle_\gamma =\langle x_{<i}\rangle_\gamma\cup \{x_i\}$ for each $i\leq n$ without losing the $\alpha${\hyp}ordering property. Indeed, suppose $x=(x_i)_{i\leq n}$ is $\alpha${\hyp}ordered and the above is true for $i<n$. Let $e=\max\{x_n\meet \gamma, x_n\meet x_0,\ldots,x_n\meet x_{n-1}\}$. By \cref{l:generated}(\ref{itm:generated 4}), $\langle x_{\leq n}\rangle_\gamma =\langle x_{<n}\rangle_\gamma\cup \{e,x_n\}$. If $e=\gamma$, we are already done. Otherwise, we claim that $x_{<n}\, e\, x_n$ is still $\alpha${\hyp}ordered. Indeed, for any $i\leq n$, $x_i\meet \gamma\leq x_n\meet \gamma\leq e\leq x_n$, so $x_i\meet \gamma \leq e\meet \gamma\leq x_n\meet\gamma$. Also, if $x_n>\gamma$, then $x_n\geq e\geq \gamma$ and, as $e\neq \gamma$, $x_n\sim_\gamma e$. Then, $\alpha(e)\sim_\gamma\alpha(x_n)\not\sim_\gamma x_i$ for $i\leq n$ and $\alpha(e)\sim_\gamma\alpha(x_n)\not\sim_\gamma x_n\sim_\gamma e$.
\end{enumerate}

Without loss of generality, we assume that $a'$ does not contain $\gamma$. Furthermore, applying the last two remarks, we assume that $a'$ is $\alpha${\hyp}ordered and $\langle a'_{\leq i}\rangle_\gamma =\langle a'_{<i}\rangle_\gamma\cup \{a'_i\}$ for each $i\leq n$. We prove by induction on $n$ that there is a tuple $a$ in $U$ with $a\equiv^{\Rr\gamma}_ca'$ such that:
\begin{enumerate}[label={\rm{(\alph*)}}, wide]
\item $a_i\indepe[\Rr\gamma]{ca_{<i}}bd\alpha(a_{<i})$.
\item $a\indepe[\Rr\gamma]{c}b$.
\item $a$ is $\alpha${\hyp}ordered.
\item If $a_i>\gamma$ and $a_i\notin \cone_\gamma(c a_{<i})$, then $\alpha(a_i),\alpha^{-1}(a_i)\notin\cone_\gamma(c a_{<i})$.
\item $a\indepe[\gamma]{d}\alpha(a)$.
\end{enumerate}
\begin{enumerate}[wide]
\item[Case $n=1$.] Applying full existence and invariance of $\indepe[\Rr\gamma]{}$, find $a_0$ with $a_0\equiv^{\Rr\gamma}_ca'_0$ such that $a_0\indepe[\Rr\gamma]{c}bd$. Applying \cref{cl:alpha ordering}, we can take $a_0$ in $U$. This gives us point {\rm{(a)}} and, by monotonicity, point {\rm{(b)}}. Since $\alpha$ is a fan, we trivially have that $a_0$ is $\alpha${\hyp}ordered, so we also have point {\rm{(c)}}.

If $a_0>\gamma$ and $a_0\notin \cone_\gamma(c)$, by \cref{l:cones and cone independence}, $a_0\notin \cone_\gamma(bdc)$. Then, by \cref{l:saturation of cones} and \cref{cl:alpha ordering}, we can take $a_0$ in $U$ with the same $\lang_{\Rr\gamma}${\hyp}type over $b,d,c$ in any cone outside $\cone_\gamma(bdc)$. In particular, we can take $a_0$ such that $\alpha(a_0),\alpha^{-1}(a_0)\notin \cone_\gamma(c)$, getting point {\rm{(d)}}.

It remains to check point {\rm{(e)}}. If $a_0\meet\gamma<\gamma$, as $\alpha$ is strictly increasing below $\gamma$, we get that $a_0\meet \gamma<\alpha(a_0\meet \gamma)=\alpha(a_0)\meet \gamma$. If $a_0\geq \gamma$, as $\alpha$ is a fan, $a_0\meet \alpha(a_0)=\gamma$. Hence, by \cref{l:easy independence}, we get $a_0\indepe[\gamma]{d}\alpha(a_0)$. \smallskip
\item[Inductive case.] Suppose it holds for $n$. Then, by induction hypothesis, we can find $a^*\coloneqq a_{<n}$ tuple in $U$ with $a^*\equiv^{\Rr\gamma}_ca'_{<n}$ satisfying points {\rm{(a--e)}}. Applying full existence of $\indepe[\Rr\gamma]{}$ and \cref{cl:alpha ordering}, we can find $\widetilde{a}_n$ in $U$ such that $a^*\widetilde{a}_n\equiv^{\Rr\gamma}_ca'_{<n}a'_n$ and 

\[\widetilde{a}_n\indepe[\Rr\gamma]{a^*c}bd\alpha(a^*).\] 

By monotony, normality and transitivity of $\indepe[\Rr\gamma]{}$, we get $a^*\widetilde{a}_n\indepe[\Rr\gamma]{c}b$. Hence, we have {\rm{(a,b)}} for the sequence $a^*\widetilde{a}_n$.  

Now, we claim that there is $a_n$ in $U$ with $a_n\equiv^{\Rr\gamma}_{bcda^*\alpha(a^*)}\widetilde{a}_n$ such that the tuple $a^*a_n$ satisfies {\rm{(c,d)}}. Since $a'$ is $\alpha${\hyp}ordered, $a'_i\meet\gamma\leq a'_n\meet \gamma$, so $a_i\meet \gamma\leq \widetilde{a}_n\meet\gamma$ by $a^*\widetilde{a}_n\equiv^{\Rr\gamma}_ca'_{<n}a'_n$. 

If $a'_n\meet \gamma<\gamma$, then $\widetilde{a}_n\meet\gamma<\gamma$, concluding that $a^*a_n$ is $\alpha${\hyp}ordered (i.e. satisfies {\rm{(c)}}) and vacuously satisfies {\rm{(d)}} with $a_n=\widetilde{a}_n$. If $a'_n>\gamma$, also $\widetilde{a}_n>\gamma$. There are two cases: 
\begin{enumerate}[wide]
\item[Case $\widetilde{a}_n\notin\cone_{\gamma}(b c d a^* \alpha(a^*))$.] By \cref{l:saturation of cones} and \cref{cl:alpha ordering}, we can find $a_n$ from $U$ in any cone disjoint to $\cone_{\gamma}(b c d a^* \alpha(a^*))$ with $a_n\equiv^{\Rr\gamma}_{bcda^*\alpha(a^*)}\widetilde{a}_n$. In particular, we can pick $a_n$ such that $\alpha(a_n),\alpha^{-1}(a_n)\notin \cone_{\gamma}(b c d a^* \alpha(a^*))$, concluding that $a^*a_n$ is $\alpha${\hyp}ordered and satisfies {\rm{(d)}}. 
\item[Case $\widetilde{a}_n\in\cone_\gamma(b c d a^* \alpha(a^*))$.] Set $a_n=\widetilde{a}_n$. By \cref{l:cones and cone independence}, $a_n\in\cone_{\gamma}(c a^*)$, getting {\rm{(d)}} vacuously. It remains to check point {\rm{(c)}}. Since $a^*a_n\equiv^{\Rr\gamma}_{c}a'_{<n}a'_n$, we conclude that $a'_n\in \cone_\gamma(c a'_{<n})$. 

If $a_n\sim_\gamma a_i$ with $i<n$, then $a'_n\sim_\gamma a'_i$. As $a'$ is $\alpha${\hyp}ordered, we get that actually $a'_n\sim_\gamma a'_{n-1}$, so $a_n\sim_\gamma a_{n-1}$. Thus, $\alpha(a_n)\sim_\gamma \alpha(a_{n-1})\not\sim_\gamma a_i$ for $i\leq n-1$ as $a^*$ is $\alpha${\hyp}ordered. Hence, $a^*a_n$ is $\alpha${\hyp}ordered. 

Now, consider the case $a_n\sim_\gamma c_0$ with $c_0$ in $c$ and, aiming a contradiction, suppose that $a^*a_n$ is not $\alpha${\hyp}ordered, i.e. $a_i\sim_\gamma \alpha(a_n)\sim_\gamma \alpha(c_0)$ for some $i<n$. Then, there is a smallest $i<n$ such that $a_i\sim_\gamma\alpha(c_0)$. As $i$ is the smallest, if $a_i\notin \cone_{\gamma}(c)$, then $a_i\notin\cone_{\gamma}(c a_{<i})$. In that case, since $a^*$ satisfies property {\rm{(d)}}, we have that $\alpha^{-1}(a_i)\notin \cone_\gamma(c)$, getting a contradiction. Thus, $a_i\in\cone_\gamma(c)$, so there is $c_1$ in $c$ such that $c_1\sim_\gamma a_i\sim_\gamma \alpha(c_0)$. Since $a^*\equiv^{\Rr\gamma}_c a'_{<n}$, we conclude that $\alpha(c_0)\sim_\gamma c_1\sim_\gamma a'_i$. On the other hand, $c_0\sim_\gamma a_n$, so $c_0\sim_\gamma a'_n$. Hence, $\alpha(a'_n)\sim_\gamma\alpha(c_0)\sim_\gamma c_1\sim_\gamma a'_i$, contradicting that $a'$ is $\alpha${\hyp}ordered.\smallskip
\end{enumerate}
Finally, it remains to check property {\rm{(e)}} for $a^*a_n$, i.e. $a^*a_n\indepe[\gamma]{d}\alpha(a^*)\alpha(a_n)$.

As $a_n\indepe[\Rr\gamma]{a^*c}bd\alpha(a^*)$, we get $a_n\indepe[\gamma]{a^*c}bd\alpha(a^*)$. By monotonicity, base monotonicity, normality and transitivity of $\indepe[\gamma]{}$, we get that $a^*a_n\indepe[\gamma]{d}\alpha(a^*)$.

Since $a^*a_n$ is $\alpha${\hyp}ordered and $\alpha$ is strictly increasing below $\gamma$, we get that $a_i\meet\gamma<\alpha(a_i\meet \gamma)\leq \alpha(a_n\meet\gamma)=\alpha(a_n)\meet\gamma$ when $a_i\meet \gamma<\gamma$. Also, if $a_i>\gamma$ for some $i$, then $a_n>\gamma$. In that case, as $a^*a_n$ is $\alpha${\hyp}ordered, $a_i\meet \alpha(a_n)\leq \gamma$ for each $i<n$. Hence, by \cref{l:easy independence}, we conclude that $a^*a_n\indepe[\gamma]{d}\alpha(a_n)$.

Now, we have that $\langle a^*a_n\rangle_{\gamma}=\langle a^*\rangle_{\gamma}\cup \{a_n\}$. Similarly, $\langle \alpha(a^*a_n)\rangle_\gamma=\langle \alpha(a^*)\rangle_\gamma\cup \{\alpha(a_n)\}$. Therefore, by \cref{l:generated}(\ref{itm:generated 4}), we get $\langle \alpha(a^*a_n)\, d\rangle_\gamma=\langle \alpha(a^*)\, d\rangle_\gamma\cup\langle \alpha(a_n)\, d\rangle_\gamma$. Hence, from $a^*a_n\indepe[\gamma]{d}\alpha(a^*)$ and $a^*a_n\indepe[\gamma]{d}\alpha(a_n)$, we conclude directly from \cref{d:cone-independence} that $a^*a_n\indepe[\gamma]{d}\alpha(a^*a_n)$.
\end{enumerate}
\item As in (1), but ``inverting'' the $\alpha${\hyp}ordering (i.e. taking the inverse of the $\alpha^{-1}${\hyp}ordering).\qedhere
\end{enumerate}
\end{proof}
\end{lem}

\begin{coro} \label{c:branched alpha ordering} Let $\Gamma$ be a branch, $\alpha$ be an automorphism of $\Tt^M$ setwise fixing $\Gamma$ and $a',b,c,d$ finite tuples with $c\subseteq d$.
\begin{enumerate}[label={\rm{(\arabic*)}}, wide]
\item If $\alpha$ is strictly increasing on $\Gamma$, then there is $a\equiv^{\Rr\Gamma}_ca'$ such that $a\indepe[\Rr\Gamma]{c}b$ and $a\indepe[\Gamma]{d}\alpha(a)$.
\item If $\alpha$ is strictly decreasing on $\Gamma$, then there is $a\equiv^{\Rr\Gamma}_ca'$ such that $b\indepe[\Rr\Gamma]{c}a$ and $\alpha(a)\indepe[\Gamma]{d}a$.
\end{enumerate}
\begin{proof} Write $T=\Tt^M$ and consider a strongly $\aleph_1${\hyp}homogeneous and $\aleph_1${\hyp}saturated elementary extension $\mathfrak{T}$ of $T$. 

By existence, extension and saturation, find $\gamma$ in $\mathfrak{T}$ such that $\gamma\indepe[\Rr]{}T\alpha(T)$. By monotonicity, $\gamma\indepe[\Rr]{}T$ and $\gamma\indepe[\Rr]{}\alpha(T)$, concluding by stationarity that there is $\gamma$ in $\mathfrak{T}$ such that $\gamma\, T\equiv^{\Rr}\gamma\, \alpha(T)$ --- where we enumerate $T$ and $\alpha(T)$ in the same way according to $\alpha$. Using the generic mix property, we can assume in addition that $\gamma>x$ for all $x\in\Gamma$. By \cref{c:quantifier free type tree}(\ref{itm:quantifier free type tree 3}), as $\qftp_{\Gamma}(T)=\qftp_{\Gamma}(\alpha(T))$, we get $\qftp_{\gamma}(T)=\qftp_\gamma(\alpha(T))$. Since $\langle T\, \gamma\rangle=T\cup \{\gamma\}$ and $\langle \alpha(T)\, \gamma\rangle=\alpha(T)\cup \{\gamma\}$  by \cref{l:generated}(\ref{itm:generated 4}), we conclude that $\qftp_{\Rr\gamma}(T)=\qftp_{\Rr\gamma}(\alpha(T))$. Hence, by quantifier elimination, we get that $T\equiv^{\Rr \gamma} \alpha(T)$. 

Using \cref{l:saturation of cones} recursively, we can find $\{y_i\}_{i\in\Z}$ above $\gamma$ from different cones such that $\widetilde{\beta}:\ T^+\rightarrow T^+$ extending $\alpha$ by $\gamma\mapsto \gamma$ and $y_i\mapsto y_{i+1}$ is a partial isomorphism, where $T^+\coloneqq T\cup\{\gamma\}\cup\{y_i\}_{i\in\Z}$. By strong $\aleph_1${\hyp}homogeneity of $\mathfrak{T}$, there is $\beta\in\Aut(\mathfrak{T})$ extending $\widetilde{\beta}$.

\begin{claim} \label{cl:omitting type} Let $A$ be a countable subset of $\mathfrak{T}$ containing $T^+$ and assume that $\{a\meet \gamma \tq a\in A\}=\Gamma$ and $\cone_\gamma(A)=\cone_\gamma(y_i\tq i\in \Z)$. Then, there is a countable elementary substructure $N\preceq \mathfrak{T}$ containing $A$ such that $N_{<\gamma}=\Gamma$ and $\cone_\gamma(N)=\cone_\gamma(y_i\tq i\in\Z)$.
\begin{proof} By the Omitting Type Theorem applied to the countable complete theory $\Th(\mathfrak{T}/A)$ and the non-isolated partial types $\pi_1(x)=(\bigwedge_{b\in\Gamma}x\neq b)\wedge x<\gamma$ and $\pi_2(x)=(\bigwedge_{i\in\Z} y_i\meet x=\gamma)\wedge x>\gamma$.
\claimqed 
\end{proof}
\end{claim}

Recursively applying \cref{cl:omitting type}, we construct a countable elementary chain $(N_i)_{i\in\N}$ of elementary substructures of $\mathfrak{T}$ with $T^+\subseteq N_0$, $\beta(N_i)\subseteq N_{i+1}$, $\beta^{-1}(N_i)\subseteq N_{i+1}$ and ${(N_i)}_{<\gamma}=\Gamma$ and $\cone_\gamma(N_i)=\cone_\gamma(y_i\tq i\in\Z)$. Take $N=\bigcup N_i$. Then, $T\preceq N\preceq\mathfrak{T}$ where $N$ is countable, $\beta(N)=N$, $N_{<\gamma}=\Gamma$ and $\cone_\gamma(N)=\cone_\gamma(y_i\tq i\in\Z)$. By $\omega${\hyp}categoricity, we get that $N\cong \Tt^M$.  Note that $\beta_{\mid N}$ is an infinite fan above $\gamma$ and has the same monotonicity in $N_{<\gamma}$ as $\alpha$. 
\begin{enumerate}[label={\rm{(\arabic*)}}, wide]
\item If $\alpha$ is strictly increasing in $\Gamma$, then $\beta_{\mid N}$ is strictly increasing in $N_{<\gamma}$. Applying \cref{l:alpha ordering}(\ref{itm:alpha ordering 1}) to the tuple $(N,T, \beta_{\mid N})$, we can find $a$ tuple in $T$ such that $a\indepe[\Rr\gamma]{c}b$, $a\indepe[\gamma]{d}\alpha(a)$ and $a\equiv^{\Rr\gamma}_ca'$. Now, by (the furthermore part of) \cref{l:semibranch and cone independence}, we get that $a\indepe[\Rr\Gamma]{c}b$ and $a\indepe[\Gamma]{d}\alpha(a)$. Also, applying \cref{c:quantifier free type tree}(\ref{itm:quantifier free type tree 3}), we conclude $a\equiv^{\Rr\Gamma}_ca'$. 
\item As before, but applying \cref{l:alpha ordering}(\ref{itm:alpha ordering 2}).\qedhere
\end{enumerate}
\end{proof}
\end{coro}

As indicated in the introduction, our strategy is similar to \cite{li2020automorphism}. The main two ingredients of that strategy are the notions from \cite{li2019automorphism} of stationary weak independence relations (see \cref{d:stationary weak independence relation}) and (almost) moving maximally automorphisms that we recall here:

\begin{defi} \label{d:move maximally} Let $\alpha:\ A\rightarrow B$ be an isomorphism between substructures and $p$ a type over a finite set $C$. We say\footnote{In \cite[Definition 1.3]{li2019automorphism}, Li uses ``almost'' moves (right/left) maximally to remark the difference with \cite[Definition 2.5]{tent2013isometry}.}  that $\alpha$ \emph{moves ($\indepe{}${\hyp})right{\hyp}maximally} $p$ if there is $a\vDash p$ on the domain of $\alpha$ such that $a\indepe{C}\alpha(a)$; we say that it moves ($\indepe{}${\hyp})right{\hyp}maximally if it moves ($\indepe{}${\hyp})right{\hyp}maximally any such type. Symmetrically, $\alpha$ \emph{moves ($\indepe{}${\hyp})left{\hyp}maximally} $p$ if there is $a\vDash p$ on the domain of $\alpha$ such that $\alpha(a)\indepe{C}a$; we say that it moves ($\indepe{}${\hyp})left{\hyp}maximally if it moves ($\indepe{}${\hyp})left{\hyp}maximally any such type.
\end{defi}

\begin{ejem} By existence, extension, invariance and monotonicity of $\indepe{}$, we have $a\indepe{C}b$ and $b\indepe{C}a$ for any $a,b,C$ with $a$ in $\acl(C)$. In particular, every automorphism moves $\indepe{}${\hyp}right{\hyp}maximally and $\indepe{}${\hyp}left{\hyp}maximally every algebraic type.
\end{ejem}

By \cref{l:alpha ordering,c:branched alpha ordering}, we have in particular the following corollary:
\begin{coro} \label{c:move maximally tree} Let $\Gamma$ be a branch and $\gamma$ a point. 
\begin{enumerate}[label={\rm{(\arabic*)}}, wide]
\item Every automorphism of $\Tt^M$ which is an infinite fan above $\gamma$ and strictly increasing below $\gamma$ moves $\indepe[\gamma]{}${\hyp}right{\hyp}maximally. 

Every automorphism of $\Tt^M$ which is an infinite fan above $\gamma$ and strictly decreasing below $\gamma$ moves $\indepe[\gamma]{}${\hyp}left{\hyp}maximally.
\item Every automorphism of $\Tt^M$ setwise fixing $\Gamma$ which is strictly increasing in $\Gamma$ moves $\indepe[\Gamma]{}${\hyp}right{\hyp}maximally. 

Every automorphism of $\Tt^M$ setwise fixing $\Gamma$ which is strictly decreasing in $\Gamma$ moves $\indepe[\Gamma]{}${\hyp}left{\hyp}maximally.
\end{enumerate}
\end{coro}

Up to this point we have only used that $\indepe[\Rr]{}$ is a stationary weak independence relation satisfying the independent generation property. In other words, we have not used the free property of $\indepe[\Rr]{}$. For instance, all the previous results still hold for the ordered universal dense meet{\hyp}tree $\Tt^{<}$.\smallskip

From now on, we assume that $\indepe[\Rr]{}$ is a free weak independence relation \cref{d:free weak independence}. Recall that freeness means:
\begin{enumerate}[wide]
\item[\rm{(Freeness)}] If $A\indepe{C}B$ and $\langle C\rangle \cap\langle AB\rangle\subseteq \langle D\rangle$ where $D\subseteq C$, then $A\indepe{D}B$.
\end{enumerate}

For instance, by \cite[Proposition 3.4]{conant2016axiomatic}, the free amalgamation independence relation in a free Fra\"{\i}ss\'{e} limit is a free independence relation. As a prototypical example, consider the Rado meet{\hyp}tree $\Tt^{\mathrm{R}}$.

To apply \cite[Corollary 2.16]{li2019automorphism}, we need to replace $\indepe[\gamma]{}$ and $\indepe[\Gamma]{}$ in \cref{c:move maximally tree} by $\indepe[\Rr\gamma]{}$ and $\indepe[\Rr\Gamma]{}$ respectively. We do this in the following lemma, which is essentially \cite[Theorem 3.2]{li2020automorphism} and relies on \cref{l:mix free independence}.

\begin{lem}\label{l:move maximally} Let $\Gamma$ be a branch and $\gamma$ a point.
\begin{enumerate}[label={\rm{(\arabic*)}}, wide]
\item Let $\alpha$ be an automorphism of $\Tt^M$ which is an infinite fan above $\gamma$ and strictly increasing below $\gamma$. Then, there is an automorphism $\beta$ of $\Tt^M$ fixing $\gamma$ such that $[\beta,\alpha]$ moves $\indepe[\Rr\gamma]{}${\hyp}right{\hyp}maximally and $[\alpha,\beta]$ moves $\indepe[\Rr\gamma]{}${\hyp}left{\hyp}maximally.
\item Let $\alpha$ be an automorphism of $\Tt^M$ setwise fixing $\Gamma$ and strictly increasing in $\Gamma$. Then, there is an automorphism $\beta$ of $\Tt^M$ setwise fixing $\Gamma$ such that $[\beta,\alpha]$ moves $\indepe[\Rr\Gamma]{}${\hyp}right{\hyp}maximally and $[\alpha,\beta]$ moves $\indepe[\Rr\Gamma]{}${\hyp}left{\hyp}maximally.
\end{enumerate}  
\begin{proof} The proofs for {\rm{(1)}} and {\rm{(2)}} are very similar, so we just expose the proof of {\rm{(1)}}. Let $(p_n)_{n\in\N}$ be an enumeration of all 
types in $\Tt^M_\gamma$ over finite sets. We define by recursion a chain $(\beta_n)_{n\in\N}$ of partial isomorphisms $\beta_n:\ X_n\rightarrow Y_n$ between finitely generated substructures of $\Tt^M_\gamma$ such that, for $i<n$, the parameters of $p_i$ are in $X_n\cap Y_n$, and $[\beta_n,\alpha]$ moves right{\hyp}maximally $p_i$ and $[\alpha,\beta_n]$ moves left{\hyp}maximally $p_i$.
\begin{enumerate}[topsep=3pt, wide]
\item[For $n=0$:] We set $\beta_0:\ \gamma\mapsto\gamma$ and $X_0=Y_0=\{\gamma\}$.
\item[Recursion case:] Suppose we already have $\beta_n$ for $n\geq 0$. Say $p_n$ is a type over $C$. By homogeneity, we extend $\beta_n$ by adding first $C\alpha(C)$ to the domain and, then, $C\alpha\beta_n(C)$ to the image. Let $\bar{x}$ enumerate $X_n$ and $\bar{y}$ be the enumeration of $Y_n$ given by $\beta_n:\ \bar{x}\mapsto \bar{y}$. 
\begin{enumerate}[label={\rm{Step {\arabic*}:}}, wide]
\item Applying \cref{l:alpha ordering}(\ref{itm:alpha ordering 1}) twice, there are $a\vDash p_n$ and $b$ with $a\bar{x}\equiv^{\Rr\gamma}b\bar{y}$ such that $a\indepe[\Rr\gamma]{C}\bar{x}$ and $a\indepe[\gamma]{\bar{x}}\alpha(a)$, and $b\indepe[\Rr\gamma]{\bar{y}}\alpha^{-1}(\bar{y})$ and $b\indepe[\gamma]{\bar{y}}\alpha(b)$. Since $a\bar{x}\equiv^{\Rr\gamma} b\bar{y}$, we get $b\indepe[\Rr\gamma]{\beta_n(C)}\bar{y}$ by invariance. By existence and extension of $\indepe[\Rr\gamma]{}$, we find $d$ with $da\bar{x}\equiv^{\Rr\gamma} \alpha(b)b\bar{y}$ such that $d\indepe[\Rr\gamma]{a\bar{x}}\alpha(a)$. By homogeneity, find $\beta'_{n+1}:\ X'_{n+1}\rightarrow Y'_{n+1}$ extending $\beta_n$ with $a\mapsto b$ and $d\mapsto \alpha(b)$, with $\alpha(a)\in X'_{n+1}$ and $\alpha^{-1}\beta'_{n+1}\alpha(a)\in Y'_{n+1}$. 

Since $b\indepe[\gamma]{\bar{y}}\alpha(b)$, we get $\langle b\bar{y}\rangle_\gamma\cap \langle \alpha(b)\bar{y}\rangle_\gamma=\langle \bar{y}\rangle_\gamma$ by \cref{o:reflexivity}. Since $b\alpha(b)\bar{y}\equiv^{\Rr\gamma}ad\bar{x}$, we conclude that $\langle a\bar{x}\rangle_{\gamma}\cap \langle d\bar{x}\rangle_{\gamma}=\langle \bar{x}\rangle_{\gamma}$. Now, $a\indepe[\gamma]{\bar{x}}\alpha(a)$ and $d\indepe[\Rr\gamma]{a\bar{x}}\alpha(a)$, so $d\indepe[\Rr\gamma]{\bar{x}}\alpha(a)$ by \cref{l:mix free independence}\footnote{Note that freeness is only used in this proof to apply \cref{l:mix free independence}}. 

Since $b\indepe[\Rr\gamma]{\beta_n(C)}\bar{y}$ and $b\indepe[\Rr\gamma]{\bar{y}}\alpha^{-1}(\bar{y})$, we get $b\indepe[\Rr\gamma]{\beta_n(C)}\alpha^{-1}(\bar{y})$ by transitivity. Applying ${\beta'}_{n+1}^{-1}\alpha$, we get \[d\indepe[\Rr\gamma]{{\beta'}^{-1}_{n+1}\alpha\beta'_{n+1}(C)}\bar{x}.\] As $d\indepe[\Rr\gamma]{\bar{x}}\alpha(a)$, we obtain $d\indepe[\Rr\gamma]{\beta'^{-1}_{n+1}\alpha\beta'_{n+1}(C)}\alpha(a)$ by transitivity. By applying $\beta'^{-1}_{n+1}\alpha^{-1}\beta'_{n+1}$, we conclude 
\[{\beta'}^{-1}_{n+1}\alpha^{-1}\beta'_{n+1}(d) \indepe[\Rr\gamma]{C} [\beta'_{n+1},\alpha](a).\]
Now, ${\beta'}^{-1}_{n+1}\alpha^{-1}\beta'_{n+1}(d)=a$, so $[\beta'_{n+1},\alpha]$ moves right{\hyp}maximally $p_n$.

\item Let $\bar{x}'$ enumerate $X'_{n+1}$ and $\bar{y}'$ be the enumeration of $Y'_{n+1}$ given by $\beta'_{n+1}:\ \bar{x}'\mapsto\bar{y}'$. 

Since $\alpha(C)\subseteq X'_{n+1}$, we get $C\subseteq  \alpha^{-1}(X'_{n+1})$. By \cref{l:alpha ordering}(\ref{itm:alpha ordering 2}) applied to $\alpha^{-1}$, there is $a'\vDash p_n$ such that $\alpha^{-1}(\bar{x}')\indepe[\Rr\gamma]{C}a'$ and $\alpha^{-1}(a')\indepe[\gamma]{\alpha^{-1}(\bar{x}')}a'$. Thus, by invariance, $\bar{x}'\indepe[\Rr\gamma]{\alpha(C)}\alpha(a')$ and $a'\indepe[\gamma]{\bar{x}'}\alpha(a')$. On the other hand, by \cref{c:move maximally tree}, there is $b'\bar{y}'\equiv^{\Rr\gamma} a'\bar{x}'$ such that $b'\indepe[\gamma]{\bar{y}'}\alpha(b')$. By existence and extension of $\indepe[\Rr\gamma]{}$, we find $d'$ with $d'a'\bar{x}'\equiv^{\Rr\gamma} \alpha(b')b'\bar{y}'$ such that $d'\indepe[\Rr\gamma]{a'\bar{x}'}\alpha(a')$. By homogeneity, find $\beta_{n+1}:\ X_{n+1}\rightarrow Y_{n+1}$ extending $\beta'_{n+1}$ with $a'\mapsto b'$ and $d'\mapsto \alpha(b')$, with $\alpha(a')\in X_{n+1}$ and $\alpha^{-1}\beta_{n+1}\alpha(a')\in Y_{n+1}$. 

Since $b'\indepe[\gamma]{\bar{y}'}\alpha(b')$, we have $\langle b'\bar{y}'\rangle_\gamma\cap\langle \alpha(b')\bar{y}'\rangle_\gamma=\langle\bar{y}'\rangle_\gamma$ by \cref{o:reflexivity}. As $b'\alpha(b')\bar{y}'\equiv^{\Rr\gamma}a'd'\bar{x}'$, we get $\langle a'\bar{x}'\rangle_\gamma \cap \langle d'\bar{x}'\rangle_\gamma=\langle\bar{x}'\rangle_\gamma$. Now, $d'\indepe[\Rr\gamma]{a'\bar{x}'}\alpha(a')$ and $a'\indepe[\gamma]{\bar{x}'}\alpha(a')$, so $d'\indepe[\Rr\gamma]{\bar{x}'}\alpha(a')$ by \cref{l:mix free independence}. 

As $\bar{x}'\indepe[\Rr\gamma]{\alpha(C)}\alpha(a')$ and $d'\indepe[\Rr\gamma]{\bar{x}'}\alpha(a')$, we conclude $d'\indepe[\Rr\gamma]{\alpha(C)}\alpha(a')$ by transitivity. Applying $\alpha^{-1}$, we get $\alpha^{-1}(d')\indepe[\Rr\gamma]{C}a'$. Now, $\alpha^{-1}(d')=\alpha^{-1}\beta^{-1}_{n+1}\alpha(b')=[\alpha,\beta_{n+1}](a')$, concluding that $[\alpha,\beta_{n+1}]$ moves left{\hyp}maximally $p_n$.
\end{enumerate}
\end{enumerate}

Taking the union $\beta=\bigcup\beta_n$, we get an automorphism of $\Tt^M_\gamma$ such that $[\beta,\alpha]$ moves right{\hyp}maximally and $[\alpha,\beta]$ moves left{\hyp}maximally.
\end{proof}
\end{lem}
\begin{defi} For a group $G$, we say that an element $a\in G$ is \emph{conormal} if it does not belong to a proper normal subgroup of $G$, i.e. $\llangle a\rrangle_G = G$ where $\llangle a \rrangle_G$ denotes the smallest normal subgroup of $G$ containing $a$.
\end{defi} 




Recall \cite[Corollary 2.16]{li2019automorphism}:

\begin{fact} \label{f:li2019automorphism} Let $N$ be a homogeneous countable structure with a stationary weak independence relation and $\alpha$ an automorphism of $N$ such that $\alpha$ moves right{\hyp}maximally and $\alpha^{-1}$ moves left{\hyp}maximally. Then, any automorphism of $N$ is the product of $16$ conjugates of $\alpha$ and $\alpha^{-1}$. In particular, $\alpha$ is conormal in $\Aut(N)$.
\end{fact}

\begin{coro}\label{c:main branched} Let $\Gamma$ be a branch. Every automorphism of $\Tt^{M}_\Gamma$ unbounded in $\Gamma$ is conormal in $\Aut(\Tt^{M}_\Gamma)$.
\begin{proof} By \cref{c:semibranch strictly increasing,l:move maximally,f:li2019automorphism}.
\end{proof}
\end{coro}
\begin{coro} \label{c:main pointed} Let $\gamma$ be a point. Every automorphism of $\Tt^M_\gamma$ permuting infinitely many cones above $\gamma$ and unbounded below $\gamma$ is conormal in $\Aut(\Tt^{M}_\gamma)$.
\begin{proof} By \cref{c:fan strictly increasing,l:move maximally,f:li2019automorphism}.
\end{proof}
\end{coro}

\begin{teo}\label{t:main} Every automorphism of $\Tt^M$ that does not pointwise fix an interval is conormal.
\begin{proof} Suppose that $\alpha$ is an autnormophism that does not fix pointwise any interval and that $\beta\in \Aut(\Tt^M)$. By \cref{l:dicotomia}, there are four cases: 
\begin{enumerate}[label={\rm Case \arabic*:}, wide]
\item $\alpha,\beta$ fix a branch setwise. By \cref{l:generic meet-tree expansions}(\ref{itm:generic meet-tree expansions 3}), all branches are conjugate, so we may assume that both setwise fix the same branch $\Gamma$. Since $\alpha$ does not pointwise fix an interval, it is unbounded in $\Gamma$, and hence conormal in $\Aut(\Tt^{M}_\Gamma)$ by \cref{c:main branched}. In particular, $\beta\in\llangle \alpha\rrangle$. 
\item $\alpha$ fixes a branch $\Gamma$ setwise while $\beta$ is a fan above $\gamma$. By \cref{l:trivial obs,c:main branched} taking products of conjugates of $\alpha,\alpha^{-1}$, we can replace $\alpha$ by an automorphism setwise fixing $\Gamma$, permuting infinitely many cones above $\gamma$ and fixing no point but $\gamma$ and $\pi(\gamma)$. In particular, $\alpha$ is unbounded below $\gamma$. By \cref{c:main pointed}, we get that $\alpha$ is conormal in $\Aut(\Tt^M_\gamma)$, so, in particular, $\beta\in\llangle \alpha\rrangle$. 
\item $\alpha$ is a fan while $\beta$ fixes a branch setwise. By \cref{l:trivial obs,c:main pointed}, taking products of conjugates of $\alpha,\alpha^{-1}$, we can replace $\alpha$ by an automorphism setwise fixing some branch and pointwise fixing no interval. Hence, by point {\rm{(1)}}, we conclude $\beta\in\llangle\alpha\rrangle$.
\item $\alpha,\beta$ are fans above some points. Similar as before, by \cref{l:trivial obs,c:main pointed}, we can replace $\alpha$ by an automorphism setwise fixing some branch and pointwise fixing no interval. Hence, by point {\rm{(2)}}, we conclude $\beta\in\llangle\alpha\rrangle$. \qedhere
\end{enumerate}
\end{proof}
\end{teo}

\begin{defi}
We say that a structure $M$ is \emph{$n${\hyp}ary}, for $n\in \N$, if, for any injective (i.e. without repetitions) tuples $a=(a_i)_{i \in I}$ and $b=(b_i)_{i\in I}$, we have that $a\equiv b$ whenever $\tp(a_i \tq i\in I')=\tp(b_i \tq i\in I')$ for every $I'\subseteq I$ with $|I'|=n$.
\end{defi}



\begin{lem} \label{l:no fixing intervals} Assume $M$ is not unary. Then, the only automorphism of $\Tt^M$ pointwise fixing an interval is the identity.
\begin{proof} We prove the contrapositive. Suppose that there is a non{\hyp}trivial automorphism $\alpha$ of $\Tt^M$ pointwise fixing an interval $I$. 

Pick any $1${\hyp}type $p$ in $\lang_\Rr$ over the empty set. By the generic mix property and \cref{l:generic mix over parameters} (and quantifier elimination), the type $p$ is realised in every interval of the tree. Since $\alpha$ is non{\hyp}trivial, there must be a realisation $a$ of $p$ which is not fixed. Indeed, suppose every realisation of $p$ is fixed and take $b$ such that $b\neq \alpha(b)$. Let $c\coloneqq b\meet\alpha(b)$. If $c\neq b$, a realisation $a$ of $p$ between $c$ and $b$ satisfies $c<a<b$ and $c<a=\alpha(a)<\alpha(b)$. Similarly, if $c\neq \alpha(b)$, a realisation $a$ of $p$ between $c$ and $\alpha(b)$, satisfies $c<a=\alpha^{-1}(a)<b$ and $c<a<\alpha(b)$. In both cases, a contradiction with the definition of $c$.

Again, by the generic mix property and \cref{l:generic mix over parameters} 
(and quantifier elimination), for any tuple $b$ disjoint to $a$ and $\alpha(a)$, there is some $b'$ in $I$ with $b'\equiv^\Rr_{a\,\alpha(a)}b$. Hence, since $\alpha$ pointwise fixes $I$, we conclude $a\equiv^\Rr_b\alpha(a)$. As $b$ is arbitrary, we conclude that $a\,\alpha(a)$ realise the partial $\lang_\Rr${\hyp}type $\tau(x,y)$ saying $\bigwedge_{\varphi\ \mathrm{formula}}\forall z (z\mathrm{\ disjoint\ to\ }x\mathrm{\ and\ }y\rightarrow (\varphi(x,z)\leftrightarrow \varphi(y,z)))$.

Pick any $a'$ realising $p$ with $a'\neq a$. By existence and extension of $\indepe[\Rr]{}$ and the generic mix property and \cref{l:generic mix over parameters} (and quantifier elimination), find $a''\in I$ with $a''\indepe[\Rr]{a}\alpha(a)$ and $a''\equiv^\Rr_aa'$. Thus, by freeness, $a''\indepe[\Rr]{}\alpha(a)$ and, applying $\alpha^{-1}$, $a''\indepe[\Rr]{}a$, so $a'\indepe[\Rr]{}a$ by invariance. In other words, any realisation of $p$ other than $a$ is $\indepe[\Rr]{}${\hyp}independent of $a$ over the empty set. In particular, by stationarity of $\indepe[\Rr]{}$, any two different realisations of $p$ have the same $\lang_\Rr${\hyp}type. Hence, any two realisations of $p$ realise $\tau$. As $p$ is arbitrary, it follows that any two elements in $M$ having the same type over the empty set have the same type over any other tuple disjoint to them. In other words, $M$ is unary: suppose $a=(a_i)_{i<n}$ and $b=(b_i)_{i<n}$ are injective tuples such that $a_i\equiv^{\Rr}b_i$ for each $i<n$ but $a \not \equiv^{\Rr} b$ with $n$ minimal. Then, $a_{<n-1}\equiv^{\Rr}b_{<n-1}$. Applying an automorphism mapping $a_{<n-1}$ to $b_{<n-1}$, $a_{n-1}$ goes to an element with the same type as $b_{n-1}$ over the empty set but different type over $b_{<n-1}$.  
\end{proof}
\end{lem}

\begin{coro} \label{c:main} The group of automorphisms of the generic (branched) meet{\hyp}tree expansion of a non{\hyp}unary infinite strong Fra\"{\i}ss\'{e} limit over a finite relational language having a free weak independence relation is simple.
\end{coro}
\begin{coro}\label{c:main rado} The group of automorphisms of the Rado (branched) meet{\hyp}tree is simple.
\end{coro}

\begin{conj} Let $\Tt^{\Rr}$ be the Rado meet{\hyp}tree and $A$ a finite subset. For $\alpha\in\Aut(\Tt^{\Rr}/A)$ and $a\in A$, let $\alpha_a$ be the induced permutation on the cones above $a$ which are disjoint to $\cone_a(A)$. Consider the group homomorphism $f:\ \Aut(\Tt^{\Rr}/A)\rightarrow \bigoplus_{a\in A} \mathbf{S}$ give by $\alpha\mapsto (\alpha_a)_{a\in A}$, where $\mathbf{S}$ is the group of infinite countable permutations. Then, the non{\hyp}trivial normal subgroups of $\Aut(\Tt^{\Rr}/A)$ are precisely the preimages by $f$ of the normal subgroups of $\bigoplus_{a\in A}\mathbf{S}$.\end{conj}

\begin{quest} Can we get rid of freeness? For instance, is the group of automorphisms of the generically ordered universal dense meet{\hyp}tree $\Tt^<$ simple?
\end{quest}


\newpage
\bibliographystyle{alpha}
\bibliography{RadoMeetTree}
\end{document}